# Intertwining operators of the quantum Teichmüller space

Filippo Mazzoli

October 11, 2016

## Abstract

In [BBL07] the authors introduced the theory of *local* representations of the quantum Teichmüller space $\mathcal{T}_S^q$ ($q$ being a fixed primitive $N$-th root of $(-1)^{N+1}$) and they studied the behaviour of the intertwining operators in this theory. One of the main results [BBL07, Theorem 20] was the possibility to select one distinguished operator (up to scalar multiplication) for every choice of a surface $S$, ideal triangulations $\lambda, \lambda'$ and isomorphic local representations $\rho, \rho'$, requiring that the whole family of operators verifies certain *Fusion* and *Composition* properties. This selection was also used to produce invariants for pseudo-Anosov diffeomorphisms and their hyperbolic mapping tori (extending to local representations what had been done in [BL07] for irreducible ones). However, by analyzing the constructions of [BBL07], we found a difficulty that we eventually fix by a slightly weaker (but actually optimal) selection procedure. In fact, for every choice of a surface $S$, ideal triangulations $\lambda, \lambda'$ and isomorphic local representations $\rho, \rho'$, we select a *finite* set of intertwining operators, naturally endowed with a structure of affine space over $H_1(S; \mathbb{Z}_N)$ ($\mathbb{Z}_N$ is the cyclic group of order $N$), in such a way that the whole family of operators verifies *augmented* Fusion and Composition properties, which incorporate the explicit behavior of the $\mathbb{Z}_N$-actions with respect to such properties. Moreover, this family is *minimal* among the collections of operators verifying the "weak" Fusion and Composition rules (in practice the ones considered in [BBL07]). In addition, we adapt the derivation of the invariants for pseudo-Anosov diffeomorphisms and their hyperbolic mapping tori made in [BBL07] and [BL07] by using our distinguished family of intertwining operators.

## Contents

# Introduction

The quantum Teichmüller space is an algebraic object associated with a punctured surface admitting an ideal triangulation. Two somewhat different versions of it have been introduced, as a quantization by deformation of the Teichmüller space of a surface, independently by Chekhov and Fock [CF99] and by Kashaev [Kas95]. As in the articles [BL07] and [BBL07], we follow the exponential version of the Chekhov-Fock approach, whose setting has been established in [Liu09]. In this way the study is focused on non-commutative algebras and their finite-dimensional representations, instead of Lie algebras and self-adjoint operators on Hilbert spaces, as in [CF99] and [Kas95].

Given $S$ a surface admitting an ideal triangulation $\lambda$, we can produce a non-commutative $\mathbb{C}$-algebra $\mathcal{T}_\lambda^q$ generated by variables $X_i^{\pm 1}$ corresponding to the edges of $\lambda$ and endowed with relations $X_i X_j = q^{2\sigma_{ij}} X_j X_i$, where $\sigma_{ij}$ is an integer number, depending on the mutual position of the edges $\lambda_i$ and $\lambda_j$ in $\lambda$, and $q \in \mathbb{C}^*$ is a complex number. The algebra $\mathcal{T}_\lambda^q$ is called the *Chekhov-Fock algebra* associated with the surface $S$ and the ideal triangulation $\lambda$. Varying $\lambda$ in the set $\Lambda(S)$ of all the ideal triangulations of $S$, we obtain a collection of algebras, whose fraction rings $\widehat{\mathcal{T}}_\lambda^q$ are related by isomorphisms $\Phi^q_{\lambda\lambda'} : \widehat{\mathcal{T}}_{\lambda'}^q \to \widehat{\mathcal{T}}_\lambda^q$. This structure allows us to consider an object realized by "gluing" all the $\widehat{\mathcal{T}}_\lambda^q$ through the maps $\Phi^q_{\lambda\lambda'}$. The result of this procedure is an intrinsic algebraic object, called the *quantum Teichmüller space* of $S$ and denoted by $\mathcal{T}^q_S$, which does not depend on the chosen ideal triangulation, but just on $S$ and on the $N$-th root $q$.

In [BL07] the authors have given a suitable notion of finite-dimensional representations of the quantum Teichmüller space and have studied its irreducible representations. A representation of the quantum Teichmüller space is a collection
$$\rho = \{\rho_\lambda : \mathcal{T}^q_\lambda \longrightarrow \mathrm{End}(V_\lambda)\}_{\lambda \in \Lambda(S)}$$
of representations of all the Chekhov-Fock algebras associated with the surface $S$, verifying a compatibility condition in terms of the isomorphisms $\Phi^q_{\lambda\lambda'}$. More precisely, two representations $\rho_\lambda \colon \mathcal{T}^q_\lambda \to \mathrm{End}(V_\lambda)$ and $\rho_{\lambda'} \colon \mathcal{T}^q_{\lambda'} \to \mathrm{End}(V_{\lambda'})$ are compatible if $\rho_\lambda \circ \Phi^q_{\lambda\lambda'}$ makes sense and it is isomorphic to $\rho_{\lambda'}$. A necessary condition for the existence of a finite-dimensional representation of any Chekhov-Fock algebra is that $q^2$ is a root of unity, hence we always assume that $q$ is a primitive $N$-th root of $(-1)^{N+1}$.

Later, Bai, Bonahon, and Liu [BBL07] introduced a new type of representations of the quantum Teichmüller space of $S$, called *local representations*, which are constructed as "fusion" of irreducible representations of the Chekhov-Fock algebras associated with the triangles composing an ideal triangulation $\lambda$. These representations have a simpler set of invariants than the irreducible ones and they can be glued in order to construct local representations of surfaces obtained by identifying couples of edges in $(S, \lambda)$. Even in this case, the authors of [BBL07] have found a good notion of local representations of the quantum Teichmüller space, defined as collections of compatible local representations of the Chekhov-Fock algebras. In addition, they have developed a theory of *intertwining operators* between couples of isomorphic local representations of the quantum Teichmüller space.

More precisely, let
$$\rho = \{\rho_\lambda \colon \mathcal{T}^q_\lambda \to \mathrm{End}(V_\lambda)\}_{\lambda \in \Lambda(S)} \qquad \rho' = \{\rho'_\lambda \colon \mathcal{T}^q_\lambda \to \mathrm{End}(V'_\lambda)\}_{\lambda \in \Lambda(S)}$$
be two isomorphic local representations of the quantum Teichmüller space $\mathcal{T}^q_S$. By definition, for every $\lambda, \lambda' \in \Lambda(S)$, the representations $\rho_\lambda \circ \Phi^q_{\lambda\lambda'}$ and $\rho_{\lambda'}$ are isomorphic and $\rho_{\lambda'}$ itself is isomorphic to $\rho'_{\lambda'}$. Therefore, there exists a linear isomorphism $L^{\rho\rho'}_{\lambda\lambda'} \colon V'_{\lambda'} \to V_\lambda$ such that
$$(\rho_\lambda \circ \Phi^q_{\lambda\lambda'})(X') = L^{\rho\rho'}_{\lambda\lambda'} \circ \rho'_{\lambda'}(X') \circ (L^{\rho\rho'}_{\lambda\lambda'})^{-1} \qquad\qquad \forall X' \in \mathcal{T}^q_{\lambda'}$$

Such a $L^{\rho\rho'}_{\lambda\lambda'}$ is called an *intertwining operator*. In general, fixed $\rho, \rho'$ and $\lambda, \lambda'$, there is not a unique intertwining operator, not even up to scalar multiplication (in the following, we denote by $\doteq$ the equality between two linear isomorphisms up to non-zero scalar multiplication). Indeed, in light of the irreducible decomposition of local representations, described in [Tou14] for surfaces with genus $g \geq 1$ and $p+1$ punctures, they admit a lot of non-trivial automorphisms (see for example the case $\rho = \rho'$ and $\lambda = \lambda'$). One of the main purposes of [BBL07] was to select a unique intertwining operator $\widehat{L}^{\rho\rho'}_{\lambda\lambda'}$ for every $\rho, \rho', \lambda, \lambda$, by requesting some additional properties on them. More precisely, one of the results stated in [BBL07] was the following theorem:

**Theorem** ([BBL07, Theorem 20]). *For every surface $S$ there exists a unique family of intertwining operators $\widehat{L}_{\lambda\lambda'}^{\rho\rho'}$, indexed by couples of isomorphic local representations of $\mathcal{T}_S^q$ and by couples of ideal triangulations $\lambda, \lambda' \in \Lambda(S)$, individually defined up to scalar multiplication, such that:*

COMPOSITION RELATION: *for every $\lambda, \lambda', \lambda'' \in \Lambda(S)$ and for every triple of isomorphic local representations $\rho, \rho', \rho''$, we have $\widehat{L}_{\lambda\lambda''}^{\rho\rho''} \doteq \widehat{L}_{\lambda\lambda'}^{\rho\rho'} \circ \widehat{L}_{\lambda'\lambda''}^{\rho'\rho''}$;*

FUSION RELATION: *let $S$ be a surface obtained from another surface $R$ by fusion, and let $\lambda, \lambda'$ be two triangulations of $S$ obtained by fusion of two triangulations $\mu, \mu'$ of $R$. If $\eta, \eta'$ are two isomorphic local representations of $\mathcal{T}_R^q$ and $\rho, \rho'$ are local representations of $\mathcal{T}_S^q$ obtained by fusion respectively from $\eta$ and $\eta'$, then we have $\widehat{L}_{\eta\eta'}^{\mu\mu'} \doteq \widehat{L}_{\lambda\lambda'}^{\rho\rho'}$.*

However, in our study of the ideas exposed in [BBL07], we have found a difficulty that compromises this statement, in particular the possibility to select a unique intertwining operator for every choice of $\rho, \rho', \lambda, \lambda'$.

Let us try to describe this obstruction. Let $\lambda$ be an ideal triangulation of $S$. Denote by $S_0$ the surface obtained by splitting $S$ along all the edges of the ideal triangulation $\lambda$ and by $\lambda_0$ its unique ideal triangulation. In [BBL07] the procedure to select the isomorphism $L_{\lambda\lambda}^{\rho\rho'}$ was the following: we fix two representatives $(\rho_j)_{j=1}^m$ and $(\rho'_j)_{j=1}^m$ respectively of $\rho_\lambda$ and $\rho'_\lambda$ that are isomorphic to each other, as local representations of the Chekhov-Fock algebra $\mathcal{T}_{\lambda_0}^q$. As observed in the proof of [BBL07, Lemma 21], a local representation of $\mathcal{T}_{\lambda_0}^q$ is irreducible, so there exists a unique isomorphism $M_{\lambda\lambda}^{\rho\rho'}$ between $(\rho_j)_{j=1}^m$ and $(\rho'_j)_{j=1}^m$, up to scalar multiplication. Then the isomorphism $\widehat{L}_{\lambda\lambda}^{\rho\rho'}$ was defined in [BBL07, Lemma 22] as $\widehat{L}_{\lambda\lambda}^{\rho\rho'} := M_{\lambda\lambda}^{\rho\rho'}$. The problem we will observe is that this choice does depend on the selected representatives $(\rho_j)_j$ and $(\rho'_j)_j$. In other words, the representation $\rho_\lambda$ has representatives that are non-trivially isomorphic to each other, so different choices of $(\rho_j)_{j=1}^m$ and $(\rho'_j)_{j=1}^m$ lead us to a (finite) collection of intertwining operators, in general not to a unique element. We will focus on this problem in the first Section and in particular in Subsubsection 1.4.2.

The main purpose of this paper is to understand this phenomenon and to try to recover a result on intertwining operators similar to the one in [BBL07, Theorem 20]. We will be able to select just a set $\mathscr{L}_{\lambda\lambda'}^{\rho\rho'}$ of intertwining operators, instead of a unique linear isomorphism. Moreover, the selected sets $\mathscr{L}_{\lambda\lambda'}^{\rho\rho'}$ will be naturally endowed with a transitive and free action of $H_1(S; \mathbb{Z}_N)$. This fact implies that, fixed $\lambda, \lambda'$, the set $\mathscr{L}_{\lambda\lambda'}^{\rho\rho'}$ is finite, but its cardinality goes to infinity by increasing the number $N \in \mathbb{N}$ and the complexity of the surface $S$. In order to provide a statement analogous to [BBL07, Theorem 20], a few steps will be necessary: we will produce the fundamental objects in Section 2, we will investigate on their properties in Section 3 and finally we will complete the procedure in Section 4, where we will prove the following Theorem:

**Theorem.** *For every surface $S$ there exists a collection $\{(\mathscr{L}_{\lambda\lambda'}^{\rho\rho'}, \psi_{\lambda\lambda'}^{\rho\rho'})\}$, indexed by couples of isomorphic local representations $\rho, \rho'$ of the quantum Teichmüller space $\mathcal{T}_S^q$ and by couples of ideal triangulations $\lambda, \lambda' \in \Lambda(S)$ such that*

INTERTWINING: for every couple of isomorphic local representations
$$\rho = \{\rho_\lambda \colon \mathcal{T}_\lambda^q \to \mathrm{End}(V_\lambda)\}_{\lambda \in \Lambda(S)} \quad \rho' = \{\rho'_\lambda \colon \mathcal{T}_\lambda^q \to \mathrm{End}(V'_\lambda)\}_{\lambda \in \Lambda(S)}$$
and for every $\lambda, \lambda' \in \Lambda(S)$, $\mathscr{L}_{\lambda\lambda'}^{\rho\rho'}$ is a set of linear isomorphisms $L_{\lambda\lambda'}^{\rho\rho'}$ from $V'_{\lambda'}$ to $V_\lambda$ such that
$$(\rho_\lambda \circ \Phi_{\lambda\lambda'}^q)(X') = L_{\lambda\lambda'}^{\rho\rho'} \circ \rho'_{\lambda'}(X') \circ (L_{\lambda\lambda'}^{\rho\rho'})^{-1} \qquad \forall X' \in \mathcal{T}_{\lambda'}^q$$
for every $X' \in \mathcal{T}_{\lambda'}^q$;

ACTION: every $\mathscr{L}_{\lambda\lambda'}^{\rho\rho'}$ is endowed with a transitive and free action $\psi_{\lambda\lambda'}^{\rho\rho'}$ of $H_1(S; \mathbb{Z}_N)$;

FUSION PROPERTY: let $R$ be a surface and $S$ obtained by fusion from $R$. Fix
$$\eta = \{\eta_\mu \colon \mathcal{T}_\mu^q \to \mathrm{End}(W_\mu)\}_{\mu \in \Lambda(R)} \quad \eta' = \{\eta'_\mu \colon \mathcal{T}_\mu^q \to \mathrm{End}(W'_\mu)\}_{\mu \in \Lambda(R)}$$
two isomorphic local representations of $\mathcal{T}_R^q$ and
$$\rho = \{\rho_\lambda \colon \mathcal{T}_\lambda^q \to \mathrm{End}(V_\lambda)\}_{\lambda \in \Lambda(S)} \quad \rho' = \{\rho'_\lambda \colon \mathcal{T}_\lambda^q \to \mathrm{End}(V'_\lambda)\}_{\lambda \in \Lambda(S)}$$
two isomorphic local representations of $\mathcal{T}_S^q$, with $\rho$ and $\rho'$ respectively obtained by fusion from $\eta$ and $\eta'$. Then for every $\mu, \mu' \in \Lambda(R)$, if $\lambda, \lambda' \in \Lambda(S)$ are the corresponding ideal triangulations on $S$, there exists a natural inclusion $j \colon \mathscr{L}_{\mu\mu'}^{\eta\eta'} \to \mathscr{L}_{\lambda\lambda'}^{\rho\rho'}$ such that for every $L \in \mathscr{L}_{\mu\mu'}^{\eta\eta'}$ the following holds
$$(j \circ \psi_{\mu\mu'}^{\eta\eta'})(c, L) = \psi_{\lambda\lambda'}^{\rho\rho'}(\pi_*(c), j(L))$$
for every $c \in H_1(R; \mathbb{Z}_N)$, where $\pi \colon R \to S$ is the projection map;

COMPOSITION PROPERTY: for every $\rho, \rho', \rho''$ isomorphic local representations of $\mathcal{T}_S^q$ and for every $\lambda, \lambda', \lambda'' \in \Lambda(S)$, the composition map
$$\begin{array}{rcl} \mathscr{L}_{\lambda\lambda'}^{\rho\rho'} \times \mathscr{L}_{\lambda'\lambda''}^{\rho'\rho''} & \longrightarrow & \mathscr{L}_{\lambda\lambda''}^{\rho\rho''} \\ (L, M) & \longmapsto & L \circ M \end{array}$$
is well defined and it verifies
$$(c \cdot L) \circ (d \cdot M) = (c + d) \cdot (L \circ M)$$

Observe that in the reviewed assertion the Fision and Composition properties have been modified incorporating the affine structure on $H_1(S; \mathbb{Z}_N)$.

Moreover, we will prove that the collection $\{\mathscr{L}_{\lambda\lambda'}^{\rho\rho'}\}$ of intertwining operators we will exhibit is minimal in the family of collections of intertwining operators verifying "weak" Fusion and Composition properties (in practice the ones of the original statement in [BBL07]).

It turns out that the basic building block of our family $\{\mathscr{L}_{\lambda\lambda'}^{\rho\rho'}\}$ is the unique intertwining operator of local representations carried by a square equipped with two triangulations related by a diagonal exchange. In Subsection 2.3 we provide explicit formulas for such basic operators.

An important application of [BBL07, Theorem 20] was the construction of invariants of pseudo-Anosov surface diffeomorphisms and their hyperbolic mapping tori, extending to local representations what had been done in [BL07] for irreducible ones (for which the uniqueness of the intertwining operator up to scalar multiples is automatic). There are good reasons to conjecture that such "local" invariants eventually are instances of *quantum hyperbolic* invariants of cusped hyperbolic 3-manifolds defined in [BB05], [BB07] and [BB15] by developing seminal Kashaev's work [Kas94]. This is indeed a reason of the interest of local representations. We will show that this construction can be adapted rather straightforwardly by replacing [BBL07, Theorem 20] with our main theorem, although the resulting invariant is a bit more complicated. This could also indicate that establishing an actual relation between "quantum Teichmüller" and quantum hyperbolic invariants for fibred cusped hyperbolic 3-manifolds can be a bit subtler than one expects.

**Aknowledgement** This article is extracted from my Master's thesis at the University of Pisa. I would like to thank my advisor Prof. R. Benedetti for his help and support all along this work.

# 1 Preliminaries

In this paper we will always assume that a surface $S$ is oriented and obtained, from a compact oriented surface $\overline{S}$ with genus $g$ and $b$ boundary components, by removing $p \geq 1$ punctures $v_1, \ldots, v_p$, with at least a puncture in each boundary components. Assume further that $2\chi(S) \leq p_\partial - 1$, where $p_\partial$ is the number of punctures lying in $\partial S$. These are the hypotheses in which $S$ admits an ideal triangulation (see below for the Definition).

Moreover, $q \in \mathbb{C}^*$ will always be a certain primitive $N$-th root of $(-1)^{N+1}$.

## 1.1 The Chekhov-Fock algebra

**Definition 1.1.** Let $S$ be a surface. We define an *ideal triangulation* $\lambda$ of $S$ as a triangulation of $\overline{S}$ having as set of vertices exactly the set of punctures $\{v_1, \ldots, v_p\}$. We identify two ideal triangulations of $S$ if they are isotopic. Also, we denote by $\Lambda(S)$ the set of all ideal triangulations of $S$.

Given $\lambda \in \Lambda(S)$, let $n$ be the number of 1-cells in the ideal triangulation $\lambda$ and $m$ the number of faces of $\lambda$. Easy calculations show that the following relations hold:

$$\begin{aligned} n &= -3\chi(\overline{S}) + 3p - p_\partial \\ &= -3\chi(S) + 2p_\partial \\ m &= -2\chi(S) + p_\partial \end{aligned}$$

In particular $m$ and $n$ do not depend on the ideal triangulation $\lambda \in \Lambda(S)$. Since $S$ is oriented, on each triangle of $\lambda$ we have a natural induced orientation. With respect to this orientation, it is reasonable to speak about a left and a right side of each spike of a triangle. We select an order $\lambda_1, \ldots, \lambda_n$ on the set of 1-cells of the triangulation $\lambda$. Given $\lambda_i$ and $\lambda_j$ 1-cells of $\lambda$, we denote by $a_{ij}$ the number of spikes of triangles in $\lambda$ in which we find $\lambda_i$ on the left side and $\lambda_j$ on the right. Now we set

$$\sigma_{ij} := a_{ij} - a_{ji}$$

**Definition 1.2.** Given $q \in \mathbb{C}^*$, we define the *Chekhov-Fock algebra* associated with the ideal triangulation $\lambda$ and the parameter $q$ as the non-commutative $\mathbb{C}$-algebra $\mathcal{T}_\lambda^q$, generated by $X_1^{\pm 1}, \ldots, X_n^{\pm 1}$ and endowed with the following relations
$$X_i X_j = q^{2\sigma_{ij}} X_j X_i$$
for every $i, j = 1, \ldots, n$.

Given $\lambda \in \Lambda(S)$, we designate the free $\mathbb{Z}$-module generated by the 1-cells of the triangulation $\lambda$ as $\mathcal{H}(\lambda; \mathbb{Z})$. A choice of an ordering on the 1-cells of $\lambda$ gives us a natural isomorphism of $\mathcal{H}(\lambda; \mathbb{Z})$ with $\mathbb{Z}^n$ and through it we can define a bilinear skew form on $\mathcal{H}(\lambda; \mathbb{Z})$ given by

$$\sigma\left(\sum_i a_i \lambda_i, \sum_j b_i \lambda_j\right) := \sum_{i,j} a_i b_j \sigma_{ij} \qquad (1)$$

Observe that $\sigma_{ij}$ is determined by the mutual positions of $\lambda_i$ and $\lambda_j$, hence $\sigma$ is indeed independent from the choice of an ordering on $\lambda$ and it is a bilinear skew form intrinsically defined on $\mathcal{H}(\lambda; \mathbb{Z})$.

Now we fix an ordering in $\lambda$ and we choose $\alpha = (\alpha_1, \ldots, \alpha_n)$ and $\beta = (\beta_1, \ldots, \beta_n)$ in $\mathbb{Z}^n \cong \mathcal{H}(\lambda; \mathbb{Z})$. We can associate with $\alpha$ a monomial in $\mathcal{T}_\lambda^q$, that we briefly denote by $X^\alpha$, defined by

$$X^\alpha := X_1^{\alpha_1} \cdots X_n^{\alpha_n}$$

Introduce also the following notation

$$\underline{X}^\alpha := q^{-\sum_{i<j} \alpha_i \alpha_j \sigma_{ij}} X^\alpha = q^{-\sum_{i<j} \alpha_i \alpha_j \sigma_{ij}} X_1^{\alpha_1} \cdots X_n^{\alpha_n}$$

**Lemma 1.3.** For every $\alpha, \beta \in \mathbb{Z}^n \cong \mathcal{H}(\lambda; \mathbb{Z})$, the following relations hold in $\mathcal{T}_\lambda^q$:

$$X^\alpha X^\beta = q^{2\sigma(\alpha,\beta)} X^\beta X^\alpha = q^{2(\sum_{i>j} \alpha_i \beta_j \sigma_{ij})} X^{\alpha+\beta} \qquad (2)$$
$$\underline{X}^\alpha \underline{X}^\beta = q^{2\sigma(\alpha,\beta)} \underline{X}^\beta \underline{X}^\alpha = q^{\sigma(\alpha,\beta)} \underline{X}^{\alpha+\beta} \qquad (3)$$

Furthermore, for every permutation $\tau \in \mathfrak{S}_l$

$$q^{-\sum_{h<k} \sigma_{i_h i_k}} X_{i_1} \cdots X_{i_l} = q^{-\sum_{h<k} \sigma_{i_{\tau(h)} i_{\tau(k)}}} X_{i_{\tau(1)}} \cdots X_{i_{\tau(l)}} \qquad (4)$$

*Proof.* The first relations easily follow by direct calculations, it is sufficient to control how the coefficients change by the appropriate permutation. We will see how to deduce the relations in 3 using 2. The first equality is obvious, it is enough to multiply both the members of 2 by $q^{-\sum_{i<j} \alpha_i \alpha_j \sigma_{ij}} q^{-\sum_{i<j} \beta_i \beta_j \sigma_{ij}}$. We would like to show that $\underline{X}^\alpha \underline{X}^\beta = q^{\sigma(\alpha,\beta)} \underline{X}^{\alpha+\beta}$ holds. By virtue of 2, it is sufficient to prove that the exponents of $q$, multiplying the elements $X^\alpha X^\beta$ and $q^{2(\sum_{i>j} \alpha_i \beta_j \sigma_{ij})} X^{\alpha+\beta}$ respectively, coincide. It is simple to show we can reduce to prove the following equality

$$-\sum_{i<j} \alpha_i \alpha_j \sigma_{ij} - \sum_{i<j} \beta_i \beta_j \sigma_{ij} = -\sum_{i<j} (\alpha+\beta)_i (\alpha+\beta)_j \sigma_{ij} + \sigma(\alpha, \beta) - 2\sum_{i>j} \alpha_i \beta_j \sigma_{ij}$$

Using the fact that $\sigma$ is skew-symmetric, we deduce

$$-\sum_{i<j}(\alpha+\beta)_i(\alpha+\beta)_j\sigma_{ij} + \sigma(\alpha,\beta) - 2\sum_{i>j}\alpha_i\beta_j\sigma_{ij} =$$
$$= -\sum_{i<j}(\alpha+\beta)_i(\alpha+\beta)_j\sigma_{ij} + \sum_{i<j}\alpha_i\beta_j\sigma_{ij} + \sum_{i>j}\alpha_i\beta_j\sigma_{ij} - \sum_{i>j}\alpha_i\beta_j\sigma_{ij} +$$
$$- \sum_{i>j}\alpha_i\beta_j\sigma_{ij}$$
$$= -\sum_{i<j}(\alpha+\beta)_i(\alpha+\beta)_j\sigma_{ij} + \sum_{i<j}\alpha_i\beta_j\sigma_{ij} + \sum_{i<j}\beta_i\alpha_j\sigma_{ij}$$
$$= -\sum_{i<j}\alpha_i\alpha_j\sigma_{ij} - \sum_{i<j}\beta_i\beta_j\sigma_{ij}$$

that concludes the proof.

For what concerns the relation 4, it is simple to prove it when $\tau$ is a transposition of consecutive elements and the general case easily follows. □

## 1.2 Representations of the Chekhov-Fock algebra

Following Definition 1.2, we see that, ordering the edges with respect to the orientation as in Figure 1, the Chekhov-Fock algebra associated with an ideal triangle is isomorphic to $\mathbb{C}[X_1, X_2, X_3]$ endowed with the relations

$$X_i X_{i+1} = q^2 X_{i+1} X_i$$

for every $i = 1, 2, 3$, where the indices are mod (3). Unless specified differently, we will always assume to be in this situation, with the edges ordered clockwise with respect to the orientation.

**Proposition 1.4.** Let $\rho\colon \mathcal{T}_T^q \to \mathrm{End}(V)$ be an irreducible representation, with $V$ a $\mathbb{C}$-vector space of dimension $d$ and with $\mathcal{T}_T^q$ denoting the Chekhov-Fock algebra associated with the ideal triangle $T$, which admits a unique ideal triangulation. Then $d = N$ and there exist $x_1, x_2, x_3, h \in \mathbb{C}^*$ such that

$$\rho(X_i^N) = x_i\, id_V \qquad \text{for } i = 1, 2, 3$$
$$\rho(H) = h\, id_V$$

where the $X_i$ denote the generators associated with the edges of $T$ and $H = q^{-1}X_1 X_2 X_3$. In addition, the following holds

$$h^N = x_1 x_2 x_3$$

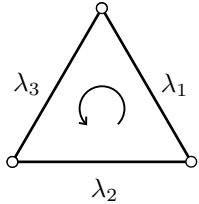

Figure 1: Standard indexing on a triangle

Moreover, two irreducible representations $\rho\colon \mathcal{T}_T^q \to \mathrm{End}(V)$ and $\rho'\colon \mathcal{T}_T^q \to \mathrm{End}(V')$ are isomorphic if and only if $x_i = x_i'$ and $h = h'$, where the $x_i, h$ are the above quantities related to $\rho$ and the $x_i', h'$ are related to $\rho'$. Furthermore, for every choice of values $x_1, x_2, x_3, h \in \mathbb{C}^*$ verifying $h^N = x_1 x_2 x_3$, there exists a representation, unique up to isomorphism, realizing them as invariants.

*Proof.* See [BBL07, Lemma 2]. $\square$

*Remark* 1.5. Following the proof of the previous Proposition, we see that, for every irreducible representation $\rho\colon \mathcal{T}_T^q \to \mathrm{End}(V)$, there exists a basis of $V$ in which the elements $\rho(X_1), \rho(X_2), \rho(X_3)$ are respectively represented by the following matrices:

$$y_1 B_1 = y_1 \begin{pmatrix} 1 & & & \\ & q^2 & & \\ & & \ddots & \\ & & & q^{2(N-1)} \end{pmatrix}$$

$$y_2 B_2 = y_2 \begin{pmatrix} 0 & \cdots & 0 & 1 \\ & & & 0 \\ & I_{N-1} & & \vdots \\ & & & 0 \end{pmatrix}$$

$$y_3 B_3 = y_3 \begin{pmatrix} 0 & q^{1-2(2-1)} & & \\ \vdots & & \ddots & \\ 0 & & & q^{1-2(N-1)} \\ q & 0 & \cdots & 0 \end{pmatrix}$$

where $y_i$ is an $N$-th root of $x_i$ for every $i = 1, 2, 3$ and $h = y_1 y_2 y_3$. When a representation $\rho\colon \mathcal{T}_T^q \to \mathrm{End}(V)$ has this form, that is $\rho(X_i) = y_i B_i \in \mathrm{End}(\mathbb{C}^N)$, we will say that $\rho$ is in *standard form*.

The standard form is not unique. More precisely, for every choice of $N$-th roots $y_i$ of the $x_i$ such that $h = y_1 y_2 y_3$, there exists a basis $\mathcal{B}$, unique up to common scalar multiplication, such that $\rho$ is represented as described above with respect to $\mathcal{B}$.

**Definition 1.6.** When two linear isomorphisms $A, B$ differ by a non-zero scalar, we will briefly write $A \doteq B$.

### 1.2.1 Fusion of surfaces

Let $S$ be a surface and let $\lambda$ be an ideal triangulation of it. We can construct, starting from $S$ and the choice of certain internal 1-cells $\lambda_{i_1}, \ldots, \lambda_{i_h}$, a surface $R$ obtained by splitting $S$ along these selected edges $\lambda_{i_1}, \ldots, \lambda_{i_h}$, i. e. by removing the identifications along $\lambda_{i_j}$ between the edges of the ideal triangles having $\lambda_{i_j}$ as edge. On $R$ we can clearly find an ideal triangulation $\mu$ and an orientation induced by $\lambda$ and by the orientation on $S$. In these circumstances, we will say that $S$ is obtained from $R$ *by fusion*, and analogously that the ideal triangulation $\lambda \in \Lambda(S)$ is obtained from $\mu \in \Lambda(R)$ *by fusion*.

Fix an indexing on the 1-cells of $\lambda$ and one on the 1-cells of $\mu$. We can construct a natural inclusion $j_{\mu\lambda} \colon \mathcal{H}(\lambda; \mathbb{Z}) \to \mathcal{H}(\mu; \mathbb{Z})$, which associates to each $e_i \in \mathbb{Z}^n \cong \mathcal{H}(\lambda; \mathbb{Z})$, corresponding to the edge $\lambda_i$, the vector $v_i = (v_{i1}, \ldots, v_{is}) \in \mathbb{Z}^s \cong \mathcal{H}(\mu; \mathbb{Z})$, defined by

$$v_{ik} := \begin{cases} 1 & \text{if } \mu_k \text{ goes by fusion in } \lambda_i \\ 0 & \text{otherwise} \end{cases}$$

It is easy to verify that map $j_{\mu\lambda}$ is indeed an inclusion. In fact, every vector $v_i$ has at most two non-zero components and, if $i \neq j$, the supports of the vectors $v_i$ and $v_j$ are disjoint. Denote by $\sigma$ and $\eta$ the skew-symmetric bilinear forms respectively on $\mathcal{H}(\lambda; \mathbb{Z})$ and $\mathcal{H}(\mu; \mathbb{Z})$, defined as in Subsection 1.1. Then the following holds

$$\sigma(v, w) = \eta(j_{\mu\lambda}(v), j_{\mu\lambda}(w))$$

In order to prove it, it is sufficient to verify the equality on a set of generators, so it is enough to prove that, for every $i, j$, we have $\sigma_{ij} = \eta(v_i, v_j)$. Recalling the definition, we can express $\sigma_{ij}$ as follows

$$\sum_{\substack{T \text{ triangle} \\ \text{of } \lambda}} \sum_{\substack{c \text{ spike} \\ \text{of } T}} s_\lambda(c, \lambda_i, \lambda_j)$$

where $s_\lambda(c, \lambda_i, \lambda_j)$ is equal to $+1$ if $c$ has $\lambda_i$ on the left and $\lambda_j$ on the right, $-1$ if $c$ has $\lambda_j$ on the left and $\lambda_i$ on the right, and is equal to $0$ otherwise. Suppose that the edge $\lambda_i$ is the result of the identification of the edges $\mu_{i_1}, \mu_{i_2}$ in $\mu$, and that analogously $\lambda_j$ is the result of the identification of the edges $\mu_{j_1}, \mu_{j_2}$. The faces of the ideal triangulation $\lambda$ are in natural bijection with the faces of $\mu$, and consequently the respective spikes too. Hence it is sufficient to prove that

$$s_\lambda(c, \lambda_i, \lambda_j) = s_\mu(c, \mu_{i_1}, \mu_{j_1}) + s_\mu(c, \mu_{i_1}, \mu_{j_2}) + s_\mu(c, \mu_{i_2}, \mu_{j_1}) + s_\mu(c, \mu_{i_2}, \mu_{j_2})$$

for every spike $c$ of the ideal triangulation. In the right member there exists at most one non-zero term and it is immediate to see that, thanks to the coherent choice of orientations on the faces of both triangulations, the equality holds.

Denote now by $X_1, \ldots, X_n$ the generators of $\mathcal{T}_\lambda^q$ associated with the edges of $\lambda$ and by $Y_1, \ldots, Y_s$ the ones of $\mathcal{T}_\mu^q$ associated with the edges of $\mu$. Define the map

$$\iota_{\mu\lambda} \colon \begin{array}{ccc} \mathcal{T}_\lambda^q & \longrightarrow & \mathcal{T}_\mu^q \\ \underline{X}^\alpha & \longmapsto & \underline{Y}^{j_{\mu\lambda}(\alpha)} \end{array}$$

Using the relation 3 and what just shown, we can verify that

$$\begin{aligned} \iota_{\mu\lambda}(\underline{X}^\alpha \underline{X}^\beta) &= \iota_{\mu\lambda}(q^{\sigma(\alpha,\beta)} \underline{X}^{\alpha+\beta}) = q^{\sigma(\alpha,\beta)} \underline{Y}^{j_{\mu\lambda}(\alpha+\beta)} \\ &= q^{\eta(j_{\mu\lambda}(\alpha), j_{\mu\lambda}(\beta))} \underline{Y}^{j_{\mu\lambda}(\alpha) + j_{\mu\lambda}(\beta)} = \underline{Y}^{j_{\mu\lambda}(\alpha)} \underline{Y}^{j_{\mu\lambda}(\beta)} \\ &= \iota_{\mu\lambda}(\underline{X}^\alpha) \iota_{\mu\lambda}(\underline{X}^\beta) \end{aligned}$$

The injectivity of $j_{\mu\lambda} \colon \mathcal{H}(\lambda; \mathbb{Z}) \to \mathcal{H}(\mu; \mathbb{Z})$ immediately implies the injectivity of $\iota_{\mu\lambda}$. Hence we have proved that, if $(S, \lambda)$ is obtained from $(R, \mu)$ by fusion, then there exists an inclusion of algebras $\iota_{\mu\lambda}\colon \mathcal{T}^q_\lambda \to \mathcal{T}^q_\mu$, carrying monomials of $\mathcal{T}^q_\lambda$ in monomials of $\mathcal{T}^q_\mu$. Moreover, these maps, for varying $S$, $R$, $\lambda$ and $\mu$ as above, verify a sort of composition property. More precisely, assume that $(S, \lambda)$, $(S', \lambda')$ and $(S'', \lambda'')$ are three surfaces endowed with ideal triangulations, with $(S, \lambda)$ obtained from $(S', \lambda')$ by fusion and with $(S', \lambda')$ obtained from $(S'', \lambda'')$ by fusion. Then, by definition, the homomorphisms $j_{\lambda''\lambda}\colon \mathcal{H}(\lambda; \mathbb{Z}) \to \mathcal{H}(\lambda''; \mathbb{Z})$ and $j_{\lambda''\lambda'} \circ j_{\lambda'\lambda}\colon \mathcal{H}(\lambda; \mathbb{Z}) \to \mathcal{H}(\lambda''; \mathbb{Z})$ coincide, so on the Chekhov-Fock algebras the following relation holds

$$\iota_{\lambda''\lambda} = \iota_{\lambda''\lambda'} \circ \iota_{\lambda'\lambda} \tag{5}$$

Now let us focus on a more specific situation. Given $S$ a surface as above and $\lambda \in \Lambda(S)$ an ideal triangulation of it, $S$ can be obtained by fusion from a surface $S_0$ realized by splitting $S$ along all its internal edges or, in other words, realized as the disjoint union of the triangles $T_i$ composing $\lambda$, for $i = 1, \ldots, m$. $S_0$ admits a unique ideal triangulation $\lambda_0$ and its Chekhov-Fock algebra associated with $\lambda_0$ is naturally isomorphic to

$$\mathcal{T}^q_{T_1} \otimes \cdots \otimes \mathcal{T}^q_{T_m}$$

In this case we will denote simply by $\iota_\lambda$ the inclusion map $\iota_{\lambda_0\lambda}\colon \mathcal{T}^q_\lambda \to \bigotimes_i \mathcal{T}^q_{T_i}$, defined as above.

Each triangle $T_i$ is endowed with an orientation, determined by the one on $S$. Order the edges $\lambda^{(i)}_1, \lambda^{(i)}_2, \lambda^{(i)}_3$ of $T_i \subseteq S_0$ clockwise, as in Figure 1, and denote respectively by $X^{(i)}_1, X^{(i)}_2, X^{(i)}_3$ the generators of $\mathcal{T}^q_\lambda$ associated with these edges. For the sake of simplicity, given $F^{(1)} \otimes \cdots \otimes F^{(m)}$ a monomial in $\bigotimes_i \mathcal{T}^q_{T_i}$, we will omit in it the tensor product by those terms verifying $F^{(i)} = 1$. Recalling the definition of $\iota_\lambda$, we want to give an explicit description of the image, under this inclusion, of the generators $X_i$ associated with the edges of the ideal triangulation $\lambda$ of $S$:

- if $\lambda_i$ is a boundary edge, then it is side of a single ideal triangle $T_{k_i}$, so we have $\iota_\lambda(X_i) = X^{(k_i)}_{a_i} \in \bigotimes_i \mathcal{T}^q_{T_i}$, where $a_i$ is the index of the edge of $T_{k_i}$ identified in $S$ with $\lambda_i$;

- if $\lambda_i$ is an internal edge and it is side of two distinct triangles $T_{l_i}$ and $T_{r_i}$, then we have $\iota_\lambda(X_i) = X^{(l_i)}_{a_i} \otimes X^{(r_i)}_{b_i} \in \bigotimes_i \mathcal{T}^q_{T_i}$, where $a_i$ and $b_i$ are the indices of the edges of $T_{l_i}$ and $T_{r_i}$, respectively, identified in $S$ with $\lambda_i$;

- if $\lambda_i$ is an internal edge and it is side of a single triangle $T_{k_i}$, then we have $\iota_\lambda(X_i) = q^{-1} X^{(k_i)}_{a_i} X^{(k_i)}_{b_i} = q\, X^{(k_i)}_{b_i} X^{(k_i)}_{a_i} \in \bigotimes_i \mathcal{T}^q_{T_i}$, where $a_i$ and $b_i$ are the indices of the edges of $T_{k_i}$ identified in $S$ with $\lambda_i$ and $\lambda^{(k_i)}_{a_i}$, $\lambda^{(k_i)}_{b_i}$ lie on the left and on the right, respectively, of their common spike.

### 1.2.2 Local representations

**Definition 1.7.** Given $S$ a surface and $\lambda$ an ideal triangulation of $S$, we denote by $\Gamma_{S,\lambda}$ the dual graph of $\lambda$. More precisely, $\Gamma_{S,\lambda}$ is a CW-complex of dimension 1, whose vertices $T^*_b$ correspond to the triangles $T_b$ of $\lambda$ and, for every $\lambda_a$ internal edge, there is a 1-cell $\lambda^*_a$ that connects the vertices corresponding to the triangles on the sides of $\lambda_a$, even if the triangles coincide.

It follows from the definition that all the vertices have valency $\leq 3$. In particular, the valency of a vertex $T_b^*$ in $\Gamma_{S,\lambda}$ is equal to the number of sides of $T_b$ that are internal edges in $\lambda$.

Given $\Gamma$ a graph, we will denote with $\Gamma^{(k)}$ the set of the $k$-cells composing $\Gamma$.

Fix $\lambda \in \Lambda(S)$ an ideal triangulation of $S$ and, for every edge $\lambda_i$ of $\lambda$, choose an arbitrary orientation on it. Now orient the edges of $\Gamma = \Gamma_{S,\lambda}$ as follows: the 1-cell $\lambda_i^*$, dual of the internal edge $\lambda_i$, is oriented in such a way that the intersection number $i(\lambda_i, \lambda_i^*)$ in $S$ is equal to $+1$ (remember that we are considering oriented surfaces). Moreover, we assume that all the vertices $T_l^*$ have positive sign.

Given $a = 1, \ldots, n$ with $\lambda_a$ an internal edge having two different triangles on its sides, we define

$$\varepsilon(a,b) := \begin{cases} +1 & \textit{if } T_b \textit{ is on the left of } \lambda_a \\ -1 & \textit{if } T_b \textit{ is on the right of } \lambda_a \\ 0 & \textit{otherwise} \end{cases}$$

for every $b = 1, \ldots, m$. If $\lambda_a$ is internal and it has the same triangle on its sides, we define $\varepsilon(a,b) = 0$ for every $b = 1, \ldots, m$. Observe that a triangle $T_b$ is on the left of $\lambda_a$ if and only if the previously fixed orientation on $\lambda_a$ coincides with the orientation determined as boundary of $T_b$. Moreover, it is immediate to see that the definition of $\varepsilon(a,b)$ can be reformulated as follows

$$\varepsilon(a,b) := \begin{cases} +1 & \textit{if } \lambda_a^* \textit{ goes towards } T_b^* \\ -1 & \textit{if } \lambda_a^* \textit{ comes from } T_b^* \\ 0 & \textit{otherwise} \end{cases}$$

if $\lambda_a^*$ has different ends, otherwise $\varepsilon(a,b) = 0$ for every $b = 1, \ldots, m$.

Fixed an orientation on an ideal triangulation $\lambda$, the definition of local representation of $\mathcal{T}_\lambda^q$ given in [BBL07] can be reformulated as follows:

**Definition 1.8.** Let $\lambda \in \Lambda(S)$ be an ideal triangulation, with triangles labelled as $T_1, \ldots, T_m$, and let $(\rho_1, \ldots, \rho_m), (\overline{\rho}_1, \ldots, \overline{\rho}_m)$ be two $m$-tuples in which, for every $j = 1, \ldots, m$, $\rho_j \colon \mathcal{T}_{T_j}^q \to \mathrm{End}(V_j)$ and $\overline{\rho}_j \colon \mathcal{T}_{T_j}^q \to \mathrm{End}(W_j)$ are irreducible representations of the triangle algebra $\mathcal{T}_{T_j}^q$. The elements $(\rho_1, \ldots, \rho_m)$, $(\overline{\rho}_1, \ldots, \overline{\rho}_m)$ are *locally equivalent* if the following hold

- for every $j = 1, \ldots, m$ the vector spaces $V_j$ and $W_j$ are equal;

- for every $i = 1, \ldots, n$, we have:

  - if $\lambda_i$ is a boundary edge, then there is a unique triangle $T_{s_i}$ that has $\lambda_i$ on its side. In this case we ask that

    $$\rho_{s_i}(X_{a_i}^{(s_i)}) = \overline{\rho}_{s_i}(X_{a_i}^{(s_i)})$$

    where $a_i$ is the index of the edge in $T_{s_i}$ that is identified to $\lambda_i$;

  - if $\lambda_i$ is an internal edge and $T_{l_i}$, $T_{r_i}$ are distinct triangles on the left and on the right, respectively, of $\lambda_i$ (recall that we have fixed orientations on the edges $\lambda_i$), then there exists $\alpha_i \in \mathbb{C}^*$ such that

$$\rho_{l_i}(X_{a_i}^{(l_i)}) = \alpha_i^{+1} \overline{\rho}_{l_i}(X_{a_i}^{(l_i)}) = \alpha_i^{\varepsilon(i,l_i)} \overline{\rho}_{l_i}(X_{a_i}^{(l_i)})$$
$$\rho_{r_i}(X_{b_i}^{(r_i)}) = \alpha_i^{-1} \overline{\rho}_{r_i}(X_{b_i}^{(r_i)}) = \alpha_i^{\varepsilon(i,r_i)} \overline{\rho}_{r_i}(X_{b_i}^{(r_i)})$$

where $a_i$ and $b_i$ are the indices of the edges in $T_{l_i}$ and $T_{r_i}$, respectively, that are identified to $\lambda_i$;

– if $\lambda_i$ is an internal edge and it has the triangle $T_{k_i}$ on both sides, then there exists $\alpha_i \in \mathbb{C}^*$ such that

$$\rho_{k_i}(X_{a_i}^{(k_i)}) = \alpha_i^{+1} \overline{\rho}_{k_i}(X_{a_i}^{(k_i)})$$
$$\rho_{k_i}(X_{b_i}^{(k_i)}) = \alpha_i^{-1} \overline{\rho}_{k_i}(X_{b_i}^{(k_i)})$$

where $a_i$ and $b_i$ are the indices of the edges in $T_{k_i}$ that are identified to $\lambda_i$ and the $a_i$-th side, unlike the $b_i$-th one, has the orientation as boundary of $T_{k_i}$ coherent with the orientation of $\lambda_i$.

A *local representation* $\rho$ of $\mathcal{T}_\lambda^q$ is a local equivalence class of $m$-tuples of representations $(\rho_1, \ldots, \rho_m)$ as above.

*Remark* 1.9. Given $[\rho_1, \ldots, \rho_m]$ a local representation of $\mathcal{T}_\lambda^q$, we can define a representation of $\mathcal{T}_\lambda^q$ as follows

$$\rho := (\rho_1 \otimes \cdots \otimes \rho_m) \circ \iota_\lambda : \mathcal{T}_\lambda^q \longrightarrow \mathrm{End}(V_1 \otimes \cdots \otimes V_m)$$

By definition of the locally equivalence relation between $m$-tuples of representations, this $\rho$ does not depend on the choice of the representative of the equivalence class $[\rho_1, \ldots, \rho_m]$.

We will confuse a representative $(\rho_1, \ldots, \rho_m)$ of a local representation $\rho$ with the obvious corresponding representation $\rho_1 \otimes \cdots \otimes \rho_m$ on the Chekhov-Fock algebra $\mathcal{T}_{\lambda_0}^q$ of the surface $S_0$, where $S_0$ is the surface obtained by splitting $S$ along all its edges and $\lambda_0$ is its ideal triangulation.

**Definition 1.10.** Let $\lambda \in \Lambda(S)$ be an ideal triangulation of $S$ and let $\rho \colon \mathcal{T}_\lambda^q \to \mathrm{End}(V)$ be a local representation of $\mathcal{T}_\lambda^q$. We denote by $\mathscr{F}_{S_0}(\rho)$ the set of representatives of $\rho$ as local representation, which are local (and irreducible, see Proposition 1.14) representations of the Chekhov-Fock algebra $\mathcal{T}_{\lambda_0}^q$ of the surface $S_0$, obtained by splitting $S$ along $\lambda$. Moreover, fixed an orientation on $\lambda$ and given $\rho_1 \otimes \cdots \otimes \rho_m$ and $\overline{\rho}_1 \otimes \cdots \otimes \overline{\rho}_m$ two elements of $\mathscr{F}_{S_0}(\rho)$, we will write

$$\overline{\rho}_1 \otimes \cdots \otimes \overline{\rho}_m \xrightarrow{\alpha_i} \rho_1 \otimes \cdots \otimes \rho_m$$

if $\rho_1 \otimes \cdots \otimes \rho_m$ and $\overline{\rho}_1 \otimes \cdots \otimes \overline{\rho}_m$ are related by the numbers $(\alpha_i)_i$ as described in Definition 1.8. The $(\alpha_i)_i$ are called the *transition constants* from $\overline{\rho}_1 \otimes \cdots \otimes \overline{\rho}_m$ to $\rho_1 \otimes \cdots \otimes \rho_m$.

It is very simple to see that, with these notations introduced, the following holds:

**Lemma 1.11.** *Let* $\zeta = \bigotimes_i \rho_i$, $\zeta' = \bigotimes_i \rho'_i$ *and* $\zeta'' = \bigotimes_i \rho''_i$ *be three representatives of a local representation* $\rho$. *Then the following properties hold:*

1. there exists a unique collection of transition constants, depending on the chosen orientation on $\lambda$, such that $\zeta \xrightarrow{\alpha_i} \zeta'$;

2. if $\zeta \xrightarrow{\alpha_i} \zeta'$ and $\zeta' \xrightarrow{\beta_i} \zeta''$, then $\zeta \xrightarrow{\alpha_i \beta_i} \zeta''$

3. if $\zeta \xrightarrow{\alpha_i} \zeta'$, then $\zeta' \xrightarrow{\alpha_i^{-1}} \zeta$. $\square$

**Definition 1.12.** Given $S$ a surface as above, $\lambda \in \Lambda(S)$ an ideal triangulation and two local representations $[\rho_1, \ldots, \rho_m], [\rho'_1, \ldots, \rho'_m]$ of $\mathcal{T}^q_\lambda$, with

$$\rho_j : \mathcal{T}^q_{T_j} \longrightarrow \mathrm{End}(V_j)$$
$$\rho'_j : \mathcal{T}^q_{T_j} \longrightarrow \mathrm{End}(V'_j)$$

we will say that $[\rho_1, \ldots, \rho_m]$ and $[\rho'_1, \ldots, \rho'_m]$ are isomorphic if there exist respectively representatives $(\rho_1, \ldots, \rho_m)$ and $(\rho'_1, \ldots, \rho'_m)$ of them and there exist linear isomorphisms $L_j \colon V_j \to V'_j$ such that, for every $j = 1, \ldots, m$ we have

$$L_j \circ \rho_j(X) \circ L_j^{-1} = \rho'_j(X) \qquad\qquad \forall X \in \mathcal{T}^q_{T_j}$$

Assume that $[\rho_1, \ldots, \rho_m]$ e $[\rho'_1, \ldots, \rho'_m]$ are isomorphic and let $(\rho_1, \ldots, \rho_m)$ and $(\rho'_1, \ldots, \rho'_m)$ be two representatives of them such that there exist linear isomorphisms $L_j \colon V_j \to V'_j$ with $L_j \circ \rho_j \circ L_j^{-1} = \rho'_j$. Then, for any other choice of a representative $(\bar{\rho}_1, \ldots, \bar{\rho}_m)$ of $[\rho_1, \ldots, \rho_m]$, the $m$-tuple of representations $\bar{\rho}'_j := L_j \circ \bar{\rho}_j \circ L_j^{-1}$ is $\sim_S$-locally equivalent to $(\rho'_1, \ldots, \rho'_m)$. From this fact immediately follows that the isomorphism relation defined above is indeed an equivalence relation.

Now we recall some results of [BBL07] that will be useful in the following analysis.

**Lemma 1.13.** Let $[\rho_1, \ldots, \rho_m]$ be a local representation of $\mathcal{T}^q_\lambda$. Then, for every generator $X_i \in \mathcal{T}^q_\lambda$ associated with the edge $\lambda_i$, the representation $\rho_\lambda := (\rho_1 \otimes \cdots \otimes \rho_m) \circ \iota_\lambda$ verifies

$$\rho(X_i^N) = x_i \, id_{V_1 \otimes \cdots \otimes V_m}$$

for a certain $x_i \in \mathbb{C}^*$. In addition, there exists $h \in \mathbb{C}^*$ such that

$$\rho(H) = h \, id_{V_1 \otimes \cdots \otimes V_m}$$

*Proof.* See [BBL07, Lemma 5]. $\square$

The following result was firstly stated in [BBL07], for the sake of completeness we provide a proof in Appendix A.

**Proposition 1.14.** Let $S$ be an ideal polygon with $p \geq 3$ vertices. Then, for every $\lambda$ triangulation of $S$ and for every local representation $[\rho_1, \ldots, \rho_m]$ of the Chekhov-Fock algebra $\mathcal{T}^q_\lambda$, the representation $\rho := (\rho_1 \otimes \cdots \otimes \rho_m) \circ \iota_\lambda$ is irreducible.

Denote by $\mathscr{R}_{loc}(\mathcal{T}^q_\lambda)$ the set of isomorphism classes, as local representations, of local representations of $\mathcal{T}^q_\lambda$.

**Theorem 1.15.** Let $S$ be a surface and fix $q \in \mathbb{C}^*$ a primitive $N$-th root of $(-1)^{N+1}$. Then the set $\mathscr{R}_{loc}(\mathcal{T}_\lambda^q)$ is in bijection with the set of elements $((x_i)_i; h)$ in $(\mathbb{C}^*)^n \times \mathbb{C}^*$ verifying
$$h^N = x_1 \cdots x_n$$
and the correspondence associates with the isomorphism class as local representations of $[\rho_1, \ldots, \rho_m]$, the element $((x_i)_i; h)$ defined by the relations
$$\rho(X_i^N) = x_i \, id_{V_1 \otimes \cdots \otimes V_m}$$
$$\rho(H) = h \, id_{V_1 \otimes \cdots \otimes V_m}$$
for $i = 1, \ldots, n$, where $\rho = (\rho_1 \otimes \cdots \otimes \rho_m) \circ \iota_\lambda$.

*Proof.* See [BBL07, Proposition 6]. □

**Corollary 1.16.** Let $[\rho_1, \ldots, \rho_m]$ and $[\rho'_1, \ldots, \rho'_m]$ be two local representations of $\mathcal{T}_\lambda^q$. Then they are isomorphic to each other as local representations if and only if $\rho := (\rho_1 \otimes \cdots \otimes \rho_m) \circ \iota_\lambda$ and $\rho' := (\rho'_1 \otimes \cdots \otimes \rho'_m) \circ \iota_\lambda$ are isomorphic as representations. □

### 1.2.3 An action of $H_1(S; \mathbb{Z}_N)$

Denote by $(C_\bullet(\Gamma; \mathbb{Z}_N), \partial_\bullet)$ the cellular chain complex of $\Gamma = \Gamma_{S,\lambda}$. Then $C_0(\Gamma; \mathbb{Z}_N)$ is the $\mathbb{Z}_N$-module freely generated by the vertices $T_b^*$ of $\Gamma$ and $C_1(\Gamma; \mathbb{Z}_N)$ is the $\mathbb{Z}_N$-module freely generated by the oriented 1-cells $\lambda_a^*$ of $\Gamma$. Because $\Gamma$ has dimension 1, all the other $C_i(\Gamma; \mathbb{Z}_N)$ are equal to zero. Thanks to what observed, we can describe the boundary $\partial_1$ in terms of the triangulation $\lambda$. Given $\lambda_a^*$ a 1-cell of $\Gamma$, the boundary $\partial_1(\lambda_a^*)$ verifies
$$\partial_1(\lambda_a^*) = \sum_{b=1}^m \varepsilon(a,b) \, T_b^* \in C_0(\Gamma; \mathbb{Z}_N)$$
where the $\varepsilon(a,b)$ are the numbers defined in 1.2.2. Hence the first group of cellular homology $H_1(\Gamma; \mathbb{Z}_N)$ is equal to the subgroup $\operatorname{Ker} \partial_1$ of $C_1(\Gamma; \mathbb{Z}_N)$, whose elements are the $\mathbb{Z}_N$-combinations $\sum_{a=1}^n c_a \lambda_a^*$ such that, for every $b = 1, \ldots, m$ the following holds:
$$\sum_{a=1}^n \varepsilon(a,b) \, c_a = 0 \in \mathbb{Z}_N$$
The dual graph $\Gamma$ is a deformation retract of $S$, so the group $H_1(\Gamma; \mathbb{Z}_N)$ can be identified to $H_1(S; \mathbb{Z}_N)$ via a certain inclusion of $\Gamma$ in $S$, which is well defined up to homotopy.

Given $\lambda \in \Lambda(S)$ an oriented ideal triangulation and $\rho \colon \mathcal{T}_\lambda^q \to \operatorname{End}(V)$ a local representation, we can define an action of $H_1(S; \mathbb{Z}_N)$ on the set $\mathscr{F}_{S_0}(\rho)$ as follows:

**Definition 1.17.** Given $c = \sum_i c_i \lambda_i^*$ an element of $H_1(\Gamma; \mathbb{Z}_N) \cong H_1(S; \mathbb{Z}_N)$ and fixed a representative $\rho_1 \otimes \cdots \otimes \rho_m$ of $\rho$, we can produce another $m$-tuple of representations $\overline{\rho}_1 \otimes \cdots \otimes \overline{\rho}_m$ of $\rho$ defined as follows: for every $l = 1, \ldots, m$

- if $T_l$ has distinct sides, labelled clockwise in $\lambda$ as $\lambda_i$, $\lambda_j$, $\lambda_k$, then $\overline{\rho}_l$ is equal to

$$\begin{aligned}
\overline{\rho}_l(X_1^{(l)}) &:= q^{2c_i\varepsilon(i,l)}\rho_l(X_1^{(l)}) \\
\overline{\rho}_l(X_2^{(l)}) &:= q^{2c_j\varepsilon(j,l)}\rho_l(X_2^{(l)}) \\
\overline{\rho}_l(X_3^{(l)}) &:= q^{2c_k\varepsilon(k,l)}\rho_l(X_3^{(l)})
\end{aligned} \qquad (6)$$

  where the edges $\lambda_i$, $\lambda_j$, $\lambda_k$ of $T_l$ respectively correspond to the variables $X_1^{(l)}, X_2^{(l)}, X_3^{(l)} \in \mathcal{T}_{T_l}^q$;

- if $T_l$ has two sides that are identified to $\lambda_i$ in $\lambda$ and the remaining one is identified with $\lambda_j$, then $\overline{\rho}_l$ is equal to

$$\begin{aligned}
\overline{\rho}_l(X_1^{(l)}) &:= q^{+2c_i}\rho_l(X_1^{(l)}) \\
\overline{\rho}_l(X_2^{(l)}) &:= q^{-2c_i}\rho_l(X_2^{(l)}) \\
\overline{\rho}_l(X_3^{(l)}) &:= q^{2c_j\varepsilon(j,l)}\rho_l(X_3^{(l)})
\end{aligned} \qquad (7)$$

  where the variables $X_1^{(l)}, X_2^{(l)}$ correspond to the sides of $T_l$ identified to $\lambda_i$ and $X_1^{(l)}$ is associated with the side having its boundary orientation coherent with the orientation of $\lambda_i$ (effectively the number $c_j$ is necessarily equal to zero, because $c$ is a cycle in $C_1(\Gamma; \mathbb{Z}_N)$). The other cases are treated in the same way.

It is immediate to see that, by construction, the representation $\overline{\rho}_1 \otimes \cdots \otimes \overline{\rho}_m$ is a representative of $\rho$ and, in the notations introduced above, the relations 6 can be summarized as

$$\rho_1 \otimes \cdots \otimes \rho_m \xrightarrow{q^{2c_i}} \overline{\rho}_1 \otimes \cdots \otimes \overline{\rho}_m$$

We will denote by $c \cdot (\rho_1 \otimes \cdots \otimes \rho_m)$ the representation $\overline{\rho}_1 \otimes \cdots \otimes \overline{\rho}_m$ constructed in this way.

Now we are able to enunciate the main result of this section, which will be the key ingredient in the following discussion concerning intertwining operators:

**Proposition 1.18.** Let $\rho$ be a local representation of $\mathcal{T}_\lambda^q$. Then there exists an action of $H_1(S; \mathbb{Z}_N)$ on $\mathscr{F}_{S_0}(\rho)$, which verifies:

- $\rho_1 \otimes \cdots \otimes \rho_m$ and $\overline{\rho}_1 \otimes \cdots \otimes \overline{\rho}_m$ are isomorphic as representations of $\mathcal{T}_{\lambda_0}^q$ if and only if there exists an element $c \in H_1(S; \mathbb{Z}_N)$ such that

$$c \cdot (\rho_1 \otimes \cdots \otimes \rho_m) = \overline{\rho}_1 \otimes \cdots \otimes \overline{\rho}_m$$

- the action is free, i. e. $c \cdot (\rho_1 \otimes \cdots \otimes \rho_m)$ is equal to $\rho_1 \otimes \cdots \otimes \rho_m$ if and only if $c = 0 \in H_1(S; \mathbb{Z}_N)$.

*Proof.* Let $\rho_1 \otimes \cdots \otimes \rho_m$ and $\overline{\rho}_1 \otimes \cdots \otimes \overline{\rho}_m$ be two representatives of $\rho$. Observe that $\rho_1 \otimes \cdots \otimes \rho_m$ and $\overline{\rho}_1 \otimes \cdots \otimes \overline{\rho}_m$, as representations of $\mathcal{T}_{\lambda_0}^q$ are isomorphic if and only if, for every $j = 1, \ldots, m$ $\rho_j$ and $\overline{\rho}_j$ are isomorphic.

Let us focus on a single triangle $T_l$ and observe how the representations $\rho_l \colon \mathcal{T}_{T_l}^q \to \mathrm{End}(V_l)$ and $\overline{\rho}_l \colon \mathcal{T}_{T_l}^q \to \mathrm{End}(V_l)$ differ. Label clockwise as $\lambda_i$, $\lambda_j$ and $\lambda_k$ the edges of $T_l$ in $\lambda$, and as $X_1^{(l)}$, $X_2^{(l)}$, $X_3^{(l)}$ the corresponding variables in $\mathcal{T}_{T_l}^q$. The relations between $\rho_l$ and $\overline{\rho}_l$ are the following

$$\rho_l(X_1^{(l)}) = \alpha_i^{\varepsilon(i,l)}\overline{\rho}_l(X_1^{(l)})$$
$$\rho_l(X_2^{(l)}) = \alpha_j^{\varepsilon(j,l)}\overline{\rho}_l(X_2^{(l)})$$
$$\rho_l(X_3^{(l)}) = \alpha_k^{\varepsilon(k,l)}\overline{\rho}_l(X_3^{(l)})$$

where we are assuming that the sides of $T_l$ are distinct in $\lambda$ and all internal, the other cases are treated in a similar way. Denote by $x_1^{(l)}$, $x_2^{(l)}$, $x_3^{(l)}$, $h^{(l)}$ the invariants of the irreducible representation $\rho_l$, and by $\overline{x}_1^{(l)}$, $\overline{x}_2^{(l)}$, $\overline{x}_3^{(l)}$, $\overline{h}^{(l)}$ the ones of $\overline{\rho}_l$. Then we deduce the following relations

$$x_1^{(l)} = \alpha_i^{N\varepsilon(i,l)}\overline{x}_1^{(l)}$$
$$x_2^{(l)} = \alpha_j^{N\varepsilon(j,l)}\overline{x}_2^{(l)}$$
$$x_3^{(l)} = \alpha_k^{N\varepsilon(k,l)}\overline{x}_3^{(l)}$$
$$h^{(l)} = \alpha_i^{\varepsilon(i,l)}\alpha_j^{\varepsilon(j,l)}\alpha_k^{\varepsilon(k,l)}\overline{h}^{(l)}$$

Now assume further that the representations $\rho_l$ and $\overline{\rho}_l$ are isomorphic. By virtue of Proposition 1.4, this is equivalent to ask that the invariants coincide. Then, for every triangle $T_l$, with edges labelled as before, we must have

$$\alpha_i^{N\varepsilon(i,l)} = 1$$
$$\alpha_j^{N\varepsilon(j,l)} = 1$$
$$\alpha_k^{N\varepsilon(k,l)} = 1$$
$$\alpha_i^{\varepsilon(i,l)}\alpha_j^{\varepsilon(j,l)}\alpha_k^{\varepsilon(k,l)} = 1$$

The first three equations can be rewritten as $\alpha_i^N = \alpha_j^N = \alpha_k^N = 1$ because $\varepsilon(a,l)$ is equal to $\pm 1$ for $a = i,j,k$. Hence, there exist $c_i, c_j, c_k \in \mathbb{Z}_N$ such that

$$\alpha_i = q^{2c_i}$$
$$\alpha_j = q^{2c_j}$$
$$\alpha_k = q^{2c_k}$$

Then the last condition can be rewritten as

$$c_i\,\varepsilon(i,l) + c_j\,\varepsilon(j,l) + c_k\,\varepsilon(k,l) = 0 \in \mathbb{Z}_N \tag{8}$$

Observe that the number $c_i\,\varepsilon(i,l) + c_j\,\varepsilon(j,l) + c_k\,\varepsilon(k,l)$ is exactly the coefficient of $T_l^*$ of the combination $\partial_1(\sum_a c_a\,\lambda_a^*)$. This relation must hold for every triangle $T_l$ in the ideal triangulation $\lambda$, so the element $\sum_a c_a\,\lambda_a^*$ is a cycle in $C_1(\Gamma; \mathbb{Z}_N)$, or equivalently it belongs to $H_1(\Gamma; \mathbb{Z}_N)$, and the representation $\rho_1 \otimes \cdots \otimes \rho_m$ coincides with $c \cdot (\overline{\rho}_1 \otimes \cdots \otimes \overline{\rho}_m)$.

Vice versa, if $\overline{\rho}_1 \otimes \cdots \otimes \overline{\rho}_m$ is equal to $c \cdot (\rho_1 \otimes \cdots \otimes \rho_m)$, then $\overline{\rho}_l$ is defined as in relation 6 or 7. Assume that the edges of $T_l$ are distinct, the other situations are analogous. Then we can easily see that $\rho_l$ is isomorphic to $\overline{\rho}_l$ for every

$l = 1, \ldots, m$. Indeed

$$\overline{\rho}_l((X_1^{(l)})^N) = \rho_l((X_1^{(l)})^N) \qquad q^{2Nc_i\varepsilon(i,l)} = 1$$
$$\overline{\rho}_l((X_2^{(l)})^N) = \rho_l((X_2^{(l)})^N) \qquad q^{2Nc_j\varepsilon(j,l)} = 1$$
$$\overline{\rho}_l((X_3^{(l)})^N) = \rho_l((X_3^{(l)})^N) \qquad q^{2Nc_k\varepsilon(k,l)} = 1$$
$$\overline{\rho}_l(H^{(l)}) = \rho_l(H^{(l)}) \qquad q^{2(c_i\varepsilon(i,l)+c_j\varepsilon(j,l)+c_k\varepsilon(k,l))} = 1$$

where $H^{(l)}$ is the central element $q^{-1}X_1^{(l)}X_2^{(l)}X_3^{(l)}$ of $\mathcal{T}_{T_l}^q$ and $c_i\,\varepsilon(i,l)+c_j\,\varepsilon(j,l)+c_k\,\varepsilon(k,l) = 0$ holds because $c$ is a cycle. These relations tell us that the invariants of $\rho_l$ and $\overline{\rho}_l$ are the same, and so that the representations are isomorphic. This concludes the proof of the first part of the assertion.

The second part is obvious because, if $c$ is not equal to zero, then the representation $c \cdot (\rho_1 \otimes \cdots \otimes \rho_m)$ is different from $\rho_1 \otimes \cdots \otimes \rho_m$. $\square$

The transition constants between two representatives of a local representation $\rho$ clearly depend on the choice of the orientation on $\lambda$, but the described action of $H_1(S; \mathbb{Z}_N)$ does not, it depends only on the orientation of $S$. The point is that, by changing the orientation of an edge $\lambda_i$, we change firstly the coefficient $c_i$ of $\lambda_i^*$ in $-c_i$, but we swap also the left with the right, so the resulting modification on the representations is the same. Therefore the action is intrinsic and does not depend on the choices we made.

*Remark* 1.19. If two elements $\overline{\rho}_1 \otimes \cdots \otimes \overline{\rho}_m$ and $\rho_1 \otimes \cdots \otimes \rho_m$ of $\mathscr{F}_{S_0}(\rho)$ are isomorphic, then there exists a linear isomorphism that leads from one to the other. We can give a quite explicit description of this application, which is unique up to scalar multiplication by virtue of Proposition 1.14.

Recall that, as seen in Remark 1.5, an irreducible representation $\rho_l \colon \mathcal{T}_{T_l}^q \to \mathrm{End}(V_l)$ of the triangle $T_l$ admits a basis $\mathcal{B}$ such that, if $L \colon V_l \to \mathbb{C}^N$ is the coordinate isomorphism induced by $\mathcal{B}$, we have

$$L \circ \rho_l(X_1^{(l)}) \circ L^{-1} = y_1^{(l)} B_1$$
$$L \circ \rho_l(X_2^{(l)}) \circ L^{-1} = y_2^{(l)} B_2$$
$$L \circ \rho_l(X_3^{(l)}) \circ L^{-1} = y_3^{(l)} B_3$$

where the $B_i$ are defined as

$$B_1 = \begin{pmatrix} 1 & & & \\ & q^2 & & \\ & & \ddots & \\ & & & q^{2(N-1)} \end{pmatrix}$$

$$B_2 = \begin{pmatrix} 0 & \cdots & 0 & 1 \\ & & & 0 \\ & I_{N-1} & & \vdots \\ & & & 0 \end{pmatrix}$$

$$B_3 = \begin{pmatrix} 0 & q^{1-2(2-1)} & & \\ \vdots & & \ddots & \\ 0 & & & q^{1-2(N-1)} \\ q & 0 & \cdots & 0 \end{pmatrix}$$

$y_i^{(l)}$ is a $N$-th root of $x_i^{(l)}$ for every $i = 1, 2, 3$ and $y_1^{(l)} y_2^{(l)} y_3^{(l)} = h^{(l)}$, the central load of $\rho_l$. Moreover, it is immediate to verify that the conjugation homomorphisms $C \mapsto A \circ C \circ A^{-1}$, with $A = B_1, B_2, B_3$, applied to $L \circ \rho \circ L^{-1}$ respectively change the $y_i^{(l)}$ as follows

$$y_1^{(l)} \longmapsto y_1^{(l)}$$
$$y_2^{(l)} \longmapsto q^2 y_2^{(l)}$$
$$y_3^{(l)} \longmapsto q^{-2} y_3^{(l)}$$

$$y_1^{(l)} \longmapsto q^{-2} y_1^{(l)}$$
$$y_2^{(l)} \longmapsto y_2^{(l)}$$
$$y_3^{(l)} \longmapsto q^2 y_3^{(l)}$$

$$y_1^{(l)} \longmapsto q^2 y_1^{(l)}$$
$$y_2^{(l)} \longmapsto q^{-2} y_2^{(l)}$$
$$y_3^{(l)} \longmapsto y_3^{(l)}$$

Then, defining $M_i := L^{-1} \circ B_i \circ L$ for $i = 1, 2, 3$, we have construct automorphisms $M_i$ of $V_l$ such that

$$M_1 \circ \rho_l(X_1^{(l)}) \circ M_1^{-1} = \rho_l(X_1^{(l)})$$
$$M_1 \circ \rho_l(X_2^{(l)}) \circ M_1^{-1} = q^2 \rho_l(X_2^{(l)})$$
$$M_1 \circ \rho_l(X_3^{(l)}) \circ M_1^{-1} = q^{-2} \rho_l(X_3^{(l)})$$

$$M_2 \circ \rho_l(X_1^{(l)}) \circ M_2^{-1} = q^{-2} \rho_l(X_1^{(l)})$$
$$M_2 \circ \rho_l(X_2^{(l)}) \circ M_2^{-1} = \rho_l(X_2^{(l)})$$
$$M_2 \circ \rho_l(X_3^{(l)}) \circ M_2^{-1} = q^2 \rho_l(X_3^{(l)})$$

$$M_3 \circ \rho_l(X_1^{(l)}) \circ M_3^{-1} = q^2 \rho_l(X_1^{(l)})$$
$$M_3 \circ \rho_l(X_2^{(l)}) \circ M_3^{-1} = q^{-2} \rho_l(X_2^{(l)})$$
$$M_3 \circ \rho_l(X_3^{(l)}) \circ M_3^{-1} = \rho_l(X_3^{(l)})$$

The isomorphisms $M_i$ are unique up to scalar multiplication, because of the irreducibility of the considered representations. Moreover, by means of compositions of these applications, we can obtain every change of parameters of the form

$$y_1^{(l)} \longmapsto q^{2k_1} y_1^{(l)}$$
$$y_2^{(l)} \longmapsto q^{2k_2} y_2^{(l)}$$
$$y_3^{(l)} \longmapsto q^{2k_3} y_3^{(l)}$$

for every $k_1, k_2, k_3 \in \mathbb{Z}_N$ such that $k_1 + k_2 + k_3 = 0$.

We observed that, given $\rho_1 \otimes \cdots \otimes \rho_m$ and $\overline{\rho}_1 \otimes \cdots \otimes \overline{\rho}_m$ two isomorphic representatives of $\rho_\lambda$, then there exists $c \in H_1(S; \mathbb{Z}_N)$ such that they are related as follows

$$\overline{\rho}_l(X_1^{(l)}) = q^{2c_i \varepsilon(i,l)} \rho_l(X_1^{(l)})$$
$$\overline{\rho}_l(X_2^{(l)}) = q^{2c_j \varepsilon(j,l)} \rho_l(X_2^{(l)})$$
$$\overline{\rho}_l(X_3^{(l)}) = q^{2c_k \varepsilon(k,l)} \rho_l(X_3^{(l)})$$

for every triangle $T_l$ of the ideal triangulation $\lambda \in \Lambda(S)$ (if the edges of $T_l$ are distinct, otherwise see relation 7). Fixed $l = 1, \ldots, m$, there exists a linear isomorphism $M^{(l)} \colon V_l \to V_l$ such that

$$M^{(l)} \circ \rho_l(X) \circ (M^{(l)})^{-1} = \overline{\rho}_l(X)$$

for every $X \in \mathcal{T}_\lambda^q$, and this map can be expressed, up to scalar multiplication, as composition of the elementary applications $M_i$ described above, because the elements $k_1 = c_i \varepsilon(i,l)$, $k_2 = c_j \varepsilon(j,l)$ and $k_3 = c_k \varepsilon(k,l)$ verify $k_1 + k_2 + k_3 = 0$, by virtue of the relation 8.

What just noticed shows that, for each $\rho_1 \otimes \cdots \otimes \rho_m$ and $\overline{\rho}_1 \otimes \cdots \otimes \overline{\rho}_m$ isomorphic representatives of a local representation $\rho$, the linear isomorphism conducing from one to the other is a tensor split isomorphism $M^{(1)} \otimes \cdots \otimes M^{(m)}$, in which every $M^{(l)}$ can be written as composition of the applications $M_1^{\pm 1}, M_2^{\pm 1}, M_3^{\pm 1}$ described above.

### 1.3 The quantum Teichmüller space

Given $\lambda \in \Lambda(S)$ an ideal triangulation, endowed with an indexing $\lambda_1, \ldots, \lambda_n$ of the edges, we can modify $\lambda$ in the following ways:

- for every permutation $\tau \in \mathfrak{S}_n$, we define $\lambda' = \tau(\lambda)$ the triangulation with the same 1-cells of $\lambda$, but with the ordering $\lambda'_i := \lambda_{\tau(i)}$;

- let $\lambda_i$ be an edge adjacent to two distinct triangles, composing a square $Q$. Then we denote by $\Delta_i(\lambda)$ the triangulation obtained from $\lambda$ by replacing the diagonal $\lambda_i$ of $Q$ with the other diagonal $\lambda'_i$. By definition, we set $\Delta_i(\lambda) = \lambda$ when the two sides of $\lambda_i$ belong to the same triangle.

These operations verify the following relations:

COMPOSITION RELATION: for every $\alpha, \beta$ in $\mathfrak{S}_n$ we have $\alpha(\beta(\lambda)) = (\beta \circ \alpha)(\lambda)$;

REFLEXIVITY RELATION: $(\Delta_i)^2 = id$;

RE-INDEXING RELATION: $\Delta_i \circ \alpha = \alpha \circ \Delta_{\alpha(i)}$ for every $\alpha \in \mathfrak{S}_n$;

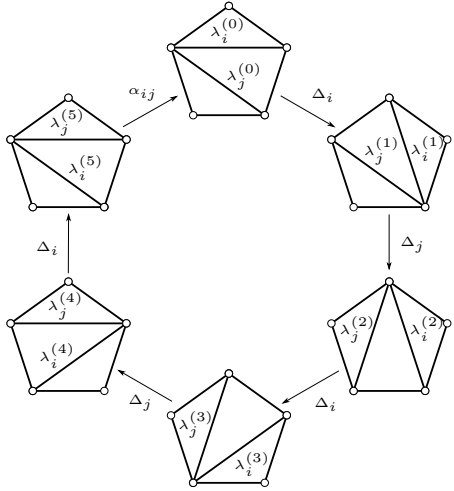

Figure 2: The Pentagon relation

DISTANT COMMUTATIVITY RELATION: if $\lambda_i$ and $\lambda_j$ do not belong to a common triangle of $\lambda$, then $(\Delta_i \circ \Delta_j)(\lambda) = (\Delta_j \circ \Delta_i)(\lambda)$;

PENTAGON RELATION: if three triangles of a triangulation $\lambda$ compose a pentagon with diagonals $\lambda_i$ and $\lambda_j$ and if we denote by $\alpha_{ij} \in \mathfrak{S}_n$ the $(ij)$ transposition, then we have

$$(\Delta_i \circ \Delta_j \circ \Delta_i \circ \Delta_j \circ \Delta_i)(\lambda) = \alpha_{ij}(\lambda)$$

The following results are due to Penner and the proofs can be found in [Pen87]:

**Theorem 1.20.** Given two ideal triangulations $\lambda, \lambda' \in \Lambda(S)$, there exists a finite sequence of ideal triangulations $\lambda = \lambda^{(0)}, \ldots, \lambda^{(k-1)}, \lambda^{(k)} = \lambda'$ such that, for every $j = 0, \ldots, k-1$ the triangulation $\lambda^{(j+1)}$ is obtained from $\lambda^{(j)}$ by a diagonal exchange or by a re-indexing of its edges.

**Theorem 1.21.** Given two ideal triangulations $\lambda, \lambda' \in \Lambda(S)$ and given two sequences $\lambda = \lambda^{(0)}, \ldots, \lambda^{(k-1)}, \lambda^{(k)} = \lambda'$ and $\lambda = \lambda'^{(0)}, \ldots, \lambda'^{(h-1)}, \lambda'^{(h)} = \lambda'$ of diagonal exchanges and re-indexing connecting $\lambda$ and $\lambda'$, then we can obtain the second sequence from the first by the applications of a finite number of the following moves or their inverses:

- using the Composition relation, replace

$$\ldots, \alpha(\lambda^{(l)}), \beta(\alpha(\lambda^{(l)})), \ldots$$

with

$$\ldots, (\alpha \circ \beta)(\lambda^{(l)}), \ldots$$

- using the Reflexivity relation, replace

$$\ldots, \lambda^{(l)}, \ldots$$

with
$$\ldots, \lambda^{(l)}, \Delta_i(\lambda^{(l)}), \lambda^{(l)}, \ldots$$

- using the Re-indexing relation, replace
$$\ldots, \lambda^{(l)}, \alpha(\lambda^{(l)}), \Delta_i(\alpha(\lambda^{(l)})), \ldots$$
with
$$\ldots, \lambda^{(l)}, \Delta_{\alpha(i)}(\lambda^{(l)}), (\alpha \circ \Delta_{\alpha(i)})(\lambda^{(l)}), \ldots$$

- using the Distant Commutativity relation, replace
$$\ldots, \lambda^{(l)}, \ldots$$
with
$$\ldots, \lambda^{(l)}, \Delta_i(\lambda^{(l)}), (\Delta_j \circ \Delta_i)(\lambda^{(l)}), \Delta_j(\lambda^{(l)}), \lambda^{(l)} \ldots$$
where $(\lambda^{(l)})_i$ and $(\lambda^{(l)})_j$ are two edges of $\lambda^{(l)}$ which do not lie in a common triangle;

- using the Pentagon relation, replace
$$\ldots, \lambda^{(l)}, \ldots$$
with
$$\ldots, \lambda^{(l)}, \Delta_i(\lambda^{(l)}), \ldots, (\Delta_j \circ \Delta_i \circ \Delta_j \circ \Delta_i)(\lambda^{(l)}), \alpha_{ij}(\lambda^{(l)}), \lambda^{(l)}, \ldots$$
where $(\lambda^{(l)})_i$ and $(\lambda^{(l)})_j$ are two diagonal of a pentagon in $\lambda^{(l)}$.

It can be observed that the Chekhov-Fock algebra is a bilateral Noetherian integral domain and so, by virtue of [Coh95, Proposition 1.3.6], it is a Øre integral domain. Therefore, we can construct $\widehat{\mathcal{T}}_\lambda^q$, the classical right quotient ring of $\mathcal{T}_\lambda^q$ (see [BL07] for the formal definition). The following Theorem, without the points related to the Fusion and Disjoint union properties, has been firstly proved in [Liu09] with a case-by-case discussion. In [BBL07] the authors have suggested a slightly simpler approach including the Fusion relation in the analysis, for the sake of completeness we have fully developed this suggestion in Appendix B.

**Theorem 1.22.** Let $S$ be a surface. Then there exists a unique collection $(\Phi_{\lambda\lambda'}^q)_{\lambda,\lambda'}$, where $\lambda$ and $\lambda'$ are varying in the set of all ideal triangulations of $S$ and $\Phi_{\lambda\lambda'}^q \colon \widehat{\mathcal{T}}_{\lambda'}^q \to \widehat{\mathcal{T}}_\lambda^q$ are algebra isomorphisms, which satisfies the following properties:

COMPOSITION RELATION: for every $\lambda, \lambda', \lambda'' \in \Lambda(S)$ ideal triangulations, we have
$$\Phi_{\lambda\lambda''}^q = \Phi_{\lambda\lambda'}^q \circ \Phi_{\lambda'\lambda''}^q$$

RE-INDEXING: if $\lambda' = \alpha(\lambda)$, then $\Phi_{\lambda\lambda'}^q(X_i) = X_{\alpha(i)}$;

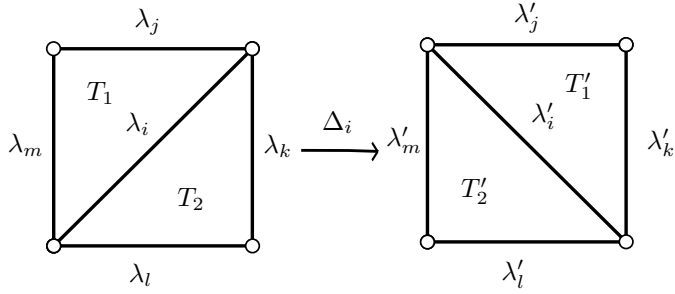

Figure 3: The ideal triangulation $\lambda, \lambda'$ of $Q$

NATURALITY: let $\varphi \colon S \to R$ be a diffeomorphism mapping the triangulations $\lambda, \lambda' \in \Lambda(S)$ in $\mu, \mu' \in \Lambda(R)$ respectively. It induces an isomorphism $\hat{\varphi}^q_{\mu\lambda} \colon \widehat{\mathcal{T}}^q_\lambda \to \widehat{\mathcal{T}}^q_\mu$, defined by sending each generator $X_i$ of $\mathcal{T}^q_\lambda$, associated with the edge $\lambda_i$ of $\lambda$, in the generator $Y_i$ of $\mathcal{T}^q_\mu$, associated with the edge $\mu_i = \varphi(\lambda_i)$ of $\mu$ and extended to the quotient rings. Analogously we define $\hat{\varphi}_{\mu'\lambda'}$. Then, the isomorphisms $\Phi^q_{\mu\mu'}$ and $\Phi^q_{\lambda\lambda'}$ verify the following relation

$$\Phi^q_{\mu\mu'} \circ \hat{\varphi}^q_{\mu'\lambda'} = \hat{\varphi}^q_{\mu\lambda} \circ \Phi^q_{\lambda\lambda'}$$

DISJOINT UNION: let $S$ be the disjoint union of $S_1$ and $S_2$, let $\lambda_1, \lambda'_1$ be two triangulations of $S_1$ and $\lambda_2, \lambda'_2$ two triangulations of $S_2$. If we denote by $\lambda := \lambda_1 \sqcup \lambda_2$ and $\lambda' = \lambda'_1 \sqcup \lambda'_2$ the induced triangulations on $S_1 \sqcup S_2$, then the isomorphism $\Phi^q_{\lambda\lambda'}$ is the extension to $\widehat{\mathcal{T}}^q_{\lambda'}$ of the algebra homomorphism

$$\mathcal{T}^q_{\lambda'} = \mathcal{T}^q_{\lambda'_1} \otimes \mathcal{T}^q_{\lambda'_2} \longrightarrow \widehat{\mathcal{T}}^q_{\lambda'_1} \otimes \widehat{\mathcal{T}}^q_{\lambda'_2} \xrightarrow{\Phi^q_{\lambda_1\lambda'_1} \otimes \Phi^q_{\lambda_2\lambda'_2}} \widehat{\mathcal{T}}^q_{\lambda_1} \otimes \widehat{\mathcal{T}}^q_{\lambda_2} \longrightarrow \widehat{\mathcal{T}}^q_\lambda$$

where the first and the third maps are the natural inclusions.

FUSION: if $S$ is obtained by fusing another surface $R$ and if $\lambda, \lambda' \in \Lambda(S)$ are obtained by fusing respectively $\mu, \mu' \in \Lambda(R)$, then

$$\hat{\iota}_{\mu\lambda} \circ \Phi^q_{\lambda\lambda'} = \Phi^q_{\mu\mu'} \circ \hat{\iota}_{\mu'\lambda'}$$

where $\hat{\iota}_{\mu\lambda}$ and $\hat{\iota}_{\mu'\lambda'}$ are the inclusions induced on the quotient rings by the maps $\iota_{\mu\lambda}$ and $\iota_{\mu'\lambda'}$ defined in 1.2.2;

DIAGONAL EXCHANGE: let $S = Q$ be the ideal square and let $\lambda, \lambda' \in \Lambda(Q)$ be the two possible ideal triangulations of $Q$, with edges labelled as in Figure 3. Then

$$\Phi^q_{\lambda\lambda'}(X'_i) = X_i^{-1}$$
$$\Phi^q_{\lambda\lambda'}(X'_j) = (1 + qX_i)X_j$$
$$\Phi^q_{\lambda\lambda'}(X'_k) = (1 + qX_i^{-1})^{-1}X_k$$
$$\Phi^q_{\lambda\lambda'}(X'_l) = (1 + qX_i)X_l$$
$$\Phi^q_{\lambda\lambda'}(X'_m) = (1 + qX_i^{-1})^{-1}X_m$$

The *quantum Teichmüller space* $\mathcal{T}_S^q$ is defined as the quotient

$$\bigsqcup_{\lambda \in \Lambda(S)} \widehat{\mathcal{T}}_\lambda^q \Big/ \sim$$

where $\sim$ is an equivalence relation that identifies two elements $X \in \widehat{\mathcal{T}}_\lambda^q$ and $X' \in \widehat{\mathcal{T}}_{\lambda'}^q$ if and only if $\Phi_{\lambda\lambda'}^q(X) = X'$. By virtue of the Composition relation, it is clear that this is an equivalence relation. We have natural bijections $i_\lambda \colon \widehat{\mathcal{T}}_\lambda^q \to \mathcal{T}_S^q$, which satisfy $i_\lambda^{-1} \circ i_{\lambda'} = \Phi_{\lambda\lambda'}^q$, for every $\lambda, \lambda' \in \Lambda(S)$. As a consequence, the set $\mathcal{T}_S^q$ can be naturally endowed with an algebra structure that make the bijections $i_\lambda$ algebra isomorphisms. Therefore, the maps $\Phi_{\lambda\lambda'}^q$ can be seen as coordinates changes, determined by the ideal triangulations $\lambda, \lambda'$, of a more intrinsic object $\mathcal{T}_S^q$.

## 1.4 Representations of the quantum Teichmüller space

Given $S$ a surface and $\lambda$ a certain triangulation of it, recall that a local representation of the Chekhov-Fock algebra $\mathcal{T}_\lambda^q$ is an equivalence class $[\rho_1, \ldots, \rho_m]$, where $\rho_j \colon \mathcal{T}_\lambda^q \to \mathrm{End}(V_{\lambda,j})$ is an irreducible representation of $\mathcal{T}_{T_j}^q$ for each triangle $T_j$ in $\lambda$. We have seen that a local representation $[\rho_1, \ldots, \rho_m]$ induces a representation of $\mathcal{T}_\lambda^q$ in the ordinary sense, by defining

$$\rho := (\rho_1 \otimes \cdots \otimes \rho_m) \circ \hat{\imath}_\lambda \colon \mathcal{T}_\lambda^q \longrightarrow \mathrm{End}(V_{\lambda,1} \otimes \cdots \otimes V_{\lambda,m}) = \mathrm{End}(V_\lambda)$$

Hereafter, with abuse, we will denote a local representation $[\rho_1, \ldots, \rho_m]$ by the only representation $\rho \colon \mathcal{T}_\lambda^q \to \mathrm{End}(V_\lambda)$.

**Definition 1.23.** Given $\lambda, \lambda' \in \Lambda(S)$ two ideal triangulations of $S$ and two representations in the same finite-dimensional vector space

$$\rho_\lambda \colon \mathcal{T}_\lambda^q \to \mathrm{End}(V) \qquad\qquad \rho_{\lambda'} \colon \mathcal{T}_{\lambda'}^q \to \mathrm{End}(V)$$

we say that $\rho_{\lambda'}$ is *compatible* with $\rho_\lambda$ and we write $\rho_{\lambda'} = \rho_\lambda \circ \Phi_{\lambda\lambda'}^q$, if, for every generator $X_i' \in \mathcal{T}_{\lambda'}^q$, the element $\Phi_{\lambda\lambda'}^q(X_i')$ can be written as $P_i Q_i^{-1} \in \widehat{\mathcal{T}}_\lambda^q$, with $P_i, Q_i \in \mathcal{T}_\lambda^q$, in such a way that $\rho_\lambda(Q_i)$ is invertible and $\rho_{\lambda'}(X_i') = \rho_\lambda(P_i)\rho_\lambda(Q_i)^{-1}$.

Observe that, by considering $\rho_{\lambda'}((X_i')^{-1})$, the element $\rho_\lambda(P_i)$ has to be invertible too.

**Lemma 1.24.** Let $\lambda^{(1)}, \lambda^{(2)}, \ldots, \lambda^{(k)}$ be a sequence of triangulations of $S$, in which $\lambda^{(i+1)}$ is obtained by re-indexing or diagonal exchange from $\lambda^{(i)}$ for every $i = 1, \ldots, k-1$ and let $\rho_i := \rho_{\lambda^{(i)}} \colon \mathcal{T}_{\lambda^{(i)}}^q \to \mathrm{End}(V)$ be a finite-dimensional representation for every $i = 1, \ldots, k$. If $\rho_i$ is compatible with $\rho_{i+1}$ for every $i = 1, \ldots, k-1$, then $\rho_1 = \rho_k \circ \Phi_{\lambda^{(k)}\lambda^{(1)}}^q$ and $\rho_k = \rho_1 \circ \Phi_{\lambda^{(1)}\lambda^{(k)}}^q$.

Consequently, the compatibility relation is symmetric, transitive and obviously reflexive.

*Proof.* See [BL07, Lemma 25] □

**Definition 1.25.** A *local representation* of the quantum Teichmüller space is a collection $\rho = \{\rho_\lambda \colon \mathcal{T}_\lambda^q \to \mathrm{End}(V_\lambda)\}_{\lambda \in \Lambda(S)}$, where

- for every $\lambda \in \Lambda(S)$ the map $\rho_\lambda\colon \mathcal{T}_\lambda^q \to \mathrm{End}(V_\lambda)$ is a local representation of $\mathcal{T}_\lambda^q$;

- for every $\lambda, \lambda' \in \Lambda(S)$ there exists a linear isomorphism $L_{\lambda\lambda'}\colon V_{\lambda'} \to V_\lambda$ such that the representation $L_{\lambda\lambda'} \circ \rho_{\lambda'}(\cdot) \circ (L_{\lambda\lambda'})^{-1}$ is compatible with $\rho_\lambda$.

**Definition 1.26.** Two local representations $\rho = \{\rho_\lambda\colon \mathcal{T}_\lambda^q \to \mathrm{End}(V_\lambda)\}_{\lambda \in \Lambda(S)}$ and $\rho' = \{\rho'_\lambda\colon \mathcal{T}_\lambda^q \to \mathrm{End}(V'_\lambda)\}_{\lambda \in \Lambda(S)}$ are said to be *isomorphic* if, for every $\lambda \in \Lambda(S)$ the representations $\rho_\lambda$ and $\rho'_\lambda$ are isomorphic.

### 1.4.1 Fusion of representations

**Definition 1.27.** If $S$ is obtained by fusion of $R$ and $\lambda \in \Lambda(S)$ is obtained by fusion of $\mu \in \Lambda(R)$, then every local representation $\rho_\mu$ of $\mathcal{T}_\mu^q$ leads to a local representation of $\mathcal{T}_\lambda^q$. More precisely, $\rho_\mu$ determines a representation $\rho_\lambda$ on $\mathcal{T}_\lambda^q$ defined by the following relation

$$\mathscr{F}_{S_0}(\rho_\mu) \subseteq \mathscr{F}_{S_0}(\rho_\lambda)$$

We will briefly say in this case that $\rho_\mu$ *represents* $\rho_\lambda$.

*Remark* 1.28. Note that requiring that $\rho_\mu$ represents $\rho_\lambda$ is stronger than saying that $\rho_\mu \circ \iota_{\mu\lambda} = \rho_\lambda$. Indeed, take $R$ equal to $S_0$, the surface obtained by splitting $S$ along every edge of $\lambda$, and $\rho_{\lambda_0}, \rho'_{\lambda_0}$ two local representations of $S_0$. If there are identified couples of edges that belong to the same triangle, the fact that $\rho_{\lambda_0}$ and $\rho'_{\lambda_0}$ can be fused to $\rho_\lambda$ does not furnish us sufficient conditions to show that they are equivalent (See the definition of local representation).

**Definition 1.29.** Let $R$ and $S$ be surfaces, with $S$ obtained by fusion from $R$. Given $\eta = \{\eta_\mu\colon \mathcal{T}_\mu^q \to \mathrm{End}(W_\mu)\}_{\mu \in \Lambda(R)}$ a local representation of $\mathcal{T}_R^q$ and $\rho = \{\rho_\lambda\colon \mathcal{T}_\lambda^q \to \mathrm{End}(V_\lambda)\}_{\lambda \in \Lambda(S)}$ a local representation of $\mathcal{T}_S^q$, $\rho$ is said to be obtained *by fusion* from $\eta$ if $\eta_\mu$ represents $\rho_\lambda$ for every ideal triangulation $\mu \in \Lambda(R)$, where $\lambda$ denotes the ideal triangulation on $S$ obtained by fusion from $\mu$.

*Remark* 1.30. There exist couples of surfaces $R, S$ with $S$ obtained by fusion from $R$, such that there exists a local representation $\eta$ of $\mathcal{T}_R^q$ not related by fusion to any local representation of $\mathcal{T}_S^q$.

For example, take $R = T_1 \sqcup T_2$ and $S = Q$, the square obtained by identifying a certain couple of edges in $T_1$ and $T_2$. $R$ admits only one ideal triangulation $\mu_0$, so a local representation of $\mathcal{T}_R^q$ is just a local representation of $\mathcal{T}_{\mu_0}^q$. Now choose a local representation $\eta_{\mu_0}$ such that its fusion $\rho_\lambda$ on $\lambda$, the induced triangulation on $Q$, has $-1$ as invariant of the diagonal in $\lambda$ of $Q$. Such a $\eta_{\mu_0}$ can be clearly constructed. Now it is evident that $\rho_\lambda$ can not be extended to a whole representation $\mathcal{T}_S^q$, because $\rho_\lambda \circ \Phi_{\lambda\lambda'}^q$ does not make sense (we are denoting by $\lambda'$ the triangulation on $Q$ obtained by diagonal exchange from $\lambda$), see Theorem 5.1. The point is that $\eta^1$, the non-quantum shadow of $\eta$, can potentially lead to a collection of non-quantum shadows that can not be extended to a non-quantum representation of $\mathcal{T}_S^1$.

### 1.4.2 The problem in [BBL07, Theorem 20]

An important consequence of Proposition 1.18 concerns the definition of the intertwining operators exposed in [BBL07]. Recall the following assertion:

**Theorem** ([BBL07, Theorem 20]). For every surface $S$ there exists a unique family of intertwining operators $\widehat{L}^{\rho\rho'}_{\lambda\lambda'}$, indexed by couples of isomorphic local representations $\rho, \rho'$ of $\mathcal{T}^q_S$ and by couples of ideal triangulation $\lambda, \lambda' \in \Lambda(S)$, individually defined up to scalar multiplication, such that:

COMPOSITION RELATION: for every $\lambda, \lambda', \lambda'' \in \Lambda(S)$ and for every triple of isomorphic local representations $\rho, \rho', \rho''$, we have $\widehat{L}^{\rho\rho''}_{\lambda\lambda''} \doteq \widehat{L}^{\rho\rho'}_{\lambda\lambda'} \circ \widehat{L}^{\rho'\rho''}_{\lambda'\lambda''}$;

FUSION RELATION: let $S$ be a surface obtained from a certain surface $R$ by fusion, and let $\lambda, \lambda'$ be two triangulations of $S$ obtained by fusion of two triangulations $\mu, \mu'$ of $R$. If $\eta, \eta'$ are two isomorphic local representations of $\mathcal{T}^q_R$ $\rho, \rho'$ are local representations of $\mathcal{T}^q_S$ obtained by fusion respectively from $\eta$ and $\eta'$, then we have $\widehat{L}^{\mu\mu'}_{\eta\eta'} \doteq \widehat{L}^{\rho\rho'}_{\lambda\lambda'}$.

In what follows we want to describe how the facts observed in Theorem 2.1 show a problem in the definition of the intertwining operators $\widehat{L}^{\rho\rho'}_{\lambda\lambda'}$ of [BBL07, Theorem 20], in particular in the case in which $\lambda = \lambda'$.

Fix $\rho = \{\rho_\lambda \colon \mathcal{T}^q_\lambda \to \mathrm{End}(V_\lambda)\}_{\lambda \in \Lambda(S)}$ and $\rho' = \{\rho'_\lambda \colon \mathcal{T}^q_\lambda \to \mathrm{End}(V_\lambda)\}_{\lambda \in \Lambda(S)}$ two isomorphic local representations of $\mathcal{T}^q_S$, where $S$ is a certain surface with $H_1(S; \mathbb{Z}_N)$ non-trivial and fix $\lambda \in \Lambda(S)$ an ideal triangulation of $S$. Since $\rho_\lambda$ and $\rho'_\lambda$ are isomorphic, we can choose representatives $\zeta = \rho_1 \otimes \cdots \otimes \rho_m$ and $\zeta' = \rho'_1 \otimes \cdots \otimes \rho'_m$ respectively of $\rho_\lambda$ and $\rho'_\lambda$ such that $\rho_j \colon \mathcal{T}^q_{T_j} \to \mathrm{End}(V_j)$ is individually isomorphic to $\rho'_j \colon \mathcal{T}^q_{T_j} \to \mathrm{End}(V'_j)$ by a linear isomorphism $L_j \colon V_j \to V'_j$. Denoting as before with $S_0$ the surface obtained by splitting $S$ along all its edges and with $\lambda_0$ its triangulation, we have that $\zeta$ and $\zeta'$ are two local representations of $\mathcal{T}^q_{\lambda_0}$ and by construction they are isomorphic by the linear transformation

$$L_1 \otimes \cdots \otimes L_m \colon V'_\lambda \longrightarrow V_\lambda$$

Moreover, this application is unique, up to scalar multiplication, because $\zeta$ and $\zeta'$ are irreducible. $S_0$ is a disjoint union of triangles, then it admits a unique ideal triangulation $\lambda_0$. This means that $\zeta$ and $\zeta'$ can be thought as local representation of the whole $\mathcal{T}^q_{S_0}$ and that $\widehat{L}^{\zeta\zeta'}_{\lambda_0\lambda_0} \doteq L_1 \otimes \cdots \otimes L_m$. Assuming that [BBL07, Theorem 20] holds, $\widehat{L}^{\rho\rho'}_{\lambda\lambda}$ must be equal to $\widehat{L}^{\zeta\zeta'}_{\lambda_0\lambda_0}$ up to scalar multiplication by virtue of the Fusion Property, so

$$\widehat{L}^{\rho\rho'}_{\lambda\lambda} \doteq L_1 \otimes \cdots \otimes L_m \tag{9}$$

Now we are going to show that different choices of representatives for $\rho_\lambda$ and $\rho'_\lambda$ produce a contradiction in relation 9. Take the same representative $\zeta$ for $\rho_\lambda$ and replace $\zeta'$ with $\overline{\zeta}' := c \cdot \zeta'$, for a certain non-trivial $c \in H_1(\Gamma; \mathbb{Z}_N) = H_1(S; \mathbb{Z}_N)$. Thanks to Proposition 1.18, the representations $\zeta'$ and $\overline{\zeta}'$ are isomorphic via an automorphism $M^{(1)} \otimes \cdots \otimes M^{(m)}$ of $V'_1 \otimes \cdots \otimes V'_1$, which is non-trivial up to scalar multiplication, because $c \neq 0$. $\overline{\zeta}'$ is a local and irreducible representation of $\mathcal{T}^q_{\lambda_0}$ and it represents $\rho'_\lambda$ on $S$, just like $\zeta'$. Moreover, the representations $\zeta$ and $\overline{\zeta}'$ are isomorphic via

$$(L_1 \circ M^{(1)}) \otimes \cdots \otimes (L_m \circ M^{(m)})$$

Applying the Fusion property as before, but on $\overline{\zeta}'$ instead of $\zeta'$, we obtain

$$\widehat{L}^{\rho\rho'}_{\lambda\lambda} \doteq \widehat{L}^{\zeta\overline{\zeta}'}_{\lambda_0\lambda_0}$$
$$\doteq (L_1 \circ M^{(1)}) \otimes \cdots \otimes (L_m \circ M^{(m)})$$

but this is in contradiction with 9, because $M^{(1)} \otimes \cdots \otimes M^{(m)}$ it is not equal to the identity up to scalar multiplication (compare with [BBL07, Lemma 22]).

The next part of the paper is devoted to the study of the consequences of Proposition 1.18 and how a result concerning intertwining operators similar to [BBL07, Theorem 20] can be recovered.

## 2 The elementary cases

The first part of our work is devoted to give the definitions of the sets $\mathscr{L}^{\rho\rho'}_{\lambda\lambda'}$ and their actions $\psi^{\rho\rho'}_{\lambda\lambda'}$ in the simplest cases, namely when $\lambda$ and $\lambda'$ differ by an elementary move. In particular the discussion will be divided in the following cases

- when $\lambda$ and $\lambda'$ are equal;
- when $\lambda$ and $\lambda'$ differ by a reindexing;
- when $\lambda$ and $\lambda'$ differ by a diagonal exchange.

In all this Section we will assume that the elements

$$\rho = \{\rho_\lambda \colon \mathcal{T}^q_\lambda \to \mathrm{End}(V_\lambda)\}_{\lambda \in \Lambda(S)} \qquad \rho' = \{\rho'_\lambda \colon \mathcal{T}^q_\lambda \to \mathrm{End}(V'_\lambda)\}_{\lambda \in \Lambda(S)}$$

are isomorphic local representations of the quantum Teichmüller space of $S$. Moreover, an intertwining operator is always thought up to scalar multiplication.

### 2.1 Same triangulation

Fixed $\lambda \in \Lambda(S)$, the maps $\rho_\lambda \colon \mathcal{T}^q_\lambda \to \mathrm{End}(V_\lambda)$ and $\rho'_\lambda \colon \mathcal{T}^q_\lambda \to \mathrm{End}(V'_\lambda)$, part of $\rho$ and $\rho'$ respectively, are two isomorphic local representations of the Chekhov-Fock algebra $\mathcal{T}^q_\lambda$. Let $S_0$ be the surface obtained by splitting $S$ along $\lambda$ and let $\lambda_0$ be its ideal triangulation. Define

$$A^{\rho\rho'}_{\lambda\lambda} := \{(\zeta, \zeta') \in \mathscr{F}_{S_0}(\rho_\lambda) \times \mathscr{F}_{S_0}(\rho'_\lambda) \mid \zeta \text{ and } \zeta' \text{ are isomorphic}\}$$

where $\zeta = \bigotimes_j \rho_j$ and $\zeta' = \bigotimes_j \rho'_j$ are thought as local (and irreducible) representations of $\mathcal{T}^q_{\lambda_0} = \bigotimes_j \mathcal{T}^q_{T_j}$. For every $(\zeta, \zeta') \in A^{\rho\rho'}_{\lambda\lambda}$ there exists a tensor-split linear isomorphism $L^{\zeta\zeta'} = L_1 \otimes \cdots \otimes L_m \colon V'_\lambda \to V_\lambda$, unique up to multiplicative constant, such that

$$L^{\zeta\zeta'} \circ \zeta'(X) \circ (L^{\zeta\zeta'})^{-1} = \zeta(X) \in \mathrm{End}(V_\lambda) \qquad\qquad \forall X \in \mathcal{T}^q_{S_0}$$

where each $L_i$ is an isomorphism between $\rho_j$ and $\rho'_j$. The uniqueness follows from the irreducibility of local representations when the surface is a disjoint union of ideal polygons (see Proposition 1.14). Now, label as $\mathscr{L}_{\lambda\lambda}^{\rho\rho'}$ the set of operators $L^{\zeta\zeta'}\colon V'_\lambda \to V_\lambda$, for varying $(\zeta,\zeta')$ in $A_{\lambda\lambda}^{\rho\rho'}$. There is an obvious surjective map $p\colon A_{\lambda\lambda}^{\rho\rho'} \to \mathscr{L}_{\lambda\lambda}^{\rho\rho'}$ that associates with a couple $(\zeta,\zeta')$ the corresponding isomorphism $L^{\zeta\zeta'}$.

Suppose that $(\zeta,\zeta')$ and $(\bar\zeta,\bar\zeta')$ go under $p$ in the same isomorphism

$$L = L^{\zeta\zeta'} = L^{\bar\zeta\bar\zeta'} = L_1 \otimes \cdots \otimes L_m : V'_\lambda \longrightarrow V_\lambda$$

In particular, $\zeta = \bigotimes_j \rho_j$ and $\bar\zeta = \bigotimes_j \bar\rho_j$ have to be equivalent, belonging both to $\mathscr{F}_{S_0}(\rho_\lambda)$. Therefore, there exist transition constants $\alpha_i \in \mathbb{C}^*$ such that $\bar\zeta \xrightarrow{\alpha_i} \zeta$. On the other hand, by hypothesis, the following hold

$$\begin{aligned} L_j \circ \rho'_j(X) \circ L_j^{-1} &= \rho_j(X) \\ L_j \circ \bar\rho'_j(X) \circ L_j^{-1} &= \bar\rho_j(X) \end{aligned} \tag{10}$$

where $\zeta' = \bigotimes_j \rho'_j$ and $\bar\zeta' = \bigotimes_j \bar\rho'_j$. By using $\bar\zeta \xrightarrow{\alpha_i} \zeta$ and the relations 10, we deduce that

- for every edge $\lambda_i$ lying in the boundary of $S$, if $T_{k_i}$ is the triangle on its side and $a_i$ is the index of the side of $T_{k_i}$ identified in $\lambda$ to $\lambda_i$, then

$$\begin{aligned} \rho'_{k_i}(X_{a_i}^{(k_i)}) &= L_{k_i}^{-1} \circ \rho_{k_i}(X_{a_i}^{(k_i)}) \circ L_{k_i} \\ &= L_{k_i}^{-1} \circ \bar\rho_{k_i}(X_{a_i}^{(k_i)}) \circ L_{k_i} \\ &= \bar\rho'_{k_i}(X_{a_i}^{(k_i)}) \end{aligned}$$

- for every internal edge $\lambda_i$ with different triangles on its sides, in the notations of Definition 1.8, we have

$$\begin{aligned} \rho'_{l_i}(X_{a_i}^{(l_i)}) &= L_{l_i}^{-1} \circ \rho_{l_i}(X_{a_i}^{(l_i)}) \circ L_{l_i} \\ &= \alpha_i\, L_{l_i}^{-1} \circ \bar\rho_{l_i}(X_{a_i}^{(l_i)}) \circ L_{l_i} \\ &= \alpha_i\, \bar\rho'_{l_i}(X_{a_i}^{(l_i)}) \end{aligned}$$

$$\begin{aligned} \rho'_{r_i}(X_{b_i}^{(r_i)}) &= L_{r_i}^{-1} \circ \rho_{r_i}(X_{b_i}^{(r_i)}) \circ L_{r_i} \\ &= \alpha_i^{-1}\, L_{r_i}^{-1} \circ \bar\rho_{r_i}(X_{b_i}^{(r_i)}) \circ L_{r_i} \\ &= \alpha_i^{-1}\, \bar\rho'_{r_i}(X_{b_i}^{(r_i)}) \end{aligned}$$

and analogously when $\lambda_i$ is an internal edge with the same triangle on its sides.

Therefore, we have shown that $p(\zeta,\zeta') = p(\bar\zeta,\bar\zeta')$ implies that the transition constants $\alpha_i$ from $\bar\zeta$ to $\zeta$ are exactly the same as those from $\bar\zeta'$ and $\zeta'$, that is

$$\begin{aligned} \bar\zeta &\xrightarrow{\alpha_i} \zeta \\ \bar\zeta' &\xrightarrow{\alpha_i} \zeta' \end{aligned}$$

We will briefly denote this phenomenon between $(\zeta,\zeta')$ and $(\bar\zeta,\bar\zeta')$ by $(\zeta,\zeta') \approx (\bar\zeta,\bar\zeta')$.

Vice versa, suppose that two couples $(\zeta, \zeta')$ and $(\bar{\zeta}, \bar{\zeta}')$ in $A^{\rho\rho'}_{\lambda\lambda}$ are in $\approx$-relation, i. e. $\bar{\zeta} \xrightarrow{\alpha_i} \zeta$ and $\bar{\zeta}' \xrightarrow{\alpha_i} \zeta'$, and label $\zeta = \bigotimes_j \rho_j$, $\zeta' = \bigotimes_j \rho'_j$, $\bar{\zeta} = \bigotimes_j \bar{\rho}_j$, $\bar{\zeta}' = \bigotimes_j \bar{\rho}'_j$. Because $(\zeta, \zeta')$ is in $A^{\rho\rho'}_{\lambda\lambda}$, there exists an isomorphism $L^{\zeta\zeta'}: V'_\lambda \to V_\lambda$ between $\zeta$ and $\zeta'$, with $L^{\zeta\zeta'} = L_1 \otimes \cdots \otimes L_m$. Then the following hold

- for every edge $\lambda_i$ lying in the boundary of $S$, if $T_{k_i}$ is the triangle on its side and $a_i$ is the index of the side of $T_{k_i}$ identified in $\lambda$ to $\lambda_i$, then

$$\begin{aligned} \bar{\rho}_{k_i}(X^{(k_i)}_{a_i}) &= \rho_{k_i}(X^{(k_i)}_{a_i}) \\ &= L_{k_i} \circ \rho'_{k_i}(X^{(k_i)}_{a_i}) \circ L^{-1}_{k_i} \\ &= L_{k_i} \circ \bar{\rho}'_{k_i}(X^{(k_i)}_{a_i}) \circ L^{-1}_{k_i} \end{aligned}$$

- for every internal edge $\lambda_i$ with different triangles on its sides, in the notations of Definition 1.8, we have

$$\begin{aligned} \bar{\rho}_{l_i}(X^{(l_i)}_{a_i}) &= \alpha_i^{-1} \rho_{l_i}(X^{(l_i)}_{a_i}) = \alpha_i^{-1} L_{l_i} \circ \rho'_{l_i}(X^{(l_i)}_{a_i}) \circ L^{-1}_{l_i} \\ &= L_{l_i} \circ \bar{\rho}'_{l_i}(X^{(l_i)}_{a_i}) \circ L^{-1}_{l_i} \\ \bar{\rho}_{r_i}(X^{(r_i)}_{b_i}) &= \alpha_i \rho_{r_i}(X^{(r_i)}_{b_i}) = \alpha_i L_{r_i} \circ \rho'_{r_i}(X^{(r_i)}_{b_i}) \circ L^{-1}_{r_i} \\ &= L_{r_i} \circ \bar{\rho}'_{r_i}(X^{(r_i)}_{b_i}) \circ L^{-1}_{r_i} \end{aligned}$$

and analogously when $\lambda_i$ has the same triangle on its sides.

Since these hold for every $i$ varying from 1 to $n$, we have shown that $L^{\zeta\zeta'}$ is an isomorphism between $\bar{\zeta}$ and $\bar{\zeta}'$, so by irreducibility $L^{\zeta\zeta'} \doteq L^{\bar{\zeta}\bar{\zeta}'}$. The equivalence relation $\approx$ on $A^{\rho\rho'}_{\lambda\lambda}$ is therefore compatible with the map $p$ and the corresponding application on the quotient $\mathscr{A}^{\rho\rho'}_{\lambda\lambda} := A^{\rho\rho'}_{\lambda\lambda} / \approx$, which we denote by

$$\widetilde{p}: \mathscr{A}^{\rho\rho'}_{\lambda\lambda} \longrightarrow \mathscr{L}^{\rho\rho'}_{\lambda\lambda}$$

is a bijection. Moreover, we can let $H_1(S, \mathbb{Z}_N)$ act on $A^{\rho\rho'}_{\lambda\lambda}$ as follows

$$c \cdot (\zeta, \zeta') := (\zeta, c \cdot \zeta') \tag{11}$$

where $c \cdot \zeta' = c \cdot (\rho'_1 \otimes \cdots \otimes \rho'_m)$ is the action defined previously, in this case on $\mathscr{F}_{S_0}(\rho'_\lambda)$. Now we want to show that the definition in 11 is compatible with the relation $\approx$ and then it leads to an action

$$\begin{aligned} \psi^{\rho\rho'}_{\lambda\lambda}: \quad H_1(S; \mathbb{Z}_N) \times \mathscr{A}^{\rho\rho'}_{\lambda\lambda} &\longrightarrow \mathscr{A}^{\rho\rho'}_{\lambda\lambda} \\ (c, [\zeta, \zeta']) &\longmapsto [\zeta, c \cdot \zeta'] \end{aligned}$$

of $H_1(S; \mathbb{Z}_N)$ on the quotient $\mathscr{A}^{\rho\rho'}_{\lambda\lambda}$ and equivalently, through $\widetilde{p}$, on $\mathscr{L}^{\rho\rho'}_{\lambda\lambda}$.

**Theorem 2.1.** The action of $H_1(S; \mathbb{Z}_N)$ on $\mathscr{A}^{\rho\rho'}_{\lambda\lambda}$, and equivalently on $\mathscr{L}^{\rho\rho'}_{\lambda\lambda}$, is well defined, transitive and free. Moreover, for every $[\zeta, \zeta'] \in \mathscr{A}^{\rho\rho'}_{\lambda\lambda}$ and for every $c \in H_1(S; \mathbb{Z}_N)$ we have

$$c \cdot [\zeta, \zeta'] = [(-c) \cdot \zeta, \zeta'] \tag{12}$$

*Proof.* In order to show the good definition, suppose that

$$\zeta = \bigotimes_i \rho_i \xrightarrow{\alpha_i} \bigotimes_i \bar{\rho}_i = \bar{\zeta}$$
$$\zeta' = \bigotimes_i \rho'_i \xrightarrow{\alpha_i} \bigotimes_i \bar{\rho}'_i = \bar{\zeta}'$$

By definition of $c\cdot$, we have

$$\zeta' \xrightarrow{q^{2c_i}} c \cdot \zeta'$$
$$\bar{\zeta}' \xrightarrow{q^{2c_i}} c \cdot \bar{\zeta}'$$

Then

$$c \cdot \zeta' \xrightarrow{q^{-2c_i}} \zeta' \xrightarrow{\alpha_i} \bar{\zeta}' \xrightarrow{q^{2c_i}} c \cdot \bar{\zeta}'$$

so $c \cdot \zeta' \xrightarrow{\alpha_i} c \cdot \bar{\zeta}'$. On the other hand, we have $\zeta \xrightarrow{\alpha_i} \bar{\zeta}$, and these facts together tell us that $c \cdot (\zeta, \zeta') := (\zeta, c \cdot \zeta') \approx (\bar{\zeta}, c \cdot \bar{\zeta}') =: c \cdot (\bar{\zeta}, \bar{\zeta}')$, as desired.

Now we will prove that the action is transitive. Let $[\zeta, \zeta'], [\bar{\zeta}, \bar{\zeta}']$ be two elements of $\mathscr{A}^{\rho\rho'}_{\lambda\lambda}$ and $(\zeta, \zeta'), (\bar{\zeta}, \bar{\zeta}')$ two representatives of them, with $\zeta = \bigotimes_j \rho_j$, $\zeta' = \bigotimes_j \rho'_j$, $\bar{\zeta} = \bigotimes_j \bar{\rho}_j$ and $\bar{\zeta}' = \bigotimes_j \bar{\rho}'_j$. Both $\zeta$ and $\bar{\zeta}$ belong to $\mathscr{F}_{S_0}(\rho_\lambda)$, then there exists a family $(\alpha_i)_i$ of transition constants such that $\zeta \xrightarrow{\alpha_i} \bar{\zeta}$. Because $(\zeta, \zeta')$ is an element of $A^{\rho\rho'}_{\lambda\lambda}$, there exist isomorphisms $L_j \colon V'_j \to V_j$ such that $L_j \circ \rho'_j \circ L_j^{-1} = \rho_j$ for every $j = 1, \ldots, m$. Now we can construct a representation $\widetilde{\zeta}' = \bigotimes_j \widetilde{\rho}'_j$ defined by

$$\widetilde{\rho}'_j(X) := L_j^{-1} \circ \bar{\rho}_j(X) \circ L_j$$

for every $j = 1, \ldots, m$. The representation $\widetilde{\zeta}'$ belongs to $\mathscr{F}_{S_0}(\rho'_\lambda)$ because by construction it can be obtained from $\zeta'$, which is an element of $\mathscr{F}_{S_0}(\rho'_\lambda)$, as

$$\zeta' \xrightarrow{\alpha_i} \widetilde{\zeta}'$$

So $(\bar{\zeta}, \widetilde{\zeta}')$ belongs to $A^{\rho\rho'}_{\lambda\lambda}$ and $\zeta', \widetilde{\zeta}'$ are related by the transition constants $(\alpha_i)_i$, just like $\zeta$ and $\bar{\zeta}$ (see the relations previously used to prove $p(\zeta, \zeta') = p(\bar{\zeta}, \bar{\zeta}') \Rightarrow (\zeta, \zeta') \approx (\bar{\zeta}, \bar{\zeta}')$). This means that the couples $(\zeta, \zeta')$ and $(\bar{\zeta}, \widetilde{\zeta}')$ are $\approx$-equivalent, i. e. $[\zeta, \zeta'] = [\bar{\zeta}, \widetilde{\zeta}']$. Moreover, both $\widetilde{\zeta}'$ and $\bar{\zeta}'$ are isomorphic to $\bar{\zeta}$ and then they are isomorphic to each other. By Proposition 1.18, there exists a unique $c \in H_1(S; \mathbb{Z}_N)$ such that $c \cdot \widetilde{\zeta}' = \bar{\zeta}'$, so

$$c \cdot [\zeta, \zeta'] = c \cdot [\bar{\zeta}, \widetilde{\zeta}'] = [\bar{\zeta}, c \cdot \widetilde{\zeta}'] = [\bar{\zeta}, \bar{\zeta}']$$

and this proves that the action is transitive.

Now suppose that there exist a $c \in H_1(S; \mathbb{Z}_N)$ and an element $[\zeta, \zeta'] \in \mathscr{A}^{\rho\rho'}_{\lambda\lambda}$ such that $[\zeta, \zeta'] = c \cdot [\zeta, \zeta']$. This means that, passing on representatives, the couples $(\zeta, \zeta')$ and $(\zeta, c \cdot \zeta')$ are $\approx$-equivalent. Because the first terms of the couples are exactly the same, they are in particular related by transition constants all equal to 1, and the same must holds for $\zeta'$ and $c \cdot \zeta'$. This means that $\zeta' = c \cdot \zeta'$ and so, by the second assertion of Proposition 1.18, we conclude.

For what concerns the equation 12, we firstly note that $[\zeta, \zeta'] = [c \cdot \zeta, c \cdot \zeta']$. Indeed, we have $\zeta \stackrel{q^{2c_i}}{\to} c \cdot \zeta$ and $\zeta' \stackrel{q^{2c_i}}{\to} c \cdot \zeta'$, which means that $(\zeta, \zeta') \approx (c \cdot \zeta, c \cdot \zeta')$. Now it is immediate to prove the relation:

$$[\zeta, c \cdot \zeta'] = [(-c) \cdot \zeta, (-c + c) \cdot \zeta'] = [(-c) \cdot \zeta, \zeta']$$

$\square$

We will denote the action of $H_1(S; \mathbb{Z}_N)$ on $\mathscr{L}_{\lambda\lambda}^{\rho\rho'}$ by $\psi_{\lambda\lambda}^{\rho\rho'}$.

## 2.2 Reindexing

In the investigation of local representations we have intentionally ignored the problems concerning the case in which the ideal triangulations $\lambda$ and $\lambda'$ differ by reindexing, i. e. $\lambda' = \gamma(\lambda)$ with $\gamma \in \mathfrak{S}_n$. We did not focus on that because all the properties of representations are basically intrinsic and does not really depend on the ordering of the edges, but only on the structure of the triangulation. Indeed, the coordinate change isomorphisms $\Phi_{\lambda\lambda'}^q$ in this case are just the maps on the fraction rings induced by the reordering applications from $\mathcal{T}_\lambda^q$ to $\mathcal{T}_{\gamma(\lambda)}^q$. Moreover, the described action of $H_1(S; \mathbb{Z}_N)$ clearly does not depend on the fixed order of the edges. We will continue to be vague on that, we want just to enunciate the fact we will use later, analogous to Theorem 2.1.

Let $\lambda, \lambda' \in \Lambda(S)$ be two ideal triangulations differing by a reindexing of the edges, with $\lambda' = \gamma(\lambda)$. Define $A_{\lambda\lambda'}^{\rho\rho'}$ as the set of couples $(\zeta_{\lambda_0}, \zeta'_{\lambda'_0})$, where $\zeta_{\lambda_0}$ is an element of $\mathscr{F}_{S_0}(\rho_\lambda)$, $\zeta'_{\lambda'_0}$ is an element of $\mathscr{F}_{S_0}(\rho'_{\lambda'})$ and $\zeta_{\lambda_0} \circ \Phi_{\lambda_0\lambda'_0}^q$ is isomorphic to $\zeta'_{\lambda'_0}$ ($\lambda$ and $\lambda'$ clearly induce the same splitted surface $S_0$, we should give details of indexing of triangulations $\lambda_0$ and $\lambda'_0$ in order give a sense to $\Phi_{\lambda_0\lambda'_0}^q$, but we omit this boring procedure). We say that $(\zeta_{\lambda_0}, \zeta'_{\lambda'_0}), (\bar\zeta_{\lambda_0}, \bar\zeta'_{\lambda'_0}) \in A_{\lambda\lambda'}^{\rho\rho'}$ are $\approx$-equivalent if

$$\zeta_{\lambda_0} \xrightarrow{\alpha_i} \bar\zeta_{\lambda_0}$$
$$\zeta'_{\lambda'_0} \xrightarrow{\beta_i} \bar\zeta'_{\lambda'_0}$$

with $\beta_i = \alpha_{\gamma(i)}$ for every $i$. As before $\mathscr{A}_{\lambda\lambda'}^{\rho\rho'}$ denotes the quotient of $A_{\lambda\lambda'}^{\rho\rho'}$ by the relation $\approx$. On $A_{\lambda\lambda'}^{\rho\rho'}$ we consider a natural map $\widetilde{p} \colon \mathscr{A}_{\lambda\lambda'}^{\rho\rho'} \to \mathrm{Hom}(V'_{\lambda'}, V_\lambda)/\doteq$, sending an element $[\zeta_{\lambda_0}, \zeta'_{\lambda'_0}]$ in the linear isomorphism between $\zeta_{\lambda_0}$ and $\zeta'_{\lambda'_0}$. The map is injective and we designate its image as $\mathscr{L}_{\lambda\lambda'}^{\rho\rho'}$. We can define on $\mathscr{A}_{\lambda\lambda'}^{\rho\rho'}$ an action $\psi_{\lambda\lambda'}^{\rho\rho'}$ of $H_1(S; \mathbb{Z}_N)$ setting $c \cdot [\zeta_{\lambda_0}, \zeta'_{\lambda'_0}] := [\zeta_{\lambda_0}, c \cdot \zeta'_{\lambda'_0}]$. In light of Theorem 2.1, it is straightforward to prove that the following holds

**Theorem 2.2.** *The action of $H_1(S; \mathbb{Z}_N)$ on $\mathscr{A}_{\lambda\lambda'}^{\rho\rho'}$, and equivalently on $\mathscr{L}_{\lambda\lambda'}^{\rho\rho'}$, is well defined, transitive and free. Moreover, for every $[\zeta_{\lambda_0}, \zeta'_{\lambda'_0}] \in \mathscr{A}_{\lambda\lambda'}^{\rho\rho'}$ and for every $c \in H_1(S; \mathbb{Z}_N)$ we have*

$$c \cdot [\zeta_{\lambda_0}, \zeta'_{\lambda'_0}] = [(-c) \cdot \zeta_{\lambda_0}, \zeta'_{\lambda'_0}]$$

## 2.3 Diagonal exchange

Assume that $\lambda$ and $\lambda'$ are ideal triangulations of $S$ differing by a diagonal exchange along $\lambda_i$. Label as $R$ the surface obtained from $S$ by splitting it along all the edges of $\lambda$ except for $\lambda_i$. $R$ is the disjoint union of an ideal square $Q$ and $m-2$ ideal triangles. In order to simplify the notation, we will assume that the triangles composing $Q$ are labelled as $T_1$ and $T_2$ and the others as $T_3, \ldots, T_m$. The triangulations $\lambda$ and $\lambda'$ induce on $R$ two ideal triangulations $\mu, \mu'$. $\mu$ is just the disjoint union of an ideal triangulation $\mu_Q$ of $Q$ and the only possible triangulation $\mu_0$ on the disjoint union of the triangles $T_j$ for $j \geq 3$. Analogously $\mu' = \mu'_Q \sqcup \mu_0$ where $\mu'_Q$ is the ideal triangulation on $Q$ obtained by diagonal exchange on $\mu_Q$. Observe that the Chekhov-Fock algebras associated with the triangulation $\mu$ and $\mu'$ on $R$ are canonically isomorphic to the tensor products

$$\mathcal{T}^q_{\mu_Q} \otimes \mathcal{T}^q_{T_3} \otimes \cdots \otimes \mathcal{T}^q_{T_m}, \quad \mathcal{T}^q_{\mu'_Q} \otimes \mathcal{T}^q_{T_3} \otimes \cdots \otimes \mathcal{T}^q_{T_m}$$

We will denote by $S_0$ the surface obtained by splitting $S$ along $\lambda$, by $S'_0$ the surface obtained by splitting $S$ along $\lambda'$ and by $\lambda_0$ and $\lambda'_0$ the respective triangulations on these surfaces.

Fix $\rho = \{\rho_\eta \colon \mathcal{T}^q_\eta \to \mathrm{End}(V_\eta)\}_{\eta \in \Lambda(S)}$ and $\rho' = \{\rho'_\eta \colon \mathcal{T}^q_\eta \to \mathrm{End}(V_\eta)\}_{\eta \in \Lambda(S)}$ two local representations of $\mathcal{T}^q_S$. We introduce the following notations

$$V_\eta = V_{\eta,1} \otimes \cdots \otimes V_{\eta,m}$$
$$V'_\eta = V'_{\eta,1} \otimes \cdots \otimes V'_{\eta,m}$$

Denote by $\mathscr{F}_{S_0}(\rho_\lambda)$ the set of local representations of $\mathcal{T}^q_{\lambda_0}$ that represent $\rho_\lambda$ on $S$ and analogously label as $\mathscr{F}_{S'_0}(\rho'_{\lambda'})$ the set of local representations of $\mathcal{T}^q_{\lambda'_0}$ that represent $\rho'_{\lambda'}$ on $S$. Given $\zeta_{\lambda_0}$ an element of $\mathscr{F}_{S_0}(\rho_\lambda)$, $\zeta_{\lambda_0}$ represents a local representation $\zeta_\mu$ of the Chekhov-Fock algebra $\mathcal{T}^q_\mu$, and in the same way a representation $\zeta'_{\lambda'_0} \in \mathscr{F}_{S'_0}(\rho'_{\lambda'})$ induces a representation $\zeta'_{\mu'}$ of $\mathcal{T}^q_{\mu'}$.

Now we define

$$A^{\rho\rho'}_{\lambda\lambda'} := \{(\zeta_{\lambda_0}, \zeta'_{\lambda'_0}) \in \mathscr{F}_{S_0}(\rho_\lambda) \times \mathscr{F}_{S'_0}(\rho'_{\lambda'}) \mid \zeta_\mu \circ \Phi^q_{\mu\mu'} \text{ is isomorphic to } \zeta'_{\mu'}\}$$

It is easy to verify that the composition $\zeta_\mu \circ \Phi^q_{\mu\mu'}$ makes sense because $\rho_\lambda \circ \Phi^q_{\lambda\lambda'}$ does (the key ingredient is that the invariant of the representation $\zeta_\mu$ associated with $\mu_i$ coincides with the one of $\rho_\lambda$ for $\lambda_i$, which is not equal to $-1$ because $\rho_\lambda \circ \Phi^q_{\lambda\lambda'}$ makes sense, being $\rho_\lambda$ part of a global representation of the quantum Teichmüller space). Given $\zeta_{\lambda_0}$ in $\mathscr{F}_{S_0}(\rho_\lambda)$, $\zeta_\mu$ is equal to the tensor product of a representation $\zeta_{\mu_Q}$ of $\mathcal{T}^q_{\mu_Q}$ and a representation $\zeta_{\mu_0}$ of $\mathcal{T}^q_{\mu_0} = \mathcal{T}^q_{T_3} \otimes \cdots \otimes \mathcal{T}^q_{T_m}$. In the same way, given $\zeta'_{\lambda'_0} \in \mathscr{F}_{S'_0}(\rho'_{\lambda'})$, $\zeta'_{\mu'}$ is the tensor product of a representation $\zeta'_{\mu'_Q}$ of $\mathcal{T}^q_{\mu'_Q}$ and a representation $\zeta'_{\mu_0}$ of $\mathcal{T}^q_{\mu_0}$. Recalling the Disjoint union property of $\Phi^q_{\lambda\lambda'}$ exposed in Theorem 1.22, the restriction of $\Phi^q_{\mu\mu'}$ on $\mathcal{T}^q_\mu = \mathcal{T}^q_{\mu_Q} \otimes \mathcal{T}^q_{\mu_0}$ coincides with $\Phi^q_{\mu_Q\mu'_Q} \otimes id$. Thus the representation $\zeta_\mu \circ \Phi^q_{\mu\mu'}$ is equal to

$$(\zeta_{\mu_Q} \circ \Phi^q_{\mu_Q\mu'_Q}) \otimes \zeta_{\mu_0} \tag{13}$$

By virtue of the irreducibility of $\zeta_\mu$ and $\zeta'_{\mu'}$ (observe that $R$ is a disjoint union of ideal polygons) there exists an isomorphism $L^{\zeta_\mu \zeta'_{\mu'}}$, unique up to scalar multiplication, such that

$$L^{\zeta_\mu \zeta'_{\mu'}} \circ \zeta'_{\mu'}(X') \circ (L^{\zeta_\mu \zeta'_{\mu'}})^{-1} = (\zeta_\mu \circ \Phi^q_{\mu\mu'})(X')$$

for every $X' \in \mathcal{T}_{\mu'}^q$. In analogy with the case $\lambda = \lambda'$, we designate as $\mathscr{L}_{\lambda\lambda'}^{\rho\rho'}$ the set of operators $L^{\zeta_\mu \zeta'_{\mu'}}$, for varying $(\zeta_\mu, \zeta'_{\mu'})$ in $A_{\lambda\lambda'}^{\rho\rho'}$. Because of relation 13, every $L^{\zeta_\mu \zeta'_{\mu'}}$ is the tensor product of an isomorphism

$$L^{\zeta_{\mu_Q} \zeta'_{\mu'_Q}} : V'_{\lambda',1} \otimes V'_{\lambda',2} \longrightarrow V_{\lambda,1} \otimes V_{\lambda,2}$$

between $\zeta_{\mu_Q} \circ \Phi^q_{\mu_Q \mu'_Q}$ and $\zeta'_{\mu'_Q}$, and of an isomorphism

$$L^{\zeta_{\mu_0} \zeta'_{\mu_0}} : V'_{\lambda',3} \otimes \cdots \otimes V'_{\lambda',m} \longrightarrow V_{\lambda,3} \otimes \cdots \otimes V_{\lambda,m}$$

between $\zeta_{\mu_0}$ and $\zeta'_{\mu_0}$, which is tensor-split.

As before, we define the map

$$p: \begin{array}{ccc} A_{\lambda\lambda'}^{\rho\rho'} & \longrightarrow & \mathscr{L}_{\lambda\lambda'}^{\rho\rho'} \\ (\zeta_{\lambda_0}, \zeta'_{\lambda'_0}) & \longmapsto & L^{\zeta_\mu \zeta'_{\mu'}} \end{array}$$

The map $p$ is tautologically surjective, we want to characterize its injective quotient. If $\zeta_{\lambda_0}$ and $\bar{\zeta}_{\lambda_0}$ belong to $\mathscr{F}_{S_0}(\rho_\lambda)$, then they both represent $\rho_\lambda$. Fixed arbitrary orientations on the edges of $\lambda$, there exist transition constants $(\alpha_j)_j$ such that $\zeta_{\lambda_0} \xrightarrow{\alpha_j} \bar{\zeta}_{\lambda_0}$. If $\zeta_\mu = \zeta_{\mu_Q} \otimes \zeta_3 \otimes \cdots \otimes \zeta_m$ and $\bar{\zeta}_{\mu_Q} \otimes \bar{\zeta}_3 \otimes \cdots \otimes \bar{\zeta}_m$ are the induced representations on $\mathcal{T}_\mu^q$, then, for every $\lambda_j$ with $j \neq i$, the following hold

- if $\lambda_j$ is on the boundary of $S$, then the two representations must coincide on the only variable in $\mathcal{T}_\mu^q$ corresponding to $\lambda_j$;

- if $\lambda_j$ is internal and it is side of two triangles $T_{l_j}$ and $T_{r_j}$, with $l_j, r_j \geq 3$, on the left and on the right respectively of $\lambda_j$, then

$$\bar{\zeta}_{l_j}(X_{a_j}^{(l_j)}) = \alpha_j \, \zeta_{l_j}(X_{a_j}^{(l_j)})$$
$$\bar{\zeta}_{r_j}(X_{b_j}^{(r_j)}) = \alpha_j^{-1} \, \zeta_{r_j}(X_{b_j}^{(r_j)})$$

where $a_j$ and $b_j$ are the indices of the sides in $T_{l_j}$ and $T_{r_j}$, respectively, identified to $\lambda_j$ in $S$ (analogously if $T_{l_j} = T_{r_j}$);

- if $\lambda_j$ is internal and it is side of a triangle $T_{k_j}$ and of the square $Q$, then

$$\bar{\zeta}_{\mu_Q}(X_{a_j}^{(Q)}) = \alpha_j^{\varepsilon(j,Q)} \, \zeta_{\mu_Q}(X_{a_j}^{(Q)})$$
$$\bar{\zeta}_{k_j}(X_{b_j}^{(k_j)}) = \alpha_j^{\varepsilon(j,k_j)} \, \zeta_{k_j}(X_{b_j}^{(k_j)})$$

where $a_j$ and $b_j$ are the indices of the sides in $Q$ and $T_{k_j}$, respectively, identified to $\lambda_j$ in $S$, $\varepsilon(j, Q)$ is equal to $+1$ if the orientation of $\lambda_j$ coincides with the boundary orientation of $Q$, $-1$ otherwise, and $\varepsilon(j, k_j)$ is equal to $+1$ if the orientation of $\lambda_j$ coincides with the boundary orientation of $T_{k_j}$, $-1$ otherwise;

and analogously in the case in which $\lambda_j$ has on both sides the square $Q$. Observe that the constant $\alpha_i$ does not appear in the discussion because we are considering the equivalence classes $\zeta_\mu, \bar{\zeta}_\mu$ of local representation of $\mathcal{T}_\mu^q$, instead of the representations $\zeta_{\lambda_0}$ and $\bar{\zeta}'_{\lambda'_0}$. We will say that the constants $\alpha_j$, for $j \neq i$, are the *transition constants* from $\zeta_\mu$ to $\zeta'_{\mu'}$, which can be thought as supported on $\lambda \setminus \lambda_i$. The same can be noted for a couple of $\zeta'_{\lambda'_0}, \bar{\zeta}'_{\lambda_0}$ in $\mathscr{F}_{S_0}(\rho'_{\lambda'})$, providing a collection of transition constants $\beta_j$ from $\zeta'_{\mu'}$ to $\bar{\zeta}'_{\mu'}$. The two collections $(\alpha_j)$ and $(\beta_j)$ can be compared in a natural way because there is a canonical correspondence between $\lambda \setminus \lambda_i$ and $\lambda' \setminus \lambda'_i$.

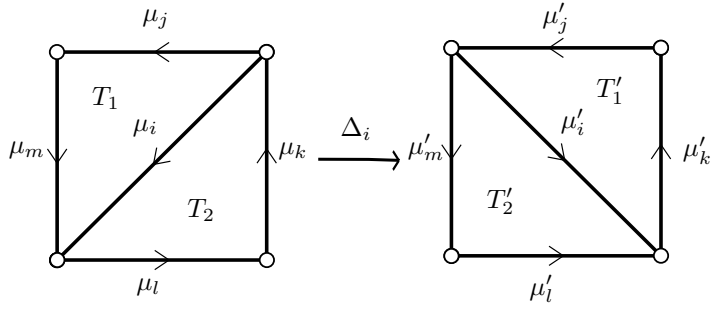

Figure 4: The ideal triangulations $\mu_Q, \mu'_Q \in \Lambda(Q)$

Given two couples $(\zeta_{\lambda_0}, \zeta'_{\lambda'_0}), (\bar{\zeta}_{\lambda_0}, \bar{\zeta}'_{\lambda'_0}) \in A_{\lambda\lambda'}^{\rho\rho'}$, they are in $\approx$-relation if the transition constants from $\zeta_\mu$ to $\bar{\zeta}_\mu$ are the same of those from $\zeta'_{\mu'}$ to $\bar{\zeta}'_{\mu'}$. Observe that the relation $\approx$ can be expressed in terms of the transition constants $\zeta_{\lambda_0} \xrightarrow{\alpha_j} \bar{\zeta}_{\lambda_0}$ and $\zeta'_{\lambda'_0} \xrightarrow{\beta_j} \bar{\zeta}'_{\lambda'_0}$ by asking $\alpha_j = \beta_j$ for every $j \neq i$, and by not requiring any restriction on $\alpha_i$ and $\beta_i$.

We want to show that, if $(\zeta_{\lambda_0}, \zeta'_{\lambda'_0}) \approx (\bar{\zeta}_{\lambda_0}, \bar{\zeta}'_{\lambda'_0})$, then $L^{\zeta_\mu \zeta'_{\mu'}} \doteq L^{\bar{\zeta}_\mu \bar{\zeta}'_{\mu'}}$. It is sufficient to prove that the following hold

$$L^{\zeta_{\mu_Q} \zeta'_{\mu'_Q}} \doteq L^{\bar{\zeta}_{\mu_Q} \bar{\zeta}'_{\mu'_Q}}$$
$$L^{\zeta_{\mu_0} \zeta'_{\mu_0}} \doteq L^{\bar{\zeta}_{\mu_0} \bar{\zeta}'_{\mu_0}}$$

The second equality can be obtained with exactly the same observations of the case $\lambda = \lambda'$ done above. Therefore, we will concentrate only on the first one, for which we must pay a little more attention because we have to manage the composition $\zeta_{\mu_Q} \circ \Phi^q_{\mu_Q \mu'_Q}$.

Assume that the edges of the square are labelled as in Figure 4 and, in order to simplify the notation, that they are oriented counter-clockwise with respect to the orientation of $Q$. Then there exist $\alpha_h \in \mathbb{C}^*$ for $h \in \{j, k, l, m\}$ such that

$$\zeta_{\mu_Q}(X_h^{(Q)}) = \alpha_h \bar{\zeta}_{\mu_Q}(X_h^{(Q)})$$

where we are denoting by $X_h^{(Q)}$ the element of the Chekhov-Fock algebra $\mathcal{T}^q_{\mu_Q}$ associated with the edge $\mu_h$. Using the fact that $\zeta_{\mu_Q}(X_i^{(Q)}) = \bar{\zeta}_{\mu_Q}(X_i^{(Q)})$ (this equality holds because both $\zeta_\mu$ and $\bar{\zeta}_\mu$ represent $\rho_\lambda$ and the edge $\mu_i$ corresponding to $\lambda_i$ is already fused in $\mu$) and the explicit formulas for $\Phi^q_{\mu_Q \mu'_Q}$, we obtain that

$$\zeta_{\mu_Q}(X_h^{(Q)}) = \alpha_h \bar{\zeta}_{\mu_Q}(X_h^{(Q)}) \Leftrightarrow (\zeta_{\mu_Q} \circ \Phi^q_{\mu_Q \mu'_Q})(Y_h^{(Q)}) = \alpha_h (\bar{\zeta}_{\mu_Q} \circ \Phi^q_{\mu_Q \mu'_Q})(Y_h^{(Q)})$$
(14)

Now the same argument of the previous case can be applied in order to conclude the desired equality, the only difference is that in this case one of the representations is on the square $Q$ instead of a triangle.

The map $p$ induces on the quotient

$$\mathscr{A}_{\lambda\lambda'}^{\rho\rho'} := A_{\lambda\lambda'}^{\rho\rho'}/_{\approx}$$

an application $\widetilde{p}$. With the same argument done in the case $\mathscr{A}_{\lambda\lambda}^{\rho\rho'}$ and using relation 14, we can prove that $\widetilde{p}$ is injective, and so bijective, since the surjectivity is obvious. Moreover, we can describe an action of $H_1(S; \mathbb{Z}_N)$ on $\mathscr{A}_{\lambda\lambda'}^{\rho\rho'}$ when $\lambda$ and $\lambda'$ differ by a diagonal exchange along $\lambda_i$, by defining

$$\begin{array}{rccc} \psi_{\lambda\lambda'}^{\rho\rho'}: & H_1(S; \mathbb{Z}_N) \times \mathscr{A}_{\lambda\lambda'}^{\rho\rho'} & \longrightarrow & \mathscr{A}_{\lambda\lambda'}^{\rho\rho'} \\ & (c, [\zeta_{\lambda_0}, \zeta'_{\lambda'_0}]) & \longmapsto & [\zeta_{\lambda_0}, c \cdot \zeta'_{\lambda'_0}] \end{array}$$

where $c \cdot \zeta'_{\lambda'_0}$ is the action of $c$ on $\zeta'_{\lambda'_0}$ in $\mathscr{F}_{S'_0}(\rho'_{\lambda'})$ as in Proposition 1.18.

**Theorem 2.3.** *The action of $H_1(S; \mathbb{Z}_N)$ on $\mathscr{A}_{\lambda\lambda'}^{\rho\rho'}$, and equivalently on $\mathscr{L}_{\lambda\lambda'}^{\rho\rho'}$, is well defined, transitive and free. Moreover, for every $[\zeta_{\lambda_0}, \zeta'_{\lambda'_0}] \in \mathscr{A}_{\lambda\lambda'}^{\rho\rho'}$ and for every $c \in H_1(S; \mathbb{Z}_N)$ we have*

$$c \cdot [\zeta_{\lambda_0}, \zeta'_{\lambda'_0}] = [(-c) \cdot \zeta_{\lambda_0}, \zeta'_{\lambda'_0}] \tag{15}$$

*Proof.* The proof will be very similar to the one of Theorem 2.1. Take two couples $(\zeta_{\lambda_0}, \zeta'_{\lambda'_0}), (\bar{\zeta}_{\lambda_0}, \bar{\zeta}'_{\lambda'_0})$ in $A_{\lambda\lambda'}^{\rho\rho'}$ that are $\approx$-equivalent. Then

$$\zeta_{\lambda_0} \xrightarrow{\alpha_j} \bar{\zeta}_{\lambda_0}$$

$$\zeta'_{\lambda'_0} \xrightarrow{\beta_j} \bar{\zeta}'_{\lambda'_0}$$

with $\alpha_j = \beta_j$ for every $j \neq i$. By definition of the action $\psi_{\lambda\lambda'}^{\rho\rho'}$ we have

$$\zeta'_{\lambda'_0} \xrightarrow{q^{2c_j}} c \cdot \zeta'_{\lambda'_0}$$

$$\bar{\zeta}'_{\lambda'_0} \xrightarrow{q^{2c_j}} c \cdot \bar{\zeta}'_{\lambda'_0}$$

Then

$$c \cdot \zeta'_{\lambda'_0} \xrightarrow{q^{-2c_j}} \zeta'_{\lambda'_0} \xrightarrow{\beta_j} \bar{\zeta}'_{\lambda'_0} \xrightarrow{q^{2c_j}} c \cdot \bar{\zeta}'_{\lambda'_0}$$

We conclude that $\zeta_{\lambda_0} \xrightarrow{\alpha_j} \bar{\zeta}_{\lambda_0}$ and $c \cdot \zeta'_{\lambda'_0} \xrightarrow{\beta_j} c \cdot \bar{\zeta}'_{\lambda'_0}$, with $\alpha_j = \beta_j$ for every $j \neq i$, hence $(\zeta_{\lambda_0}, c \cdot \zeta'_{\lambda'_0}) \approx (\bar{\zeta}_{\lambda_0}, c \cdot \bar{\zeta}'_{\lambda'_0})$, which proves the good definition of the action.

In order to prove the transitivity, fix two elements $[\zeta_{\lambda_0}, \zeta'_{\lambda'_0}], [\bar{\zeta}_{\lambda_0}, \bar{\zeta}'_{\lambda'_0}]$ in $\mathscr{A}_{\lambda\lambda'}^{\rho\rho'}$ and two respective representatives $(\zeta_{\lambda_0}, \zeta'_{\lambda'_0}), (\bar{\zeta}_{\lambda_0}, \bar{\zeta}'_{\lambda'_0})$. Denote by $(\alpha_j)_j$ the transition constants from $\zeta_{\lambda_0}$ to $\bar{\zeta}_{\lambda_0}$. The element $(\zeta_{\lambda_0}, \zeta'_{\lambda'_0})$ belongs to $A_{\lambda\lambda'}^{\rho\rho'}$, hence the induced local representations $\zeta_\mu \circ \Phi^q_{\mu\mu'}$ and $\zeta'_{\mu'}$ are isomorphic via an isomorphism

$$L^{\zeta_\mu \zeta'_{\mu'}} = L^{\zeta_{\mu_Q} \zeta'_{\mu'_Q}} \otimes L^{\zeta_{\mu_0} \zeta'_{\mu'_0}} : V'_{\lambda'} \longrightarrow V_\lambda$$

We can construct a representation $\widetilde{\zeta}'_{\lambda'_0}$ such that $(\bar{\zeta}_{\lambda_0}, \widetilde{\zeta}'_{\lambda'_0})$ belongs to $A^{\rho\rho'}_{\lambda\lambda'}$ and $(\zeta_{\lambda_0}, \zeta'_{\lambda'_0}) \approx (\bar{\zeta}_{\lambda_0}, \widetilde{\zeta}'_{\lambda'_0})$, simply by defining $\widetilde{\zeta}'_{\lambda'_0}$ as the representation verifying

$$\zeta'_{\lambda'_0} \xrightarrow{\alpha_j} \widetilde{\zeta}'_{\lambda'_0}$$

where the $\alpha_j$ are the transition constants from $\zeta_{\lambda_0}$ to $\bar{\zeta}_{\lambda_0}$ (it is not important which is the transition constant in the edge $\lambda'_i$). Because $\zeta_\mu$ and $\zeta'_{\mu'}$ verify

$$L^{\zeta_\mu \zeta'_{\mu'}} \circ \zeta'_{\mu'}(X') \circ (L^{\zeta_\mu \zeta'_{\mu'}})^{-1} = (\zeta_\mu \circ \Phi^q_{\mu\mu'})(X') \qquad \forall X' \in \mathcal{T}^q_{\mu'}$$

and because $\zeta_{\lambda_0} \xrightarrow{\alpha_j} \bar{\zeta}_{\lambda_0}$, then we have also

$$L^{\zeta_\mu \zeta'_{\mu'}} \circ \widetilde{\zeta}'_{\mu'}(X') \circ (L^{\zeta_\mu \zeta'_{\mu'}})^{-1} = (\bar{\zeta}_\mu \circ \Phi^q_{\mu\mu'})(X')$$

The proof can be done using relation 14 and the irreducibility of the considered representations. This justifies the fact that $(\bar{\zeta}_{\lambda_0}, \widetilde{\zeta}'_{\lambda'_0})$ belongs to $A^{\rho\rho'}_{\lambda\lambda'}$ and $(\zeta_{\lambda_0}, \zeta'_{\lambda'_0}) \approx (\bar{\zeta}_{\lambda_0}, \widetilde{\zeta}'_{\lambda'_0})$. Both the representations $\widetilde{\zeta}'_{\mu'}$ and $\bar{\zeta}'_{\mu'}$ are isomorphic to $\bar{\zeta}_\mu \circ \Phi^q_{\mu\mu'}$, so they are isomorphic to each other. Possibly by changing $\widetilde{\zeta}'_{\lambda'_0}$ in its class of local representation of $\mathcal{T}^q_\mu$, we can assume that $\bar{\zeta}'_{\lambda'_0}$ and $\widetilde{\zeta}'_{\lambda'_0}$ are isomorphic. Indeed, change $\widetilde{\zeta}'_{\lambda'_0}$ in its class of local representation of $\mathcal{T}^q_\mu$ is equivalent to take a representation $\widetilde{\eta}'_{\lambda'_0}$ defined by $\widetilde{\zeta}'_{\lambda'_0} \xrightarrow{\gamma_j} \widetilde{\eta}'_{\lambda'_0}$ with $\gamma_j = 1$ for every $j \neq i$. We assert that $\gamma_i$ can be chosen so that $\widetilde{\eta}'_{\lambda'_0}$ is isomorphic to $\bar{\zeta}'_{\lambda'_0}$. We know that the central loads $\widetilde{h}$ and $\bar{h}$, associated with $Q$, of $\widetilde{\zeta}'_{\mu'}$ and $\bar{\zeta}'_{\mu'}$ are the product of the central loads of $\widetilde{\zeta}'_{\lambda'_0}$ and $\bar{\zeta}'_{\lambda'_0}$ on the triangles composing $Q$. Moreover, being $\widetilde{\zeta}'_{\mu'}$ and $\bar{\zeta}'_{\mu'}$ isomorphic, we have $\widetilde{h} = \bar{h}$. Denoting by $\widetilde{h}^1$ and $\widetilde{h}^2$ the central loads of $\widetilde{\zeta}'_{\lambda'_0}$ on the triangles in $Q$ and by $\bar{h}^1, \bar{h}^2$ the ones of $\bar{\zeta}'_{\lambda'_0}$, we observe that the relative central loads of $\widetilde{\eta}'_{\lambda'_0}$ change like $\gamma_i \widetilde{h}^1$ and $\gamma_i^{-1} \widetilde{h}^2$. Then there exists a unique $\gamma_i$ such that $\bar{h}^1 = \gamma_i \widetilde{h}^1$. Doing this choice, the following holds

$$\bar{h}^2 = \frac{\bar{h}}{\bar{h}^1} = \frac{\widetilde{h}}{\gamma_i \widetilde{h}^1} = \gamma_i^{-1} \widetilde{h}^2$$

Now the central loads of $\bar{\zeta}'_{\lambda'_0}$ and $\widetilde{\eta}'_{\lambda'_0}$ on the triangles are equal. It remains to prove that the invariants on the couple of edges in $\lambda'_0$ corresponding to $\lambda_i$ are equal, but this is clear because we already know that, for both the triangles in $Q$, their representations have the same central loads and two invariants of edges coinciding. Recalling that $h^N = x_1 x_2 x_3$ the assertion follows.

Hence we can assume that $\bar{\zeta}'_{\lambda'_0}$ and $\widetilde{\zeta}'_{\lambda'_0}$ are isomorphic. Via Proposition 1.18, there exists a $c \in H_1(S; \mathbb{Z}_N)$ such that $c \cdot \widetilde{\zeta}'_{\lambda'_0} = \bar{\zeta}'_{\lambda'_0}$, hence

$$c \cdot [\zeta_{\lambda_0}, \zeta'_{\lambda'_0}] = c \cdot [\bar{\zeta}_{\lambda_0}, \widetilde{\zeta}'_{\lambda'_0}] = [\bar{\zeta}_{\lambda_0}, c \cdot \widetilde{\zeta}'_{\lambda'_0}] = [\bar{\zeta}_{\lambda_0}, \bar{\zeta}'_{\lambda'_0}]$$

and so the transitivity is proved.

Now suppose that there exist a $c \in H_1(S; \mathbb{Z}_N)$ and an element $[\zeta_{\lambda_0}, \zeta'_{\lambda'_0}] \in \mathscr{A}^{\rho\rho'}_{\lambda\lambda'}$ such that $[\zeta_{\lambda_0}, \zeta'_{\lambda'_0}] = c \cdot [\zeta_{\lambda_0}, \zeta'_{\lambda'_0}]$. This means that, passing on representatives, the couples $(\zeta_{\lambda_0}, \zeta'_{\lambda'_0})$ and $(\zeta_{\lambda_0}, c \cdot \zeta'_{\lambda'_0})$ are $\approx$-equivalent. Because the first terms of the couples are exactly the same, they are in particular related by transition constants all equal to 1, and the same must hold for $\zeta'_{\lambda'_0}$ and $c \cdot \zeta'_{\lambda'_0}$, except possibly along the edge $\lambda_i$, where there are not restrictions. But $(\lambda'_i)^*$ has distinct vertices and, because there are not elements of $H_1(S; \mathbb{Z}_N)$ that act only in the edge $\lambda_i$, we conclude that $c$ must be trivial. Therefore the action is free.

Finally, the relation 15 can be proved exactly in the same way of relation 12 in Theorem 2.1. $\square$

### 2.3.1 An explicit calculation

The previous discussion shows us that the elements in $\mathscr{L}^{\rho\rho'}_{\lambda\lambda'}$ are tensor products of a tensor split isomorphism $L^{\zeta_{\mu_0}\zeta'_{\mu_0}}$ and an isomorphism $L^{\zeta_{\mu_Q}\zeta'_{\mu'_Q}}$ between two irreducible representations $\zeta_{\mu_Q} \circ \Phi^q_{\mu_Q\mu'_{Q'}}$ and $\zeta_{\mu'_{Q'}}$ on the square $Q$. In what follows we want to give an explicit description of this linear isomorphism $L^{\zeta_{\mu_Q}\zeta'_{\mu'_Q}}$.

Redefine the notations: let $Q$ be an ideal square and let $\lambda, \lambda'$ be its ideal triangulations, with edges labelled as in Figure 5. Given $\rho = \{\rho_\lambda, \rho_{\lambda'}\}$ a local representation of $\mathcal{T}^q_Q$, we know that there exists a linear isomorphism $L^{\rho\rho}_{\lambda\lambda'} \colon V_{\lambda'} \to V_\lambda$, unique up to scalar multiplication, such that

$$L^{\rho\rho}_{\lambda\lambda'} \circ \rho_{\lambda'}(X') \circ (L^{\rho\rho}_{\lambda\lambda'})^{-1} = (\rho_\lambda \circ \Phi^q_{\lambda\lambda'})(X')$$

Let us describe in a explicit way this linear isomorphism. We firstly reduce to the standard situation, which means that $\rho_\lambda$ and $\rho_{\lambda'}$ are represented by the tensor product of standard irreducible representations of the triangle algebras (here standard means that the representation sends each generator $X_s$ of $\mathcal{T}^q_T$ in a $N \times N$ matrix, which is a multiple of the $B_i$ described in Proposition 1.4). In order to determine a standard representation of the triangle algebra, we need the following data: a clockwise indexing of the edges of each triangle and the choice of $N$-th roots of the invariants on the edges of the square. We will order the edges of each triangle as described in Figure 5 by red numbers (the square on the left represents the ideal triangulation $\lambda$ and the indexing on its triangles $T_1$ and $T_2$, the square on the right represents the ideal triangulation $\lambda'$ and the indexing on its triangles $T'_1$ and $T'_2$). The representations $\rho_\lambda \colon \mathcal{T}^q_\lambda \to \mathrm{End}(\mathbb{C}^N \otimes \mathbb{C}^N)$ and $\rho_{\lambda'} : \mathcal{T}^q_{\lambda'} \to \text{End}(\mathbb{C}^N \otimes \mathbb{C}^N)$ have the following form

$$\rho_\lambda(X_i) = y_i\, B_2 \otimes B_1$$
$$\rho_\lambda(X_j) = y_j\, B_1 \otimes I$$
$$\rho_\lambda(X_k) = y_k\, I \otimes B_2$$
$$\rho_\lambda(X_l) = y_l\, I \otimes B_3$$
$$\rho_\lambda(X_m) = y_m\, B_3 \otimes I$$

$$\rho_{\lambda'}(X'_i) = v_i\, B_3 \otimes B_2$$
$$\rho_{\lambda'}(X'_j) = v_j\, B_1 \otimes I$$
$$\rho_{\lambda'}(X'_k) = v_k\, B_2 \otimes I$$
$$\rho_{\lambda'}(X'_l) = v_l\, I \otimes B_3$$
$$\rho_{\lambda'}(X'_m) = v_m\, I \otimes B_1$$

where the numbers $y_i, y_j, y_k, y_l, y_m$ are $N$-th roots of $x_i, x_j, x_k, x_l, x_m$ and the numbers $v_i, v_j, v_k, v_l, v_m$ are $N$-th roots of $x_i^{-1}, (1+x_i)x_j, (1+x_i^{-1})^{-1}x_k, (1+x_i)x_l, (1+x_i^{-1})^{-1}x_m$. Moreover, the product of the $y_s$ and the product of the $v_s$ are equal and coincide with the central load of the representations $\rho_\lambda$ and $\rho_{\lambda'}$. Denote by $z_i$ the number in $\mathbb{Z}_N$ such that $v_i = y_i^{-1} q^{2z_i}$ (recall that $\Phi^q_{\lambda\lambda'}(X'_i) = X_i^{-1}$).

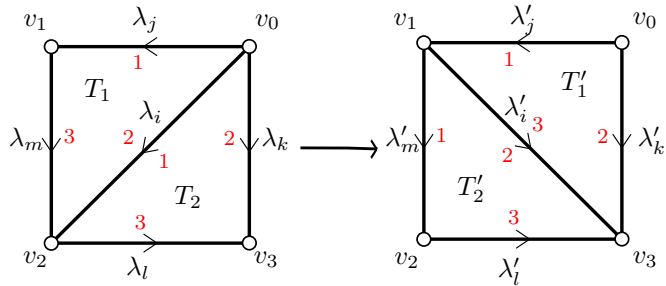

Figure 5: Useful notations

Because of the expression of $\Phi^q_{\lambda\lambda'}$, the representation $\rho_\lambda \circ \Phi^q_{\lambda\lambda'}$ has the following behaviour

$$(\rho_\lambda \circ \Phi^q_{\lambda\lambda'})(X'_i) = y_i^{-1}\, B_2^{-1} \otimes B_1^{-1}$$
$$(\rho_\lambda \circ \Phi^q_{\lambda\lambda'})(X'_j) = y_j\, (I \otimes I + q y_i\, B_2 \otimes B_1) B_1 \otimes I$$
$$(\rho_\lambda \circ \Phi^q_{\lambda\lambda'})(X'_k) = y_k\, (I \otimes I + q y_i^{-1}\, B_2^{-1} \otimes B_1^{-1})^{-1} I \otimes B_2$$
$$(\rho_\lambda \circ \Phi^q_{\lambda\lambda'})(X'_l) = y_l\, (I \otimes I + q y_i\, B_2 \otimes B_1) I \otimes B_3$$
$$(\rho_\lambda \circ \Phi^q_{\lambda\lambda'})(X'_m) = y_m\, (I \otimes I + q y_i^{-1}\, B_2^{-1} \otimes B_1^{-1})^{-1} B_3 \otimes I$$

Now we define, for $a \in \mathbb{N}$, the following function

$$f(a) := \begin{cases} \left(\frac{y_j}{v_j}\right)^a \prod_{u=1}^a (1 + y_i q^{1-2(u+z_i)}) & \text{if } a \neq 0 \\ 1 & \text{if } a = 0 \end{cases}$$

Observe that, because $q$ is a primitive $N$-th root of $(-1)^{N+1}$, we have

$$\prod_{u=1}^N (1 + y_i q^{1-2(u+z_i)}) = 1 + y_i^N = 1 + x_i$$

On the other hand, $\frac{y_j}{v_j}$ is an $N$-th root of $(1+x_i)^{-1}$, indeed

$$\left(\frac{y_j}{v_j}\right)^N = \frac{x_j}{(1+x_i)x_j} = (1+x_i)^{-1}$$

These facts imply immediately that, for every $a \in \mathbb{N}$, the following relation holds

$$f(a+N) = f(a)$$

Moreover, we define the following polynomial

$$p(x) := \sum_{d=0}^{N-1} \left( \frac{y_l y_k y_i}{v_l v_k} \right)^d x^d$$

Observe that $\frac{y_l y_k y_i}{v_l v_k}$ is an $N$-th root of unity, indeed:

$$\left( \frac{y_l y_k y_i}{v_l v_k} \right) = \frac{x_l x_k x_i}{x'_l x'_k} = \frac{x_l x_k x_i}{(1+x_i) x_l (1+x_i^{-1})^{-1} x_k} = 1$$

Now we have all the tools required for the description of the isomorphism $L^{\rho\rho}_{\lambda\lambda'}$, up to scalar multiplication.

**Proposition 2.4.** The map $L^{\rho\rho}_{\lambda\lambda'} \colon \mathbb{C}^N \otimes \mathbb{C}^N \to \mathbb{C}^N \otimes \mathbb{C}^N$ has components $(L^{\rho\rho}_{\lambda\lambda'})^{b,c}_{s,t}$ equal to

$$q^{-s^2 + 2z_i(b-c-z_i) + 2bc} \left( \frac{y_j y_k y_i}{v_j v_k} \right)^{c+z_i} p(q^{2(s+t-c-z_i)}) \left( \sum_{a=0}^{N-1} q^{2a(b-s)} f(a) \right) \quad (16)$$

for varying $s, t, b, c \in \{0, \ldots, N-1\}$, where the indices of $e_{-r, r-(s+t+z_i)}$ are thought as elements of $\mathbb{Z}_N = \{[0], \ldots, [N-1]\}$ and

$$L^{\rho\rho}_{\lambda\lambda'}(e_s \otimes e_t) = \sum_{b,c=0}^{N-1} (L^{\rho\rho}_{\lambda\lambda'})^{b,c}_{s,t} \, e_b \otimes e_c$$

with $\{e_0, \ldots, e_{N-1}\}$ the canonical basis of $\mathbb{C}^N$.

Recall that $L^{\rho\rho}_{\lambda\lambda'}$ is defined up to scalar multiplication, so by multiplying the relations above, for varying $b, c, s, t$, by a common scalar, we obtain another linear isomorphism verifying the property

$$L^{\rho\rho}_{\lambda\lambda'} \circ \rho_{\lambda'}(X') \circ (L^{\rho\rho}_{\lambda\lambda'})^{-1} = (\rho_\lambda \circ \Phi^q_{\lambda\lambda'})(X') \qquad \forall X' \in \mathcal{T}^q_{\lambda'}$$

*Proof of Proposition 2.4.* In what follows we will describe the strategy conducing to the relation 16. Firstly, we choose an indexing on the edges of the triangles different from the one of Figure 5, but more appropriate in order to do an explicit calculation. In particular, we have chosen the indexing in Figure 6. Denote by $\bar{\rho}_\lambda$ and $\bar{\rho}_{\lambda'}$ the standard representations determined by this choice of indexing and by the $N$-th roots $y_s$ and $v_s$ respectively, the same that we have chosen above. In particular, in this case we have the following relations

$$\bar{\rho}_\lambda(X_i) = y_i \, B_1 \otimes B_1$$
$$\bar{\rho}_\lambda(X_j) = y_j \, B_3 \otimes I$$
$$\bar{\rho}_\lambda(X_k) = y_k \, I \otimes B_2$$
$$\bar{\rho}_\lambda(X_l) = y_j \, I \otimes B_3$$
$$\bar{\rho}_\lambda(X_m) = y_m \, B_2 \otimes I$$

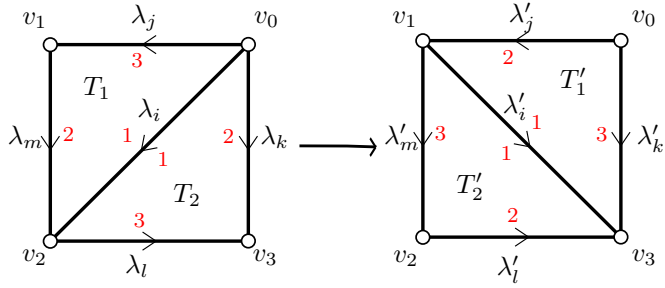

Figure 6: Another indexing

$$\bar{\rho}_{\lambda'}(X'_i) = v_i\, B_1 \otimes B_1$$
$$\bar{\rho}_{\lambda'}(X'_j) = v_j\, B_2 \otimes I$$
$$\bar{\rho}_{\lambda'}(X'_k) = v_k\, B_3 \otimes I$$
$$\bar{\rho}_{\lambda'}(X'_l) = v_l\, I \otimes B_2$$
$$\bar{\rho}_{\lambda'}(X'_m) = v_m\, I \otimes B_3$$

As previously done, we will denote by $z_i \in \mathbb{Z}_N$ the number verifying $v_i = q^{2z_i} y_i^{-1}$. We denote the vector $e_s \otimes e_t$ briefly by $e_{s,t}$. The representation $\bar{\rho}_\lambda \circ \Phi^q_{\lambda\lambda'}$ has the following behaviour

$$(\bar{\rho}_\lambda \circ \Phi^q_{\lambda\lambda'})(X'_i) = y_i^{-1}\, B_1^{-1} \otimes B_1^{-1}$$
$$(\bar{\rho}_\lambda \circ \Phi^q_{\lambda\lambda'})(X'_j) = y_j\, (I \otimes I + q y_i\, B_1 \otimes B_1) B_3 \otimes I$$
$$(\bar{\rho}_\lambda \circ \Phi^q_{\lambda\lambda'})(X'_k) = y_k\, (I \otimes I + q y_i^{-1}\, B_1^{-1} \otimes B_1^{-1})^{-1} I \otimes B_2$$
$$(\bar{\rho}_\lambda \circ \Phi^q_{\lambda\lambda'})(X'_l) = y_l\, (I \otimes I + q y_i\, B_1 \otimes B_1) I \otimes B_3$$
$$(\bar{\rho}_\lambda \circ \Phi^q_{\lambda\lambda'})(X'_m) = y_m\, (I \otimes I + q y_i^{-1}\, B_1^{-1} \otimes B_1^{-1})^{-1} B_2 \otimes I$$

The main benefits of this choice are

- the matrices $I \otimes I + q y_i\, B_1 \otimes B_1$ and $(I \otimes I + q y_i^{-1}\, B_1^{-1} \otimes B_1^{-1})^{-1}$ are diagonal;

- it is very simple to find the family of isomorphisms $\psi$ verifying

$$\psi^{-1} \circ (\bar{\rho}_\lambda \circ \Phi^q_{\lambda\lambda'})(X'_i) \circ \psi = \bar{\rho}_{\lambda'}(X'_i)$$

Indeed, they both are diagonal matrices, so it is sufficient to ask that $\psi$ carries the $\alpha$-eigenspace of $\bar{\rho}_{\lambda'}(X_i)$ in the $\alpha$-eigenspace of $(\bar{\rho}_\lambda \circ \Phi^q_{\lambda\lambda'})(X'_i)$ for every $\alpha$. In particular, such a $\psi$ has the following behaviour

$$e_{s,t} \longmapsto \sum_{r=0}^{N-1} a_{r,s,t}\, e_{-r,\,r-(s+t+z_i)}$$

where the indices of $e_{-r,r-(s+t+z_i)}$ must be thought as elements of $\mathbb{Z}_N = \{[0], \ldots, [N-1]\}$.

Now we require that $\psi$ carries the whole representation $\bar{\rho}_\lambda \circ \Phi^q_{\lambda\lambda'}$ in $\bar{\rho}_{\lambda'}$ by conjugation and we find equations in the constants $a_{r,s,t}$, which determine these elements $a_{r,s,t}$ up to a common multiplicative scalar. In particular, we obtain the following equations

$$\begin{cases} y_j(1 + y_i q^{-2(s+t+z_i+1)+1}) q^{1+2(r-1)} a_{r-1,s,t} = v_j\, a_{r,s+1,t} \\ y_k(1 + y_i^{-1} q^{2(s+t+z_i-1)+1})^{-1} a_{r,s,t} = v_k q^{1-2s}\, a_{r,s-1,t} \\ y_l(1 + y_i q^{-2(s+t+z_i+1)+1}) q^{1+2(s+t+z_i-r)} a_{r,s,t} = v_l\, a_{r,s,t+1} \\ y_m(1 + y_i^{-1} q^{2(s+t+z_i-1)+1})^{-1} a_{r+1,s,t} = v_m\, a_{r,s,t-1} \end{cases}$$

conducing to the following expression

$$a_{r,s,t} = \left(\frac{v_k}{y_k y_i}\right)^s \left(\frac{y_j y_k y_i}{v_j v_k}\right)^r \left(\frac{y_l}{v_l}\right)^t q^{(t-r)^2 + 2(s+t-r)z_i + 2st} \prod_{u=1}^{s+t}(1 + y_i q^{1-2(u+z_i)})\, a_{0,0,0}$$

Now, defining $\xi \colon \mathbb{C}^N \to \mathbb{C}^N$ the isomorphism

$$\begin{array}{rcl} \xi \colon \mathbb{C}^N & \longrightarrow & \mathbb{C}^N \\ e_k & \longmapsto & \frac{1}{\sqrt{N}} \sum_{h=0}^{N-1} q^{2hk+h^2}\, e_h \end{array}$$

we observe that the following relations hold

$$(\xi \otimes \xi^{-1}) \circ \bar{\rho}_{\lambda'} \circ (\xi^{-1} \otimes \xi) = \rho_{\lambda'}$$
$$(\xi^{-1} \otimes I) \circ \bar{\rho}_{\lambda} \circ (\xi \otimes I) = \rho_{\lambda}$$

The point is that $\xi$ is the linear isomorphism that change a standard triangular representation in another standard triangular representation simply by rotate the indexing. More precisely, if $X_1, X_2, X_3$ are generators of $\mathcal{T}_T^q$, corresponding to edges $\lambda_1, \lambda_2, \lambda_3$, ordered clockwise, and if $\eta$ is a representation of $\mathcal{T}_T^q$ defined by

$$\eta(X_i) = u_i\, B_i$$

for every $i \in \{1, 2, 3\}$, then $\xi^{-1} \circ \eta(\cdot) \circ \xi$ verifies

$$\xi^{-1} \circ \eta(X_i) \circ \xi = u_i\, B_{i+1}$$

where the indices are in $\mathbb{Z}_3 = \{[1], [2], [3]\}$. So we obtain that the composition

$$(\xi^{-1} \otimes I) \circ \psi \circ (\xi^{-1} \otimes \xi)$$

verifies the property defining $L^{\rho\rho}_{\lambda\lambda'}$. The relation 16 can be found developing the composition, where we have chosen $a_{0,0,0} = N\sqrt{N}$. □

## 3 The elementary properties

In this Section we will focus on the properties verified by the objects $(\mathscr{L}^{\rho\rho'}_{\lambda\lambda'}, \psi^{\rho\rho'}_{\lambda\lambda'})$ just defined. In particular, we will investigate on the relations that will conduce to the Fusion and Composition properties in the general case. The first part is dedicated to the "baby" version of the Fusion property. The second Subsection will request some efforts and will conduce us to the proof of a technical Lemma that will be useful in the last Subsection, where we will explicit the elementary version of the Composition property (here elementary means restricted to the case in which $\lambda, \lambda$ differ by an elementary move).

## 3.1 Elementary Fusion property

If $S$ is obtained by fusion from $R$ along some boundary components, there is a natural map of projection $\pi\colon R \to S$. Given $\mu \in \Lambda(R)$ an ideal triangulation and $\lambda \in \Lambda(S)$ the induced ideal triangulation on $S$, the map $\pi$ induces an identification of $\Gamma_{R,\mu}$ with a subgraph of $\Gamma_{S,\lambda}$. Moreover, thinking to $\Gamma_{R,\mu}$ as a deformation retract of $R$, the map $\pi_*\colon H_1(R;\mathbb{Z}_N) \to H_1(S;\mathbb{Z}_N)$ is injective, because the map obtained from the inclusion of $\Gamma_{R,\mu}$ in $\Gamma_{S,\lambda}$ on $H_1(\,\cdot\,;\mathbb{Z}_N)$ is.

**Lemma 3.1.** Let $R$ be a surface as above and $S$ be obtained by fusion from $R$. Fix $\eta = \{\eta_\mu\colon \mathcal{T}_\mu^q \to \mathrm{End}(W_\mu)\}_{\mu \in \Lambda(R)}$, $\eta' = \{\eta'_\mu\colon \mathcal{T}_\mu^q \to \mathrm{End}(W'_\mu)\}_{\mu \in \Lambda(R)}$ two isomorphic local representations of $\mathcal{T}_R^q$ and $\rho = \{\rho_\lambda\colon \mathcal{T}_\lambda^q \to \mathrm{End}(V_\lambda)\}_{\lambda \in \Lambda(S)}$, $\rho' = \{\rho'_\lambda\colon \mathcal{T}_\lambda^q \to \mathrm{End}(V'_\lambda)\}_{\lambda \in \Lambda(S)}$ two isomorphic local representations of $\mathcal{T}_S^q$, with $\rho$ and $\rho'$ obtained by fusion from $\eta$ and $\eta'$, respectively. Then, for every $\mu, \mu' \in \Lambda(R)$ that differ by diagonal exchange or a re-indexing, if $\lambda, \lambda' \in \Lambda(S)$ are the corresponding ideal triangulations on $S$, there exists a natural inclusion $j\colon \mathcal{L}_{\mu\mu'}^{\eta\eta'} \to \mathcal{L}_{\lambda\lambda'}^{\rho\rho'}$ such that, for every $L \in \mathcal{L}_{\mu\mu'}^{\eta\eta'}$, the following holds

$$j(c \cdot L) = \pi_*(c) \cdot j(L)$$

for every $c \in H_1(R;\mathbb{Z}_N)$.

*Proof.* We will prove only the case in which $\lambda$ and $\lambda'$ differ by a diagonal exchange, the other situations are analogous. Moreover, we will use the sets $\mathcal{A}_{\lambda\lambda'}^{\rho\rho'}$ instead of $\mathcal{L}_{\lambda\lambda'}^{\rho\rho'}$, in order to describe the relations with the action in a more explicit way. On the $\mathcal{L}$-level, the map $j$ will be just the inclusion as subsets of $\mathrm{Hom}(V'_{\lambda'}, V_\lambda)/\doteq$.

Fix $\mu, \mu' \in \Lambda(R)$ and $\lambda, \lambda' \in \Lambda(S)$ as in the statement, with $\lambda = \Delta_i(\lambda')$ and $\mu = \Delta_i(\mu')$. It is clear that the surfaces $S_0$ and $R_0$, obtained by splitting $S$ and $R$ along $\lambda$ and $\mu$ respectively, can be identified and analogously for $S'_0$ and $R'_0$, obtained by splitting $S$ and $R$ along $\lambda'$ and $\mu'$. If $\eta_\mu$, local representation of $\mathcal{T}_\mu^q$, represents $\rho_\lambda$, then $\mathscr{F}_{R_0}(\eta_\mu) \subseteq \mathscr{F}_{S_0}(\rho_\lambda)$, by definition. From this fact we deduce an inclusion $i$ of $A_{\mu\mu'}^{\eta\eta'}$ in $A_{\lambda\lambda'}^{\rho\rho'}$. The map $j$ will be the application induced by $i$ from $\mathcal{A}_{\mu\mu'}^{\eta\eta'}$ to $\mathcal{A}_{\lambda\lambda'}^{\rho\rho'}$. We need to prove the good definition of $j$ and the injectivity.

Let $(\zeta, \zeta') \in \mathscr{F}_{R_0}(\eta_\mu) \times \mathscr{F}_{R_0}(\eta'_{\mu'})$ be an element in $A_{\mu\mu'}^{\eta\eta'}$ and denote by $[\zeta, \zeta']_S$ its image in $\mathcal{A}_{\lambda\lambda'}^{\rho\rho'}$. Take $(\bar{\zeta}, \bar{\zeta}')$ another representative of $[\zeta, \zeta']_R \in \mathcal{A}_{\mu\mu'}^{\eta\eta'}$, the equivalence class of $(\zeta, \zeta')$ in $A_{\mu\mu'}^{\eta\eta'}$. Then there exist transition constants $\alpha_j$ and $\beta_j$, one for each internal edge of $\mu$, such that $\zeta \xrightarrow{\alpha_j} \bar{\zeta}$ and $\zeta' \xrightarrow{\beta_j} \bar{\zeta}'$, with $\alpha_j = \beta_j$ for every $j \neq i$. Note that the representations $\zeta$ and $\bar{\zeta}$ need to coincide on the variables corresponding to the boundary edges of $R$. In particular this means that the elements $\zeta$ and $\bar{\zeta}$, as representations in $\mathscr{F}_{S_0}(\rho_\lambda)$, have transition constants equal to 1 for every $\lambda_j$ that is the result of the identification of a couple of boundary components in $R$, and equal to $\alpha_j \in \mathbb{C}^*$ otherwise. The same must hold for $\zeta'$ and $\bar{\zeta}'$, so in particular $[\zeta, \zeta']_S = [\bar{\zeta}, \bar{\zeta}']_S$, which proves the good definition of $j$.

Now take $(\zeta, \zeta'), (\bar{\zeta}, \bar{\zeta}') \in A_{\mu\mu'}^{\eta\eta'}$ and assume that $[\zeta, \zeta']_S = [\bar{\zeta}, \bar{\zeta}']_S$. This means that there exist transitions constants $\alpha_j$ and $\beta_j$, one for each internal edge of $\lambda$, such that $\zeta \xrightarrow{\alpha_j} \bar{\zeta}$ and $\zeta' \xrightarrow{\beta_j} \bar{\zeta}'$, with $\alpha_j = \beta_j$ if $j \neq i$. On the other hand, the representations $\zeta$ and $\bar{\zeta}$ belong to $\mathscr{F}_{R_0}(\eta_\mu)$, so they must coincide on the variables corresponding to the boundary edges of $R$. This implies that $\alpha_j = 1$ for every $j$ such that $\lambda_j^* \notin \Gamma_{R,\mu}$ (we are identifying $\Gamma_{R,\mu}$ with its image in $\Gamma_{S,\lambda}$ under $\pi$). In the same way we can see that $\beta_j = 1$ for every $j$ such that $(\lambda_j')^* \notin \Gamma_{R,\mu'}$, and these observations lead to the equality $[\zeta, \zeta']_R = [\bar{\zeta}, \bar{\zeta}']_R$. Hence we have concluded the proof of the injectivity.

Finally observe that, for every $c \in H_1(R; \mathbb{Z}_N)$, we have

$$\begin{aligned} j(c \cdot [\zeta, \zeta']_R) &= j([\zeta, c \cdot \zeta']_R) \\ &= [\zeta, \pi_*(c) \cdot \zeta']_S \\ &= \pi_*(c) \cdot [\zeta, \zeta']_S \\ &= \pi_*(c) \cdot j([\zeta, \zeta']_R) \end{aligned}$$

□

## 3.2 A technical Lemma

In the previous Section we have given a presentation of the elements in $\mathscr{A}_{\lambda\lambda'}^{\rho\rho'}$ in terms of equivalence classes of representations on the surfaces $S_0$ or $S_0'$. In this Subsection we are going to prove a Lemma that will give an alternative construction of the sets $\mathscr{A}_{\lambda\lambda'}^{\rho\rho'}$ in terms of equivalence classes of local representations in a intermediate common level $R'$ between $S$ and the surfaces $S_0$, $S_0'$.

In other words, we want to represent local representations on $S$ with local representations on another surface $R'$, which can be obtained by fusion from both $S_0$ and $S_0'$ and which is a splitting of $S$ along certain edges of an ideal triangulation of $S$. In addition, we will require that $R'$ is a disjoint union of ideal polygons. The first important example of this situation is the surface $R$ that we have introduced in the Subsection 2.3.

In the setting that this technical Lemma will allow us to introduce, the proof of the Elementary Composition property in the next Section will be simpler and more expressive.

### 3.2.1 Diagonal exchange

Let $\lambda, \lambda' \in \Lambda(S)$ be two ideal triangulations that differ by a diagonal exchange along $\lambda_i$. We have denoted by $S_0$ and $S_0'$ the surfaces, obtained by splitting $S$ along $\lambda$ and $\lambda'$ respectively, endowed with the triangulations $\lambda_0$ and $\lambda_0'$. Moreover, we have defined $R$ as the surface obtained by splitting $S$ along all the edges except for $\lambda_i$, on which we have the ideal triangulations $\mu = \mu_Q \sqcup \mu_0$ and $\mu' = \mu_Q' \sqcup \mu_0$ induced by $\lambda$ and $\lambda'$.

The triangulations $\lambda_0$ and $\lambda_0'$ are the result of the splitting of $\mu$ and $\mu'$ along $\mu_i$ and $\mu_i'$, diagonals of the square $Q$ in $R$. Now we take an intermediate surface $R'$ between $R$ and $S$, that is a surface obtained by splitting $S$ along certain $\lambda_j$, with $j \neq i$, with induced ideal triangulations $\nu$ and $\nu'$. Furthermore, we assume that $R'$ is the disjoint union of ideal polygons. We represent the situation of surfaces and triangulations related by fusion in the diagram on the right, where an arrow from $A$ to $B$ means that $B$ is obtained by fusion from $A$, and on the sides there are the relative triangulations.

Fixed

$$\rho = \{\rho_\lambda \colon \mathcal{T}_\lambda^q \to \mathrm{End}(V_\lambda)\}_{\lambda \in \Lambda(S)}$$
$$\rho' = \{\rho'_\lambda \colon \mathcal{T}_\lambda^q \to \mathrm{End}(V'_\lambda)\}_{\lambda \in \Lambda(S)}$$

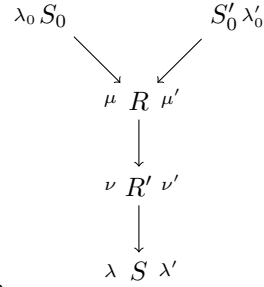

two local representations of $\mathcal{T}_S^q$, we define the following sets

$$\mathscr{F}_{R'}(\rho_\lambda) := \{\zeta_\nu \in \mathrm{Repr}_{loc}(\mathcal{T}_\nu^q, V_\lambda) \mid \zeta_\nu \text{ represents } \rho_\lambda\}$$
$$\mathscr{F}_{R'}(\rho'_{\lambda'}) := \{\zeta'_{\nu'} \in \mathrm{Repr}_{loc}(\mathcal{T}_{\nu'}^q, V'_{\lambda'}) \mid \zeta'_{\nu'} \text{ represents } \rho'_{\lambda'}\}$$

Now we introduce a set $B_{\lambda\lambda'}^{\rho\rho'}$ that will perform the role of $A_{\lambda\lambda'}^{\rho\rho'}$ in the new setting:

$$B_{\lambda\lambda'}^{\rho\rho'} := \{(\zeta_\nu, \zeta'_{\nu'}) \in \mathscr{F}_{R'}(\rho_\lambda) \times \mathscr{F}_{R'}(\rho'_{\lambda'}) \mid \zeta_\nu \circ \Phi_{\nu\nu'}^q \text{ is isomorphic to } \zeta'_{\nu'}\}$$

Fix two elements $(\zeta_\nu, \zeta'_{\nu'}), (\bar\zeta_\nu, \bar\zeta'_{\nu'})$ of $B_{\lambda\lambda'}^{\rho\rho'}$. For any choice of representatives $\eta_{\lambda_0}$ in $\mathscr{F}_{S_0}(\zeta_\nu)$, $\bar\eta_{\lambda_0}$ in $\mathscr{F}_{S_0}(\bar\zeta_\nu)$, both $\eta_{\lambda_0}$ and $\bar\eta_{\lambda_0}$ represents $\rho_\lambda$ on $S$. This means that there exist transition constants $\alpha_j$ such that $\eta_{\lambda_0} \xrightarrow{\alpha_j} \bar\eta_{\lambda_0}$, with $\alpha_j$ that corresponds to the edge $\lambda_j^*$ in the dual graph $\Gamma_{S,\lambda}$. In the same way, fixed $\eta'_{\lambda'_0}$ in $\mathscr{F}_{S'_0}(\zeta'_{\nu'})$ and $\bar\eta'_{\lambda'_0}$ in $\mathscr{F}_{S'_0}(\bar\zeta'_{\nu'})$, there exist transition constants $\beta_j$, indexed by the edges of $\Gamma_{S,\lambda'}$, such that $\eta'_{\lambda'_0} \xrightarrow{\beta_j} \bar\eta'_{\lambda'_0}$. The transition constants $\alpha_j$ and $\beta_j$ can be naturally compared for every $j \neq i$. We will say that $(\zeta_\nu, \zeta'_{\nu'})$ and $(\bar\zeta_\nu, \bar\zeta'_{\nu'})$ are $\approx_{R'}$-equivalent if $\alpha_j = \beta_j$ for every $j$ such that $\lambda_j^*$ does not belong to the subgraph $\Gamma_{R',\nu} \subset \Gamma_{S,\lambda}$, for any choice of the representatives $\eta_{\lambda_0}$, $\bar\eta_{\lambda_0}$, $\eta'_{\lambda'_0}$ and $\bar\eta'_{\lambda'_0}$. This definition does not depend on the choices of representatives, because we are not comparing the transition constants on the edges of $\Gamma_{R',\nu}$ and $\Gamma_{R',\nu'}$, which are the only ones that can be modified by a different choice of representatives. We will denote by $\mathscr{B}_{\lambda\lambda'}^{\rho\rho'}$ the quotient of $B_{\lambda\lambda'}^{\rho\rho'}$ by the equivalence relation $\approx_{R'}$.

Analogously to what previously done, we define a map

$$\bar p \colon B_{\lambda\lambda'}^{\rho\rho'} \longrightarrow \mathrm{Hom}(V'_{\lambda'}, V_\lambda)/\doteq$$

carrying an element $(\zeta_\nu, \zeta'_{\nu'})$ in the isomorphism $M^{\zeta_\nu, \zeta'_{\nu'}}$ that verifies

$$(\zeta_\nu \circ \Phi_{\nu\nu'}^q)(X') = M^{\zeta_\nu, \zeta'_{\nu'}} \circ \zeta'_{\nu'}(X') \circ (M^{\zeta_\nu, \zeta'_{\nu'}})^{-1}$$

Note that there is a natural map $f \colon \mathscr{A}_{\lambda\lambda'}^{\rho\rho'} \to \mathscr{B}_{\lambda\lambda'}^{\rho\rho'}$, which sends a couple $[\zeta_{\lambda_0}, \zeta'_{\lambda'_0}] \in \mathscr{A}_{\lambda\lambda'}^{\rho\rho'}$ in the couple $[\zeta_\nu, \zeta'_{\nu'}] \in \mathscr{B}_{\lambda\lambda'}^{\rho\rho'}$, where $\zeta_\nu$ is represented by $\zeta_{\lambda_0}$ and $\zeta'_{\nu'}$ is represented by $\zeta'_{\lambda'_0}$. It is very easy to check that this map is well defined.

Now we are able to give the statement of the announced technical lemma:

**Lemma 3.2.** In the above notations, the following hold:

- the map $f \colon \mathscr{A}_{\lambda\lambda'}^{\rho\rho'} \to \mathscr{B}_{\lambda\lambda'}^{\rho\rho'}$ is a bijection;
- the following diagram is commutative

$$\begin{array}{ccc}
\mathscr{A}^{\rho\rho'}_{\lambda\lambda'} & \xrightarrow{f} & \mathscr{B}^{\rho\rho'}_{\lambda\lambda'} \\
& \searrow \widetilde{p} \quad \swarrow \bar{p} & \\
& \mathrm{Hom}(V'_{\lambda'}, V_\lambda)/\doteq &
\end{array}$$

In particular $\mathscr{L}^{\rho\rho'}_{\lambda\lambda'} = \mathrm{Im}\,\widetilde{p} = \mathrm{Im}\,\bar{p}$.

Before dealing with the proof, we want to remark the consequences of this fact. Thanks to this statement, it is not important to represent an element of $\mathscr{A}^{\rho\rho'}_{\lambda\lambda'}$ as an equivalence class of irreducible representations on the algebras

$$\mathcal{T}^q_{T_1} \otimes \cdots \otimes \mathcal{T}^q_{T_m} \qquad\qquad \mathcal{T}^q_{T'_1} \otimes \cdots \otimes \mathcal{T}^q_{T'_m}$$

but it is sufficient to choose a certain surface $R'$, which is a disjoint union of polygons, and take couples of local representations on $R'$, with a proper equivalence relation that generalizes the one defined in original construction of $\mathscr{A}^{\rho\rho'}_{\lambda\lambda'}$. In other words, in order to obtain all the intertwining operators in $\mathscr{L}^{\rho\rho'}_{\lambda\lambda'}$ via the action of $H_1(S; \mathbb{Z}_N)$, it is not important to split $S$ in all the ideal triangles that compose it but is sufficient to split the surface in simple connected pieces.

*Proof of Lemma 3.2.* Firstly we will prove the surjectivity of $f$. Fixed $[\zeta_\nu, \zeta'_{\nu'}]$ in $\mathscr{B}^{\rho\rho'}_{\lambda\lambda'}$, we want to find an element $(\eta_{\lambda_0}, \eta'_{\lambda'_0}) \in A^{\rho\rho'}_{\lambda\lambda'}$ such that the local representations $\eta_\nu$ and $\eta'_{\nu'}$, represented by $\eta_{\lambda_0}$ and $\eta'_{\lambda'_0}$ on $\mathcal{T}^q_\nu$ and $\mathcal{T}^q_{\nu'}$, respectively, verify $[\eta_\nu, \eta'_{\nu'}]_{R'} = [\zeta_\nu, \zeta'_{\nu'}]_{R'} \in \mathscr{B}^{\rho\rho'}_{\lambda\lambda'}$. Take a representative $\zeta_{\lambda_0}$ of the local representation $\zeta_\nu$ and analogously $\zeta'_{\lambda'_0}$ of $\zeta'_{\nu'}$. We denote by $\zeta_\mu$ and $\zeta'_{\mu'}$ the corresponding representations on $\mathcal{T}^q_\mu$ and $\mathcal{T}^q_{\mu'}$. Take $\eta_{\lambda_0} := \zeta_{\lambda_0} \in \mathscr{F}_{S_0}(\rho_\lambda)$. We want to change the element $\zeta'_{\lambda'_0}$ with an $\eta'_{\lambda'_0} \in \mathscr{F}_{S'_0}(\zeta'_{\nu'})$ in such a way that the corresponding $\eta'_{\mu'}$ is isomorphic to $\eta_\mu \circ \Phi^q_{\mu\mu'} = \zeta_\mu \circ \Phi^q_{\mu\mu'}$ (observe that $\zeta_\mu \circ \Phi^q_{\mu\mu'}$ is well defined because $\zeta_{\lambda_0}$ is a representative of $\rho_\lambda$, and $\rho_\lambda \circ \Phi^q_{\lambda\lambda'}$ makes sense). In other words, we need to find transition constants $\alpha_j$, one for every edge $(\nu'_j)^*$ in the graph $\Gamma_{R',\nu'}$, such that, if $\eta'_{\lambda'_0}$ verifies the following relation

$$\zeta'_{\lambda'_0} \xrightarrow{\alpha_j} \eta'_{\lambda'_0}$$

as elements of $\mathscr{F}_{S'_0}(\zeta'_{\nu'})$, then the local representation $\eta'_{\mu'}$, represented by $\eta'_{\lambda'_0}$ on $\mathcal{T}^q_{\mu'}$, has the same invariants of the ones of $\zeta_\mu \circ \Phi^q_{\mu\mu'}$. Exhibiting such a $\eta'_{\lambda'_0}$, we will find a couple $(\eta_{\lambda_0}, \eta'_{\lambda'_0})$ that belongs to $A^{\rho\rho'}_{\lambda\lambda'}$, because by construction $\eta'_{\mu'}$ is isomorphic to $\eta_\mu \circ \Phi^q_{\mu\mu'}$, and such that $\eta_\nu = \zeta_\nu$, $\eta'_{\nu'} = \zeta'_{\nu'}$.

Remember that $\zeta_\nu \circ \Phi^q_{\nu\nu'}$ and $\zeta'_{\nu'}$ are isomorphic, so the invariants of all the edges of $(R', \nu')$ and the central loads of every component of $R'$ must coincide. In particular, for every $j$ such that the edge $\mu'_j$ in $\partial R$ goes in $\partial R'$ through the fusion, the invariant $\zeta'_{\mu'}((X'_j)^N)$ is the same of $\zeta_\mu \circ \Phi^q_{\mu\mu'}$. Hence we have to find the $\alpha_s$ in order to make coincide the central loads of the connected components of $R$ and the invariants of $\eta'_{\mu'}$ associated with the edges of $\mu'$ that are fused in $\nu'$.

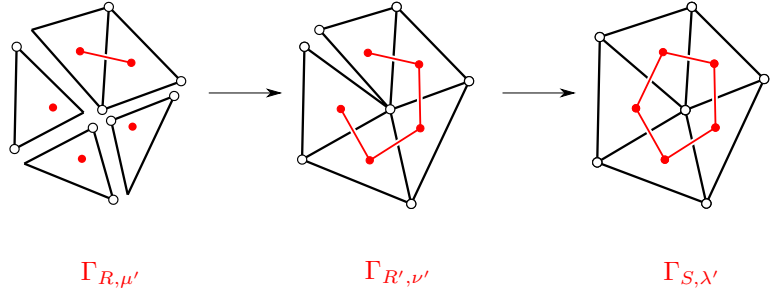

Figure 7: An example of $R \to R' \to S$ and their dual graphs

We will focus at the moment on the component $R'_1$ of $R'$ containing the edge $\nu_i$, along which we make a diagonal exchange. The same procedure that we are going to describe can be applied to each component and will lead to the conclusion. Take the graph $\Gamma_{R',\eta'}$ and denote by $\Gamma_0$ the component of $\Gamma_{R',\eta'}$ corresponding to $R'_1$. $\Gamma_0$ is a tree because is a deformation retract of $R'_1$, which is simply connected by hypothesis. We are going to describe a recursive procedure, with

INPUT: a sub-tree $\Gamma_k$ of $\Gamma_0$, containing $(\nu'_i)^*$, and a transition constant $\alpha_j$ for every $(\nu'_j)^*$ that is in $\Gamma_0^{(1)} \setminus \Gamma_k^{(1)}$ verifying: if $(\zeta'_{\lambda'_0})^{(k)}$ denotes the representation in $\mathscr{F}_{S'_0}(\zeta'_{\nu'})$ defined by

$$\zeta'_{\lambda'_0} \xrightarrow{\beta_j} (\zeta'_{\lambda'_0})^{(k)}$$

where $\beta_j = \alpha_j$ when $(\nu'_j)^* \in \Gamma_0^{(1)}$, $\beta_j = 1$ otherwise, then $(\zeta'_{\mu'})^{(k)}$, the local representation of $\mathcal{T}_{\mu'}^q$ represented by $(\zeta'_{\lambda'_0})^{(k)}$, has the same invariants of $\zeta_\mu \circ \Phi^q_{\mu\mu'}$ on all the couples of edges corresponding to $\Gamma_0^{(1)} \setminus \Gamma_k^{(1)}$ and the same central loads on all the triangles corresponding to vertices in $\Gamma_0^{(0)} \setminus \Gamma_k^{(0)}$;

OUTPUT: a sub-tree $\Gamma_{k+1}$ of $\Gamma_k$, obtained by removing a certain edge $(\nu'_n)^*$, $n \neq i$ and a vertex corresponding to a triangle $T'_s$, and a new transition constant $\alpha_n$, associated with $(\nu'_n)^*$ such that the conditions in the input are verified by $\Gamma_{k+1}$ instead of $\Gamma_{k+1}$.

The algorithm ends when $\Gamma_k$ is composed of the only edge $(\nu'_i)^*$ and its ends. Before describing the procedure, we want to convince ourselves that the final transition constants $(\alpha_j)_j$ provide the desired representation. By construction, the resulting $\eta'_{\lambda'_0} := (\zeta'_{\lambda'_0})^{(k)}$ leads to a representation $\eta'_{\mu'}$ that has the proper central loads and edge invariants on every triangle. The last thing we need to check is that the invariants on the square $Q$ are correct too.

Recall that the central load of a fusion is the product of the central loads of the glued terms. We already know that the central load on $R'_1$ of $\zeta_\nu \circ \Phi^q_{\nu\nu'}$ is equal to the one of $\zeta'_{\nu'}$ and we have constructed a representation $\zeta_{\mu'}$ that has the same loads of $\zeta_\mu \circ \Phi^q_{\mu\mu'}$ on all the triangles not contained in $Q$, so it is immediate to check that the same holds on the square $Q$, by the Fusion property. With analogous observations we can check that also the invariants on the boundary of $Q$ have to be the ones of $\zeta_\mu \circ \Phi_{\mu\mu'}$. In order to conclude the proof of the surjectivity, it is sufficient to repeat the procedure on the other components of $\Gamma_{R',\nu'}$, removing from the conditions on the input the restrictions on $(\nu'_i)^*$. In these cases the procedure ends with transition constants that conduce to a representation with the proper invariants on all the triangles composing the fixed component.

Now we can describe the algorithm. $\Gamma_k$ is a tree, so we can select a leaf in it, i. e. a vertex with valence equal to 1. Assume that the vertex corresponds to the triangle $T'_s$ of the triangulation of $R'_1$. By hypothesis this vertex is on the side of a unique cell $(\nu'_n)^* \in \Gamma_k^{(1)}$, dual of the edge $\nu'_n$. If $n = i$ and there are not any other leaves, then the tree $\Gamma_k$ is the graph of the only square $Q$, and so the algorithm ends. Otherwise, replace the first leaf considered with this one.

Because $(T'_s)^*$ is a leaf, the $(\alpha_j)_j$ selected in the previous steps lead to a representation $(\zeta'_{\lambda'_0})^{(k)}$ that has the correct invariants on two of the three sides of $T'_s$, the ones different from $\nu'_n$. Now we want to select a transition constant $\alpha_n$ in order to make correct also the central load of $T'_s$ and the invariant of the edge of $T'_s$ corresponding to $\nu'_n$. Suppose that the sides of $T'_s$ in $\mu'$ are labelled as $\mu'_l, \mu'_m$ and $\mu'_n$, and that the invariants prescribed by $\zeta_\mu \circ \Phi^q_{\mu\mu'}$ are $y_l, y_m, y_n$ and $h$. Moreover, denote by $\bar{y}_l, \bar{y}_m, \bar{y}_n, \bar{h}$ the invariants of the tensor term of $(\zeta'_{\lambda'_0})^{(k)}$ related to $T'_s$. By hypothesis, we already know that the following hold

$$y_l = \bar{y}_l$$
$$y_m = \bar{y}_m$$

If $\alpha_n$ is the transition constant associated with $\nu'_n$, it is immediate to check that the suitable multiplication by $\alpha_n$ of the tensor term of $(\zeta'_{\lambda'_0})^{(k)}$ associated with $T'_s$ modifies the invariants $\bar{y}_u, \bar{h}$ as follows

$$\bar{y}_l \longrightarrow \bar{y}_l$$
$$\bar{y}_m \longrightarrow \bar{y}_m$$
$$\bar{y}_n \longrightarrow \alpha_n^{N\varepsilon(n,s)} \bar{y}_n$$
$$\bar{h} \longrightarrow \alpha_n^{\varepsilon(n,s)} \bar{h}$$

Now fix a certain $N$-th root $\beta$ of $\left(\frac{y_n}{\bar{y}_n}\right)^{\varepsilon(n,s)}$ and define $\alpha_n := q^{2t}\beta$, with $t \in \mathbb{Z}_N$

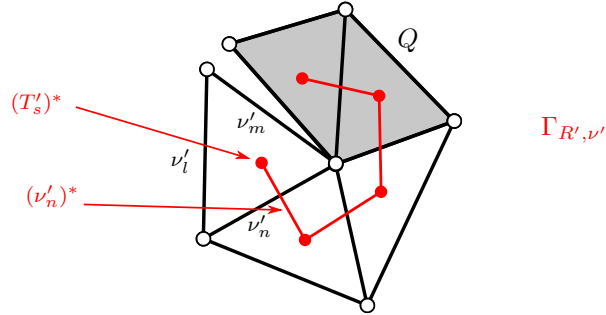

Figure 8: The first step of the algorithm

to be determined. By construction

$$\bar{y}_n \longrightarrow (q^{2t}\beta)^{N\varepsilon(n,s)}\bar{y}_n = q^{2Nt\varepsilon(n,s)}\left(\frac{y_n}{\bar{y}_n}\right)^{\varepsilon(n,s)^2}\bar{y}_n$$

$$= y_n$$

Now we want to choose $t \in \mathbb{Z}_N$ in order to send $\bar{h}$ in $h$. Recall that $h^N = y_l y_m y_n$ and $\bar{h}^N = \bar{y}_l \bar{y}_m \bar{y}_n$, so with every choice of $t$ we have $h^N = \bar{h}^N$. On the other hand

$$\bar{h} \longrightarrow q^{2t\varepsilon(n,s)}\beta\bar{h}$$

and this implies that we can realize, by changing $t$, all the possible $N$-th roots of $\bar{y}_l\bar{y}_m\bar{y}_n = y_l y_m y_n = h^N$, and then there exists a $\bar{t}$ such that $\bar{h} \to h$.

Denote by $\Gamma_{k+1}$ the tree obtained by removing $(\nu'_n)^*$ and $(T'_s)^*$ from $\Gamma_k$. Let us verify that $\Gamma_{k+1}$ has all the properties in order to repeat the algorithm on it: by construction the representation $(\zeta'_{\lambda'_0})^{(k+1)}$ has the same central loads of $\zeta_\mu \circ \Phi^q_{\mu\mu'}$ on all the triangles corresponding to vertices in $\Gamma^{(0)}_0 \setminus \Gamma^{(0)}_{k+1}$. Indeed, we do not have modified $(\zeta'_{\lambda'_0})^{(k)}$ on the triangles associated with the vertices in $\Gamma^{(0)}_0 \setminus \Gamma^{(0)}_k$ and we have chosen $\alpha_n$ in order to have the same central load on $T'_s$ too. The only thing we need to check is that the invariants on the couples of edges corresponding to the elements of $\Gamma^{(1)}_0 \setminus \Gamma^{(1)}_{k+1}$ are correct. The edge $\nu'_n$ is the result of the fusion in $\mu'$ of the edge $\mu'_n$ of $T'_s$ and of another edge $\mu'_u$ of a certain triangle $T'_v$. We need to prove that the invariant $y_u$ of $(\zeta'_{\mu'})^{(k+1)}$ on the edge $\mu'_u$ is also correct. But this easily follows from the fact that, labelled as $\bar{y}_n, \bar{y}_u$ the invariants of $\zeta_\mu \circ \Phi^q_{\mu\mu'}$ on $\mu_n$ and $\mu_u$, the products $y_n y_u$ and $\bar{y}_n \bar{y}_u$ are equal because $\rho'_{\lambda'}$ is isomorphic to $\rho_\lambda \circ \Phi^q_{\lambda\lambda'}$ and, thanks to the choice made for $\alpha_n$, we also have $y_n = \bar{y}_n$. With this last observation we conclude the proof of the surjectivity of $f$.

Now we have to deal with the injectivity. Fix $[\zeta_{\lambda_0}, \zeta'_{\lambda'_0}]$, $[\bar{\zeta}_{\lambda_0}, \bar{\zeta}'_{\lambda'_0}]$ two elements of $\mathscr{A}^{\rho\rho'}_{\lambda\lambda'}$ and suppose that their images $[\zeta_\nu, \zeta'_{\nu'}]$, $[\bar{\zeta}_\nu, \bar{\zeta}'_{\nu'}]$ in $\mathscr{B}^{\rho\rho'}_{\lambda\lambda'}$ coincide. As usual we are denoting by $\zeta_\nu, \bar{\zeta}_\nu$ the local representations on $\mathcal{T}^q_\nu$ represented by $\zeta_{\lambda_0}, \bar{\zeta}_{\lambda_0}$ and analogously for $\zeta'_{\nu'}, \bar{\zeta}'_{\nu'}$. Because $[\zeta_\nu, \zeta'_{\nu'}] = [\bar{\zeta}_\nu, \bar{\zeta}'_{\nu'}]$, we have

$$\zeta_{\lambda_0} \xrightarrow{\alpha_j} \bar{\zeta}_{\lambda_0}$$

$$\zeta'_{\lambda'_0} \xrightarrow{\beta_j} \bar{\zeta}'_{\lambda'_0}$$

with $\alpha_j = \beta_j$ for every $j$ such that $\lambda^*_j \in \Gamma_{S,\lambda} \setminus \Gamma_{R',\nu'}$. We must show that $\alpha_j = \beta_j$ for every $j \neq i$. Similarly to what done before, we take $\Gamma_0$ the sub-tree of $\Gamma_{R',\nu'}$ related to a connected component of $R'$, and we prove that on all the edges $(\nu'_j)^*$ of $\Gamma_0$, except for $(\lambda'_i)^* = (\nu'_i)^*$ possibly, we have $\alpha_j = \beta_j$. Select a leaf of $\Gamma_0$, with vertex $T'_s$ and edge $(\nu'_n)^*$ as before. If $n = i$ then we look for another leaf: if it exists we replace $\nu'_n$ with it in the following procedure; if it does not, then this component is dual of the only square $Q$, hence we can skip it and focus on a different connected component. Assume that $n \neq i$. Because $[\zeta_{\lambda_0}, \zeta'_{\lambda'_0}]$ belongs to $\mathscr{A}^{\rho\rho'}_{\lambda\lambda'}$, the representation $\zeta'_{\mu'}$ is isomorphic to $\zeta_\mu \circ \Phi^q_{\mu\mu'}$ and analogously $\bar{\zeta}'_{\mu'}$ is isomorphic to $\bar{\zeta}_\mu \circ \Phi^q_{\mu\mu'}$. In particular, because $T'_s$ is not contained in $Q$ and thanks to the Fusion property of $\Phi^q_{\mu\mu'}$, the invariants of the edges and the central load on $T'_s$ of $\zeta'_{\lambda'_0}$ are the same of those of $\zeta_{\lambda_0}$, and analogously the invariants of the edges and the central load on $T'_s$ of $\bar\zeta'_{\lambda'_0}$ are the same of those of $\bar\zeta_{\lambda_0}$. Denoting by $\nu'_l, \nu'_m, \nu'_n$ the edges of $T'_s$, we already know that $\alpha_l = \beta_l$ and $\alpha_m = \beta_m$ because $(T'_h)^*$ is a leaf. Now we see that, if $\alpha_n \neq \beta_n$, then it can not happen in the same moment that $\zeta_\mu \circ \Phi^q_{\mu\mu'}$ is isomorphic to $\zeta'_{\mu'}$ and $\bar\zeta_\mu \circ \Phi^q_{\mu\mu'}$ is isomorphic to $\bar\zeta'_{\mu'}$, by inspection of the invariants. Indeed, if $\alpha_n^N \neq \beta_n^N$, then the invariants on the edge $\nu'_n$ can not be equal in both cases. If $\alpha_n^N = \beta_n^N$ but $\alpha_n \neq \beta_n$, then the invariants on the edges coincide in both cases, but not the central loads of $T'_s$. This concludes the proof of the main part of the Lemma.

To see that $\bar p \circ f = \widetilde p$, it is sufficient to observe that, if $L$ is an isomorphism between $\zeta_\mu \circ \Phi_{\mu\mu'}$ and $\zeta'_{\mu'}$, then $L$ is an isomorphism between the fusions $\zeta_\nu \circ \Phi^q_{\nu\nu'}$ and $\zeta'_{\nu'}$ too. $\square$

If $\zeta_\nu, \bar\zeta_\nu$ are two elements of $\mathscr{F}_{R'}(\rho_\lambda)$, then we can define a notion of transition constants like in the case of $R' = S_0$. Indeed, taken two representatives $\zeta_{\lambda_0}$ and $\bar\zeta_{\lambda_0}$ of $\zeta_\nu$ and $\bar\zeta_\nu$, then for every $\lambda_j$ in $\lambda$ there exists a constant $\alpha_j$ such that $\zeta_{\lambda_0} \xrightarrow{\alpha_j} \bar\zeta_{\lambda_0}$. The constants $\alpha_j$ depend in general on the chosen representatives, but only those $\alpha_j$ such that $\lambda_j^*$ is an edge of $\Gamma_{R',\nu} \subset \Gamma_{S,\lambda}$. Therefore we can define the transition constants from $\zeta_\nu$ and $\bar\zeta_\nu$ as the collection of the $\alpha_j$ that correspond to edges in $\Gamma_{S,\lambda} \setminus \Gamma_{R',\nu}$, fixed a certain couple of representatives. We will briefly write this as
$$\zeta_\nu \xrightarrow[R']{\alpha_j} \bar\zeta_\nu$$

Moreover, given $\zeta_{\lambda_0}$ a representative of $\zeta_\nu$ and such a collection of transitions constants, we can extend it arbitrarily to a set with one $\alpha_j$ for each $\lambda_j$ and define a new representation $\bar\zeta_{\lambda_0}$ of $\mathcal{T}^q_{\lambda_0}$ by taking $\zeta_{\lambda_0} \xrightarrow{\alpha_j} \bar\zeta_{\lambda_0}$. It is immediate to see that the local representation $\bar\zeta_\nu$ on $\mathcal{T}^q_\nu$ represented by $\bar\zeta_{\lambda_0}$ does not depend on the way we extended the set $(\alpha_j)_j$. In conclusion, given $\zeta_\nu \in \mathscr{F}_{R'}(\rho_\lambda)$ and a set of $\alpha_j \in \mathbb{C}^*$, one for each edge of $\Gamma_{S,\lambda} \setminus \Gamma_{R',\nu}$, there is a unique local representation $\bar\zeta_\nu \in \mathscr{F}_{R'}(\rho_\lambda)$ such that $\zeta_\nu \xrightarrow[R']{\alpha_j} \bar\zeta_\nu$.

### 3.2.2 The other cases

It is not difficult to deduce a Lemma for the case in which $\lambda = \alpha(\lambda')$, with $\alpha \in \mathfrak{S}_n$, analogous to Lemma 3.2. For the sake of simplicity, we will deal with the case $\lambda = \lambda'$, but the same holds in the case of a generic reindexing. Fix a surface $R'$, disjoint union of ideal polygons, such that $S$ is obtained from $R'$ by fusion and denote by $\nu \in \Lambda(R')$ its ideal triangulation corresponding to $\lambda \in \Lambda(S)$.

Now define $B^{\rho\rho'}_{\lambda\lambda}$ as the set of couples $(\zeta_\nu, \zeta'_\nu)$ in the product $\mathscr{F}_{R'}(\rho_\lambda) \times \mathscr{F}_{R'}(\rho'_\lambda)$ such that $\zeta_\mu$ and $\zeta'_\mu$ are isomorphic. The equivalence relation on $B^{\rho\rho'}_{\lambda\lambda}$ leading to $\mathscr{B}^{\rho\rho'}_{\lambda\lambda}$ is defined as follows: $(\zeta_\nu, \zeta'_\nu)$ is $\approx_{R'}$-equivalent to $(\bar\zeta_\nu, \bar\zeta'_\nu)$ if
$$\zeta_\nu \xrightarrow[R']{\alpha_j} \bar\zeta_\nu$$
$$\zeta'_\nu \xrightarrow[R']{\beta_j} \bar\zeta'_\nu$$

with $\alpha_j = \beta_j$ for every $j$ such that $\lambda_j^* \in \Gamma_{S,\lambda} \setminus \Gamma_{R',\nu}$. As in the previous case, we have natural maps $f\colon \mathscr{A}_{\lambda\lambda}^{\rho\rho'} \to \mathscr{B}_{\lambda\lambda}^{\rho\rho'}$ and $\bar{p}\colon \mathscr{B}_{\lambda\lambda}^{\rho\rho'} \to \mathrm{Hom}(V_\lambda', V_\lambda)/\doteq$, defined in the same way. With the same procedure as in the proof of Lemma 3.2 the following fact can be shown:

**Lemma 3.3.** In the above notations, the following hold:

- the map $f\colon \mathscr{A}_{\lambda\lambda}^{\rho\rho'} \to \mathscr{B}_{\lambda\lambda}^{\rho\rho'}$ is a bijection;

- the following diagram is commutative

$$\begin{array}{ccc} \mathscr{A}_{\lambda\lambda}^{\rho\rho'} & \xrightarrow{f} & \mathscr{B}_{\lambda\lambda}^{\rho\rho'} \\ {\scriptstyle \widetilde{p}} \searrow & & \swarrow {\scriptstyle \bar{p}} \\ & \mathrm{Hom}(V_\lambda', V_\lambda)/\doteq & \end{array}$$

In particular $\mathscr{L}_{\lambda\lambda'}^{\rho\rho'} = \mathrm{Im}\,\widetilde{p} = \mathrm{Im}\,\bar{p}$.

Let $\lambda, \lambda' \in \Lambda(S)$ be two ideal triangulations differing by a diagonal exchange or a reindexing. Then we can describe an action of $H_1(S;\mathbb{Z}_N)$ on the set $\mathscr{B}_{\lambda\lambda'}^{\rho\rho'}$ just introduced. Given $c \in H_1(S;\mathbb{Z}_N)$ and $[\zeta_\nu, \zeta'_{\nu'}] \in \mathscr{B}_{\lambda\lambda'}^{\rho\rho'}$, we take a representative $\zeta'_{\lambda_0'}$ of $\zeta'_{\nu'}$ and we define

$$c \cdot [\zeta_\nu, \zeta'_{\nu'}] := [\zeta_\nu, c \cdot \zeta'_{\nu'}]$$

where we are denoting by $c \cdot \zeta'_{\nu'}$ the local representation on $(R', \nu')$ induced by $c \cdot \zeta'_{\lambda_0'}$. It is clear that this action corresponds via $f$ to the usual action of $H_1(S;\mathbb{Z}_N)$ on $\mathscr{A}_{\lambda\lambda'}^{\rho\rho'}$ hence, thanks to Lemmas 3.2 and 3.3, the properties of the action that we observed on $\mathscr{A}_{\lambda\lambda'}^{\rho\rho'}$, like transitivity and freeness, hold also on $\mathscr{B}_{\lambda\lambda'}^{\rho\rho'}$.

## 3.3 Elementary Composition property

Let $\lambda$ be an ideal triangulation of a surface $S$ and assume that there exist in $S$ three triangles composing a pentagon with diagonals $\lambda_i$ and $\lambda_j$, possibly not embedded in $S$. We enumerate the sequence of triangulations appearing in the Pentagon relation as follows:

$$\lambda^{(0)} := \lambda$$
$$\lambda^{(1)} := \Delta_i(\lambda)$$
$$\vdots$$
$$\lambda^{(4)} := (\Delta_j \circ \Delta_i \circ \Delta_j \circ \Delta_i)(\lambda)$$
$$\lambda^{(5)} := \alpha_{ij}(\lambda)$$
$$\lambda^{(6)} := \lambda$$

Label as $R$ the surface obtained from $S$ by splitting it along all the edges except for $\lambda_i$ and $\lambda_j$. Then $R$ is the disjoint union of an ideal pentagon $P$ and $m-3$ triangles $T_4, \ldots, T_m$. The ideal triangulations $\lambda^{(k)}$ lift to a sequence $\mu^{(k)}$ of triangulations on $R$, which are related by diagonal exchanges along $\mu_i$ and $\mu_j$, the diagonals of $P$. Each ideal triangulation $\mu^{(k)}$ can be naturally presented as the disjoint union of a triangulation $\mu_P^{(k)}$ of the pentagon and the only possible triangulation $\bar{\mu}$ on $T_4 \sqcup \cdots \sqcup T_m$.

By definition of the sets $\mathscr{L}_{\lambda^{(k)}\lambda^{(k+1)}}^{\rho\rho}$, it makes sense to compose the intertwining operators as follows

$$\begin{array}{rcl} \prod_{k=0}^{5} \mathscr{L}_{\lambda^{(k)}\lambda^{(k+1)}}^{\rho\rho} & \longrightarrow & \mathrm{End}(V_\lambda) \\ (L_0, \ldots, L_5) & \longmapsto & L_0 \circ \cdots \circ L_5 \end{array}$$

**Lemma 3.4.** Let $\rho = \{\rho_\lambda \colon \mathcal{T}_\lambda^q \to \mathrm{End}(V_\lambda)\}_{\lambda \in \Lambda(S)}$ be a local representation of $\mathcal{T}_S^q$ and $\lambda^{(k)}$ a sequence of triangulations as described above. Then the composition

$$\begin{array}{rcl} \prod_{k=0}^{5} \mathscr{L}_{\lambda^{(k)}\lambda^{(k+1)}}^{\rho\rho} & \longrightarrow & \mathscr{L}_{\lambda\lambda}^{\rho\rho} \\ (L_0, \ldots, L_5) & \longmapsto & L_0 \circ \cdots \circ L_5 \end{array}$$

is well defined and it verifies

$$(c_0 \cdot L_0) \circ \cdots \circ (c_5 \cdot L_5) = \left(\sum_{k=0}^{5} c_k\right) \cdot (L_0 \circ \cdots \circ L_5)$$

*Proof.* By virtue of Lemma 3.2, each $\mathscr{L}_{\lambda^{(k)}\lambda^{(k+1)}}^{\rho\rho}$ is in bijection with the set $\mathscr{B}_{\lambda^{(k)}\lambda^{(k+1)}}^{\rho\rho}$, defined as

$$\left\{(\zeta^{(k)}, \zeta^{(k+1)}) \in \mathscr{F}_R(\rho_{\lambda^{(k)}}) \times \mathscr{F}_R(\rho_{\lambda^{(k+1)}}) \mid \zeta^{(k)} \circ \Phi_{\mu^{(k)}\mu^{(k+1)}}^q \text{ isom to } \zeta^{(k+1)}\right\}\Big/_{\approx}$$

where $R$ is the surface described above, which is a disjoint union of an ideal pentagon and triangles, in particular a disjoint union of polygons. The first step of the proof will be the following: we want to translate the composition map on the $\mathscr{L}_{\lambda^{(k)}\lambda^{(k+1)}}^{\rho\rho}$ in an application defined on the product of the sets $\mathscr{B}_{\lambda^{(k)}\lambda^{(k+1)}}^{\rho\rho}$, in order to have a better control of the behaviour of the actions of $H_1(S; \mathbb{Z}_N)$. All the efforts spent to prove Lemma 3.2 allow us to manage local representations $\zeta^{(k)}$ defined on Chekhov-Fock algebras of the same surface $R$ associated with the ideal triangulations $\mu^{(k)}$.

We will denote an element of $\prod_{k=0}^{5} \mathscr{B}_{\lambda^{(k)}\lambda^{(k+1)}}^{\rho\rho}$ by $([\zeta_0^{(k)}, \zeta_1^{(k)}])_{k=0}^{5}$, where $\zeta_0^{(k)}$ is a local representation of $\mathcal{T}_{\mu^{(k)}}^q$ that represents $\rho_{\lambda^{(k)}}$ on $S$ and $\zeta_1^{(k)}$ is a local representation of $\mathcal{T}_{\mu^{(k+1)}}^q$ that represents $\rho_{\lambda^{(k+1)}}$ on $S$. Suppose that the 6-tuple corresponding to $([\zeta_0^{(k)}, \zeta_1^{(k)}])_{k=0}^{5}$ in $\prod_{k=0}^{5} \mathscr{L}_{\lambda^{(k)}\lambda^{(k+1)}}^{\rho\rho}$ is $(L_0, \ldots, L_5)$. We would like to understand if there exists an element in $\mathscr{B}_{\lambda\lambda}^{\rho\rho}$ corresponding to $L_0 \circ \cdots \circ L_5$ and how can be described. The elements $[\zeta_0^{(k)}, \zeta_1^{(k)}]$ belong to $\mathscr{B}_{\lambda^{(k)}\lambda^{(k+1)}}^{\rho\rho}$, so for every $k \in \{0, \ldots, 5\}$ and $i \in \{0, 1\}$ the local representation $\zeta_i^{(k)}$ represents $\rho_{\lambda^{(k+i)}}$, i. e. it is an element of $\mathscr{F}_R(\rho_{\lambda^{(k+i)}})$. Then, for every $k = 1, \ldots, 5$ there exist transition constants $\alpha_h^{(k)}$ such that

$$\zeta_0^{(k)} \xrightarrow[R]{\alpha_h^{(k)}} \zeta_1^{(k-1)}$$

where $\alpha_h^{(k)}$ is associated with an edge $\lambda_h^* \in \Gamma_{S,\lambda^{(k)}} \setminus \Gamma_{R,\mu^{(k)}}$. Observe that there is a natural bijection

$$\Gamma_{S,\lambda^{(k)}} \setminus \Gamma_{R,\mu^{(k)}} \longleftrightarrow \Gamma_{S,\lambda^{(k+1)}} \setminus \Gamma_{R,\mu^{(k+1)}}$$

for every $k = 0, \ldots, 5$. Indeed, $\lambda^{(k)}$ and $\lambda^{(k+1)}$ differ by a diagonal exchange, so we have a canonical correspondence between all the edges of them except for the ones on which we do diagonal exchange; in particular on all the edges different from $\lambda_i^{(k)}$ and $\lambda_j^{(k)}$ that compose $\Gamma_{R,\mu^{(k)}}$.

We can change representative of $[\zeta_0^{(1)}, \zeta_1^{(1)}]$ by taking $(\zeta_1^{(0)}, \bar{\zeta}_1^{(1)})$, where

$$\zeta_0^{(1)} \xrightarrow[R]{\alpha_h^{(1)}} \zeta_1^{(0)}$$

$$\zeta_1^{(1)} \xrightarrow[R]{\alpha_h^{(1)}} \bar{\zeta}_1^{(1)}$$

Analogously we can construct a representative $(\bar{\zeta}_1^{(1)}, \bar{\zeta}_1^{(2)})$ of $[\zeta_0^{(2)}, \zeta_1^{(2)}]$ defined by

$$\zeta_0^{(2)} \xrightarrow[R]{\alpha_h^{(2)}} \zeta_1^{(1)} \xrightarrow[R]{\alpha_h^{(1)}} \bar{\zeta}_1^{(1)}$$

$$\zeta_1^{(2)} \xrightarrow[R]{\alpha_h^{(2)}} \widetilde{\zeta}_1^{(2)} \xrightarrow[R]{\alpha_h^{(1)}} \bar{\zeta}_1^{(2)}$$

In the same way, for every $k = 3, 4, 5$ we choose the representative $(\bar{\zeta}_1^{(k-1)}, \bar{\zeta}_1^{(k)})$ of $[\zeta_0^{(k)}, \zeta_1^{(k)}]$ obtained as follows

$$\zeta_0^{(k)} \xrightarrow[R]{\alpha_h^{(k)}} \zeta_1^{(k-1)} \xrightarrow[R]{\alpha_h^{(k-1)}} \cdots \xrightarrow[R]{\alpha_h^{(1)}} \bar{\zeta}_1^{(k-1)}$$

$$\zeta_1^{(k)} \xrightarrow[R]{\alpha_h^{(k)}} \widetilde{\zeta}_1^{(k)} \xrightarrow[R]{\alpha_h^{(k-1)}} \cdots \xrightarrow[R]{\alpha_h^{(1)}} \bar{\zeta}_1^{(k)}$$

In this way we obtain a 6-tuple of representatives of $([\zeta_0^{(k)}, \zeta_1^{(k)}])_{k=0}^{5}$ that looks like

$$([\zeta_0^{(0)}, \zeta_1^{(0)}], [\zeta_1^{(0)}, \bar{\zeta}_1^{(1)}], [\bar{\zeta}_1^{(1)}, \bar{\zeta}_1^{(2)}], [\bar{\zeta}_1^{(2)}, \bar{\zeta}_1^{(3)}], [\bar{\zeta}_1^{(3)}, \bar{\zeta}_1^{(4)}], [\bar{\zeta}_1^{(4)}, \bar{\zeta}_1^{(5)}])$$

Despite their terrible appearance, these representatives have the good property that the elements appearing that belong to $\mathscr{F}_R(\rho_{\lambda^{(k)}})$ are both equal to $\bar{\zeta}_1^{(k-1)}$ for $k = 2, 3, 4$, or to $\zeta_1^{(0)}$ when $k = 1$. We claim that the local representations $\zeta_0^{(0)}$ and $\bar{\zeta}_1^{(5)}$, both elements of $\mathscr{F}_R(\rho_\lambda)$, are isomorphic via $L_0 \circ \cdots \circ L_5$, so $[\zeta_0^{(0)}, \bar{\zeta}_1^{(5)}]$ is an element of $\mathscr{B}_{\lambda\lambda}^{\rho\rho}$ whose image is $L_0 \circ \cdots \circ L_5$.

Because we have chosen representatives of the classes $[\zeta_0^{(i)}, \zeta_1^{(i)}]$ corresponding to the linear isomorphisms $L_i$, we have that

$$(\zeta_0^{(0)} \circ \Phi_{\mu^{(0)}\mu^{(1)}}^q)(X^{(1)}) = L_0 \circ \zeta_1^{(0)}(X^{(1)}) \circ L_0^{-1} \qquad \forall X^{(1)} \in \mathcal{T}_{\mu^{(1)}}^q$$

$$(\zeta_1^{(0)} \circ \Phi_{\mu^{(1)}\mu^{(2)}}^q)(X^{(2)}) = L_1 \circ \bar{\zeta}_1^{(1)}(X^{(2)}) \circ L_1^{-1} \qquad \forall X^{(2)} \in \mathcal{T}_{\mu^{(1)}}^q$$

$$(\bar{\zeta}_1^{(k-1)} \circ \Phi_{\mu^{(k)}\mu^{(k+1)}}^q)(X^{(k+1)}) = L_k \circ \bar{\zeta}_1^{(k)}(X^{(k+1)}) \circ L_k^{-1} \quad \forall X^{(k+1)} \in \mathcal{T}_{\mu^{(k+1)}}^q$$

for $k = 2, \ldots, 5$. Hence we deduce that, for every $X \in \mathcal{T}_\lambda^q$

$$\begin{aligned}
\zeta_0^{(0)}(X) &= (\zeta_0^{(0)} \circ \Phi^q_{\mu^{(0)}\mu^{(6)}})(X) \\
&= (\zeta_0^{(0)} \circ \Phi^q_{\mu^{(0)}\mu^{(1)}} \circ \cdots \circ \Phi^q_{\mu^{(5)}\mu^{(6)}})(X) \\
&= L_0 \circ (\zeta_1^{(0)} \circ \Phi^q_{\mu^{(1)}\mu^{(2)}} \circ \cdots \circ \Phi^q_{\mu^{(5)}\mu^{(6)}})(X) \circ L_0^{-1} \\
&\vdots \\
&= (L_0 \circ \cdots \circ L_5) \circ \bar\zeta_1^{(5)}(X) \circ (L_0 \circ \cdots \circ L_5)^{-1}
\end{aligned}$$

where we have used the Pentagon relation of $\Phi^q_{\lambda\lambda'}$. This proves the claim.

Denote by $\odot$ the operation from $\prod_{k=0}^{5} \mathscr{B}^{\rho\rho}_{\lambda^{(k)}\lambda^{(k+1)}}$ to $\mathscr{B}^{\rho\rho}_{\lambda\lambda}$ corresponding to the composition of $\mathscr{L}^{\rho\rho}_{\lambda^{(k)}\lambda^{(k+1)}}$ just described. Then we have proved the following equality

$$[\zeta_0^{(0)}, \zeta_1^{(0)}] \odot [\zeta_1^{(0)}, \bar\zeta_1^{(1)}] \odot [\bar\zeta_1^{(1)}, \bar\zeta_1^{(2)}] \odot \cdots \odot [\bar\zeta_1^{(4)}, \bar\zeta_1^{(5)}] = [\zeta_0^{(0)}, \bar\zeta_1^{(5)}] \qquad (17)$$

Observe that we do not have to verify any dependence on the chosen representatives, because the composition map is obviously well defined and we have just rewrite it on the sets $\mathscr{B}^{\rho\rho}_{\lambda^{(k)}\lambda^{(k+1)}}$.

Now the conclusion is just a verification. Define

$$s_k := \sum_{i=0}^{k} c_i \in H_1(S; \mathbb{Z}_N)$$

for every $k = 0, \ldots, 5$. If $\bar p \colon \mathscr{B}^{\rho\rho}_{\lambda\lambda} \to \mathscr{L}^{\rho\rho}_{\lambda\lambda}$ is the bijection described in Lemma 3.3, we observe that

$$\begin{aligned}
\bar p^{-1}((c_0 \cdot L_0) \circ \cdots \circ (c_5 \cdot L_5)) &= \\
&= (c_0 \cdot [\zeta_0^{(0)}, \zeta_1^{(0)}]) \odot (c_1 \cdot [\zeta_1^{(0)}, \bar\zeta_1^{(1)}]) \odot \cdots \odot (c_5 \cdot [\bar\zeta_1^{(4)}, \bar\zeta_1^{(5)}]) \\
&= [\zeta_0^{(0)}, c_0 \cdot \zeta_1^{(0)}] \odot [\zeta_1^{(0)}, c_1 \cdot \bar\zeta_1^{(1)}] \odot \cdots \odot [\bar\zeta_1^{(4)}, c_5 \cdot \bar\zeta_1^{(5)}] \\
&= [\zeta_0^{(0)}, s_0 \cdot \zeta_1^{(0)}] \odot [\zeta_1^{(0)}, (s_1 - s_0) \cdot \bar\zeta_1^{(1)}] \odot \cdots \odot [\bar\zeta_1^{(4)}, (s_5 - s_4) \cdot \bar\zeta_1^{(5)}] \\
&= [\zeta_0^{(0)}, s_0 \cdot \zeta_1^{(0)}] \odot [s_0 \cdot \zeta_1^{(0)}, s_1 \cdot \bar\zeta_1^{(1)}] \odot \cdots \odot [s_4 \cdot \bar\zeta_1^{(4)}, s_5 \cdot \bar\zeta_1^{(5)}] \\
&= [\zeta_0^{(0)}, s_5 \cdot \bar\zeta_1^{(5)}] \\
&= s_5 \cdot [\zeta_0^{(0)}, \bar\zeta_1^{(5)}] \\
&= s_5 \cdot ([\zeta_0^{(0)}, \zeta_1^{(0)}] \odot [\zeta_1^{(0)}, \bar\zeta_1^{(1)}] \odot \cdots \odot [\bar\zeta_1^{(4)}, \bar\zeta_1^{(5)}]) \\
&= s_5 \cdot \bar p^{-1}(L_0 \circ \cdots \circ L_5) \\
&= \bar p^{-1}(s_5 \cdot (L_0 \circ \cdots \circ L_5))
\end{aligned}$$

where we are using the relation 17 and the equality $[\zeta, \zeta'] = [c \cdot \zeta, c \cdot \zeta'] \in \mathscr{B}^{\rho\rho}_{\lambda^{(k)}\lambda^{(k+1)}}$. Finally, by applying $\bar p$ to the first and the last members we obtain

$$(c_0 \cdot L_0) \circ \cdots \circ (c_5 \cdot L_5) = s_5 \cdot (L_0 \circ \cdots \circ L_5) = \left(\sum_{k=0}^{5} c_k\right) \cdot (L_0 \circ \cdots \circ L_5)$$

as desired. $\square$

In the same way, we can prove analogous results with respect to the other relations between the elementary operations on the ideal triangulations. We limit ourselves to the enunciations of these properties, their proof can be obtained with procedures analogous to the Pentagon relation case, by considering the surface $R$, result of the splitting of $S$ along all the edges except for the ones along we do diagonal exchange, if there is any.

**Lemma 3.5.** Let $\rho = \{\rho_\lambda \colon \mathcal{T}_\lambda^q \to \mathrm{End}(V_\lambda)\}_{\lambda \in \Lambda(S)}$ be a local representation of $\mathcal{T}_S^q$ and let $\lambda \in \Lambda(S)$ be an ideal triangulation.

COMPOSITION RELATION: given $\alpha, \beta \in \mathfrak{S}_n$, consider the following path of ideal triangulations
$$\lambda^{(0)} := \lambda, \quad \lambda^{(1)} := \alpha(\lambda), \quad \lambda^{(2)} := (\alpha \circ \beta)(\lambda)$$

Then the composition
$$\begin{array}{rcl} \mathscr{L}^{\rho\rho}_{\lambda^{(0)}\lambda^{(1)}} \times \mathscr{L}^{\rho\rho}_{\lambda^{(1)}\lambda^{(2)}} & \longrightarrow & \mathscr{L}^{\rho\rho}_{\lambda^{(0)}\lambda^{(2)}} \\ (L_0, L_1) & \longmapsto & L_0 \circ L_1 \end{array}$$

is well defined and it verifies
$$(c_0 \cdot L_0) \circ (c_1 \cdot L_1) = (c_0 + c_1) \cdot (L_0 \circ L_1)$$

for every $c_i \in H_1(S; \mathbb{Z}_N)$;

REFLEXIVITY RELATION: given $\lambda_i$ a diagonal of a certain square in $\lambda$, consider the following path of ideal triangulations
$$\lambda^{(0)} := \lambda, \quad \lambda^{(1)} := \Delta_i(\lambda), \quad \lambda^{(2)} := \lambda$$

Then the composition
$$\begin{array}{rcl} \mathscr{L}^{\rho\rho}_{\lambda^{(0)}\lambda^{(1)}} \times \mathscr{L}^{\rho\rho}_{\lambda^{(1)}\lambda^{(2)}} & \longrightarrow & \mathscr{L}^{\rho\rho}_{\lambda\lambda} \\ (L_0, L_1) & \longmapsto & L_0 \circ L_1 \end{array}$$

is well defined and it verifies
$$(c_0 \cdot L_0) \circ (c_1 \cdot L_1) = (c_0 + c_1) \cdot (L_0 \circ L_1)$$

for every $c_i \in H_1(S; \mathbb{Z}_N)$;

RE-INDEXING RELATION: given $\lambda_i$ a diagonal of a certain square in $\lambda$ and $\alpha \in \mathfrak{S}_n$, consider the following path of ideal triangulations
$$\begin{aligned} \lambda^{(0)} &:= \lambda \\ \lambda^{(1)} &:= \alpha(\lambda) \\ \lambda^{(2)} &:= \Delta_i(\alpha(\lambda)) = \alpha(\Delta_{\alpha(i)}(\lambda)) \\ \lambda^{(3)} &:= \Delta_{\alpha(i)}(\lambda) \\ \lambda^{(4)} &:= \lambda \end{aligned}$$

Then the composition
$$\begin{array}{rcl} \prod_{i=0}^{3} \mathscr{L}^{\rho\rho}_{\lambda^{(k)}\lambda^{(k+1)}} & \longrightarrow & \mathscr{L}^{\rho\rho}_{\lambda\lambda} \\ (L_0, L_1, L_2, L_3) & \longmapsto & L_0 \circ L_1 \circ L_2 \circ L_3 \end{array}$$

is well defined and it verifies

$$(c_0 \cdot L_0) \circ (c_1 \cdot L_1) \circ (c_2 \cdot L_2) \circ (c_3 \cdot L_3) = \left(\sum_{k=0}^{3} c_k\right) \cdot (L_0 \circ L_1 \circ L_2 \circ L_3)$$

for every $c_i \in H_1(S; \mathbb{Z}_N)$. The same holds for the inverse sequence $\bar{\lambda}^{(i)}$, with $\bar{\lambda}^{(i)} := \lambda^{(4-i)}$ for $i = 0, \ldots, 4$;

DISTANT COMMUTATIVITY RELATION: given $\lambda_i$ and $\lambda_j$ diagonals in $\lambda$ that do not belong to a common triangle, consider the following path of ideal triangulations

$$\lambda^{(0)} := \lambda$$
$$\lambda^{(1)} := \Delta_i(\lambda)$$
$$\lambda^{(2)} := (\Delta_j \circ \Delta_i)(\lambda) = (\Delta_i \circ \Delta_j)(\lambda)$$
$$\lambda^{(3)} := \Delta_j(\lambda)$$
$$\lambda^{(4)} := \lambda$$

Then the composition

$$\begin{array}{ccc} \prod_{i=0}^{3} \mathscr{L}^{\rho\rho}_{\lambda^{(k)}\lambda^{(k+1)}} & \longrightarrow & \mathscr{L}^{\rho\rho}_{\lambda\lambda} \\ (L_0, L_1, L_2, L_3) & \longmapsto & L_0 \circ L_1 \circ L_2 \circ L_3 \end{array}$$

is well defined and it verifies

$$(c_0 \cdot L_0) \circ (c_1 \cdot L_1) \circ (c_2 \cdot L_2) \circ (c_3 \cdot L_3) = \left(\sum_{k=0}^{3} c_k\right) \cdot (L_0 \circ L_1 \circ L_2 \circ L_3)$$

for every $c_i \in H_1(S; \mathbb{Z}_N)$;

PENTAGON RELATION: given $\lambda_i$ and $\lambda_j$ diagonals of a common pentagon in $\lambda$, consider the following path of ideal triangulations

$$\lambda^{(0)} := \lambda$$
$$\lambda^{(1)} := \Delta_i(\lambda)$$
$$\vdots$$
$$\lambda^{(4)} := (\Delta_j \circ \Delta_i \circ \Delta_j \circ \Delta_i)(\lambda)$$
$$\lambda^{(5)} := \alpha_{ij}(\lambda)$$
$$\lambda^{(6)} := \lambda$$

Then the composition

$$\begin{array}{ccc} \prod_{k=0}^{5} \mathscr{L}^{\rho\rho}_{\lambda^{(k)}\lambda^{(k+1)}} & \longrightarrow & \mathscr{L}^{\rho\rho}_{\lambda\lambda} \\ (L_0, \ldots, L_5) & \longmapsto & L_0 \circ \cdots \circ L_5 \end{array}$$

is well defined and it verifies

$$(c_0 \cdot L_0) \circ \cdots \circ (c_5 \cdot L_5) = \left(\sum_{k=0}^{5} c_k\right) \cdot (L_0 \circ \cdots \circ L_5)$$

for every $c_i \in H_1(S; \mathbb{Z}_N)$.

The relations between the actions exposed in Lemma 3.5 and the transitivity of the actions proved in Theorem 2.1 imply that the compositions maps are surjective in every case exposed in Lemma 3.5.

**Lemma 3.6.** Let

$$\rho = \{\rho_\lambda \colon \mathcal{T}_\lambda^q \to \mathrm{End}(V_\lambda)\}_{\lambda \in \Lambda(S)}$$
$$\rho' = \{\rho'_\lambda \colon \mathcal{T}_\lambda^q \to \mathrm{End}(V'_\lambda)\}_{\lambda \in \Lambda(S)}$$
$$\rho'' = \{\rho''_\lambda \colon \mathcal{T}_\lambda^q \to \mathrm{End}(V''_\lambda)\}_{\lambda \in \Lambda(S)}$$

be three isomorphic local representations of $\mathcal{T}_S^q$, and let $\lambda, \lambda'$ be two ideal triangulations of $S$ that differ by a diagonal exchange or a re-indexing. Then the compositions

$$\circ \colon \begin{array}{ccc} \mathscr{L}_{\lambda\lambda'}^{\rho\rho'} \times \mathscr{L}_{\lambda'\lambda'}^{\rho'\rho''} & \longrightarrow & \mathscr{L}_{\lambda\lambda'}^{\rho\rho''} \\ (L, M) & \longmapsto & L \circ M \end{array}$$

$$\circ \colon \begin{array}{ccc} \mathscr{L}_{\lambda\lambda}^{\rho\rho'} \times \mathscr{L}_{\lambda\lambda'}^{\rho'\rho''} & \longrightarrow & \mathscr{L}_{\lambda\lambda'}^{\rho\rho''} \\ (L, M) & \longmapsto & L \circ M \end{array}$$

are well defined and they verify

$$(c \cdot L) \circ (d \cdot M) = (c + d) \cdot (L \circ M)$$

for every $c, d \in H_1(S; \mathbb{Z}_N)$.

*Proof.* Assume that the triangulations differ by diagonal exchange along the edge $\lambda_i$. Now take $R$ the surface obtained by splitting $S$ along all the edges except for $\lambda_i$. Then we can represent, as in the proof of Lemma 3.2, the compositions on the sets $\mathscr{B}_{\lambda\lambda'}^{\rho\rho'} \times \mathscr{B}_{\lambda'\lambda'}^{\rho'\rho''} \longrightarrow \mathscr{B}_{\lambda\lambda'}^{\rho\rho''}$ and $\mathscr{B}_{\lambda\lambda}^{\rho\rho'} \times \mathscr{B}_{\lambda\lambda'}^{\rho'\rho''} \longrightarrow \mathscr{B}_{\lambda\lambda'}^{\rho\rho''}$ respectively, where the set $\mathscr{B}$ are defined as quotients of sets of local representations on $R$. Now the proof can be achieved with the same ideas of Lemma 3.2. □

**Lemma 3.7.** Given $\lambda, \lambda' \in \Lambda(S)$ that differ by an elementary move, i. e. a diagonal exchange or a re-indexing of the edges, the map

$$(\cdot)^{-1} \colon \begin{array}{ccc} \mathscr{L}_{\lambda\lambda'}^{\rho\rho} & \longrightarrow & \mathscr{L}_{\lambda'\lambda}^{\rho\rho} \\ L & \longmapsto & L^{-1} \end{array}$$

verifies

$$(c \cdot L)^{-1} = (-c) \cdot L^{-1}$$

where $c$ is an element of $H_1(S; \mathbb{Z}_N)$, the action in the first member is on $\mathscr{L}_{\lambda\lambda'}^{\rho\rho}$ and the action in the second member is on $\mathscr{L}_{\lambda'\lambda}^{\rho\rho}$.

*Proof.* Thanks to Lemma 3.5, the proof is immediate:

$$((-c) \cdot L^{-1}) \circ (c \cdot L) = (-c + c) \cdot (L^{-1} \circ L) = id$$

□

## 4 The complete definition

In the previous paragraphs we have studied the elementary properties of the sets $\mathscr{L}_{\lambda\lambda'}^{\rho\rho'}$ endowed with certain actions $\psi_{\lambda\lambda'}^{\rho\rho'}$ of $H_1(S;\mathbb{Z}_N)$, but we have defined these objects only in the case in which $\lambda$ and $\lambda'$ differ by an elementary move. Now we have achieved all the elements to deal with the general construction.

Let $\lambda, \lambda' \in \Lambda(S)$ be two ideal triangulations of $S$. Thanks to Theorem 1.20 there exists a sequence of elementary moves that leads from $\lambda$ to $\lambda'$. Label the ideal triangulations we pass through as follows

$$\lambda = \lambda^{(0)}, \lambda^{(1)}, \ldots, \lambda^{(h)}, \lambda^{(h+1)} = \lambda'$$

Then we define $\mathscr{L}_{\lambda\lambda'}^{\rho\rho'}$ as the image of the composition map

$$\begin{array}{rcl} \prod_{k=0}^{h} \mathscr{L}_{\lambda^{(k)}\lambda^{(k+1)}}^{\rho\rho} \times \mathscr{L}_{\lambda'\lambda'}^{\rho\rho'} & \longrightarrow & \mathrm{Hom}(V_{\lambda'}', V_\lambda)/\doteq \\ (L_0, \ldots, L_h, L_{h+1}) & \longmapsto & L_0 \circ \cdots \circ L_{h+1} \end{array}$$

Moreover, denoting by $\pi \colon H_1(S;\mathbb{Z}_N)^{h+2} \to H_1(S;\mathbb{Z}_N)$ the map

$$\pi(c_0, \ldots, c_{h+1}) := \sum_{i=0}^{h+1} c_i$$

we can fix a section $s \colon H_1(S;\mathbb{Z}_N)^{h+2} \to H_1(S;\mathbb{Z}_N)$ of $\pi$, in other words a map that verifies $\pi \circ s = id$, not necessarily a homomorphism, and define the action of $c \in H_1(S;\mathbb{Z}_N)$ on an element $L = L_0 \circ \cdots \circ L_{h+1}$ in $\mathscr{L}_{\lambda\lambda'}^{\rho\rho}$ as follows

$$c \cdot L := (c_0 \cdot L_0) \circ \cdots \circ (c_{h+1} \cdot L_{h+1})$$

where $s(c) = (c_0, \ldots, c_{h+1})$ and we have chosen an element $(L_0, \ldots, L_{h+1})$ in the fibre of $L$ under the composition map.

In this definition of $(\mathscr{L}_{\lambda\lambda'}^{\rho\rho'}, \psi_{\lambda\lambda'}^{\rho\rho'})$ we have done some arbitrary choices:

- the path of triangulations between $\lambda$ and $\lambda'$;
- the section $s$ of $\pi$;
- the decomposition of $L \in \mathscr{L}_{\lambda\lambda'}^{\rho\rho'}$ as image under the composition map of a certain $(h+2)$-tuple $(L_0, \ldots, L_{h+1})$.

We need to prove that the object $(\mathscr{L}_{\lambda\lambda'}^{\rho\rho'}, \psi_{\lambda\lambda'}^{\rho\rho'})$ does not depend on the choices made. We start from the last one: take $(L_i)_i$ and $(L_i')_i$ such that

$$L_0 \circ \cdots \circ L_{h+1} = L_0' \circ \cdots \circ L_{h+1}'$$

We want to show that, for every $(c_0, \ldots, c_{h+1}) \in H_1(S;\mathbb{Z}_N)^{h+2}$, the following holds

$$(c_0 \cdot L_0) \circ \cdots \circ (c_{h+1} \cdot L_{h+1}) = (c_0 \cdot L_0') \circ \cdots \circ (c_{h+1} \cdot L_{h+1}')$$

Firstly observe that, because of the transitivity of the action of $H_1(S;\mathbb{Z}_N)$, for every $i = 0, \ldots, h+1$ there exists a $d_i \in H_1(S;\mathbb{Z}_N)$ such that $d_i \cdot L_i = L_i'$. By hypothesis, we have

$$(L_0 \circ \cdots \circ L_{h+1}) \circ (L_0' \circ \cdots \circ L_{h+1}')^{-1} = id \in \mathrm{End}(V_\lambda)$$

On the other hand

$$
\begin{aligned}
(L_0 \circ \cdots \circ L_{h+1}) &\circ (L'_0 \circ \cdots \circ L'_{h+1})^{-1} = \\
&= L_0 \circ \cdots \circ L_{h+1} \circ (d_{h+1} \cdot L_{h+1})^{-1} \circ \cdots \circ (d_0 \cdot L_0)^{-1} \\
&= L_0 \circ \cdots \circ L_{h+1} \circ ((-d_{h+1}) \cdot L_{h+1}^{-1}) \circ \cdots \circ ((-d_0) \cdot L_0^{-1}) \quad \text{Lemma 3.7} \\
&= L_0 \circ \cdots \circ L_h \circ ((-d_{h+1}) \cdot id) \circ ((-d_h) \cdot L_h^{-1}) \circ \cdots \circ ((-d_0) \cdot L_0^{-1}) \\
&\qquad\qquad\qquad\qquad\qquad\qquad\qquad\qquad\qquad\qquad\qquad \text{Lemma 3.5} \\
&= L_0 \circ \cdots \circ L_h \circ ((-d_{h+1} - d_h) \cdot L_h^{-1}) \circ \cdots \circ ((-d_0) \cdot L_0^{-1}) \quad \text{Lemma 3.6} \\
&\;\;\vdots \\
&= \left(-\sum_{i=0}^{h+1} d_i\right) \cdot id
\end{aligned}
$$

We observed that the action of $H_1(S; \mathbb{Z}_N)$ on $\mathscr{L}^{\rho\rho}_{\lambda\lambda}$ is free, so $\sum_{i=0}^{h+1} d_i$ must be equal to 0. It is simple to see that the steps of the relation above prove also that the following relation holds

$$
\begin{aligned}
((c_0 \cdot L_0) \circ \cdots \circ (c_{h+1} \cdot L_{h+1})) &\circ ((c'_0 \cdot L_0) \circ \cdots \circ (c'_{h+1} \cdot L_{h+1}))^{-1} = \\
&= \left(\sum_{i=0}^{h+1}(c_i - c'_i)\right) \cdot id
\end{aligned}
\qquad (18)
$$

In particular, we have

$$
\begin{aligned}
(c_0 \cdot L_0) &\circ \cdots \circ (c_{h+1} \cdot L_{h+1}) \circ ((c_0 \cdot L'_0) \circ \cdots \circ (c_{h+1} \cdot L'_{h+1}))^{-1} = \\
&= (c_0 \cdot L_0) \circ \cdots \circ (c_{h+1} \cdot L_{h+1}) \circ (((c_0 + d_0) \cdot L_0) \circ \cdots \circ ((c_{h+1} + d_{h+1}) \cdot L_{h+1}))^{-1} \\
&= \left(\sum_{i=0}^{h+1} c_i - \sum_{i=0}^{h+1}(d_i + c_i)\right) \cdot id \\
&= id
\end{aligned}
$$

and this concludes the proof of the independence of $(c_0 \cdot L_0) \circ \cdots \circ (c_{h+1} \cdot L_{h+1})$ from the $(h+2)$-tuple $(L_i)_i$.

Fixed $s, s'$ two sections of $\pi$ and $c$ an element of $H_1(S; \mathbb{Z}_N)$, we have

$$\sum_{i=0}^{h+1} c_i = \sum_{i=0}^{h+1} c'_i$$

where $s(c) = (c_1, \ldots, c_{h+1})$ and $s'(c) = (c'_1, \ldots, c'_{h+1})$. Hence relation 18 proves the independence of $c \cdot L$ from the fixed section too. In addition, by selecting an homomorphism as section, which clearly exists, we conclude that $(c, L) \mapsto c \cdot L$ is indeed an action. The last step in order to obtain the good definition of $(\mathscr{L}^{\rho\rho'}_{\lambda\lambda'}, \psi^{\rho\rho'}_{\lambda\lambda'})$ is the independence from the choice of the path of ideal triangulations between $\lambda$ and $\lambda'$.

Consider two sequences of ideal triangulations from $\lambda$ to $\lambda'$. By Theorem 1.21 we know that these sequences are connected by a path of certain moves. It is sufficient to prove that, starting from a sequence

$$\lambda = \lambda^{(0)}, \lambda^{(1)}, \ldots, \lambda^{(k)}, \ldots, \lambda^{(h)}, \lambda^{(h+1)} = \lambda'$$

and by applying any move described in Theorem 1.21, the sets $(\mathscr{L}^{\rho\rho'}_{\lambda\lambda'}, \psi^{\rho\rho'}_{\lambda\lambda'})$ defined through the different sequences are the same. In what follows, we will show that this fact is a simple consequence of Lemma 3.5. Take a sequence obtained by modifying the original one, which will look like

$$\lambda = \lambda^{(0)}, \lambda^{(1)}, \ldots, \lambda^{(k)} = \bar{\lambda}^{(0)}, \bar{\lambda}^{(1)}, \ldots, \bar{\lambda}^{(n+1)} = \lambda^{(k+1)}, \ldots, \lambda^{(h)}, \lambda^{(h+1)} = \lambda'$$

where $\lambda^{(k)} = \bar{\lambda}^{(0)}, \bar{\lambda}^{(1)}, \ldots, \bar{\lambda}^{(n)}$ is one of the sequences appearing in the assertion of Lemma 3.5. Then we need to compare the images of the following composition maps, labelled as $\varphi_1$ and $\varphi_2$ respectively:

$$\prod_{i=0}^{k-1} \mathscr{L}^{\rho\rho}_{\lambda^{(i)}\lambda^{(i+1)}} \times \prod_{j=0}^{n} \mathscr{L}^{\rho\rho}_{\bar{\lambda}^{(j)}\bar{\lambda}^{(j+1)}} \times \prod_{l=k+1}^{h} \mathscr{L}^{\rho\rho}_{\lambda^{(l)}\lambda^{(l+1)}} \times \mathscr{L}^{\rho\rho'}_{\lambda'\lambda'} \xrightarrow{\varphi_1} \operatorname{Hom}(V'_{\lambda'}, V_\lambda)/\doteq$$

$$\prod_{i=0}^{h} \mathscr{L}^{\rho\rho}_{\lambda^{(i)}\lambda^{(i+1)}} \times \mathscr{L}^{\rho\rho'}_{\lambda'\lambda'} \xrightarrow{\varphi_2} \operatorname{Hom}(V'_{\lambda'}, V_\lambda)/\doteq$$

Because of Lemma 3.5, the composition of an $(n+1)$-tuple in $\prod_{j=0}^{n} \mathscr{L}^{\rho\rho}_{\bar{\lambda}^{(j)}\bar{\lambda}^{(j+1)}}$ provides an element in $\mathscr{L}^{\rho\rho}_{\lambda^{(k)}\lambda^{(k+1)}}$, so the inclusion $\operatorname{Im} \varphi_1 \subseteq \operatorname{Im} \varphi_2$ is obvious. Moreover, the composition is surjective, so also the inverse inclusion holds and therefore $\operatorname{Im} \varphi_1 = \operatorname{Im} \varphi_2$. Hence the sets $\mathscr{L}^{\rho\rho'}_{\lambda\lambda'}$ are well defined, it remains to prove the good definition of the action $\psi^{\rho\rho'}_{\lambda\lambda'}$.

Fixed $L \in \mathscr{L}^{\rho\rho'}_{\lambda\lambda'}$, we can write it as an element in the image of $\varphi_1$ and $\varphi_2$, respectively, as follows

$$\begin{aligned} L &= L_0 \circ \cdots \circ L_{k-1} \circ \bar{L}_0 \circ \cdots \circ \bar{L}_n \circ L_{k+1} \circ \cdots \circ L_h \circ L_{h+1} \\ &= L'_0 \circ \cdots \circ L'_{k-1} \circ L'_k \circ L'_{k+1} \circ \cdots \circ L'_h \circ L'_{h+1} \end{aligned} \quad (19)$$

By virtue of the transitivity of the action on the terms $\mathscr{L}^{\rho\rho}_{\lambda^{(i)}\lambda^{(i+1)}}$ for every $i \neq k$, there exist $c_i \in H_1(S; \mathbb{Z}_N)$ such that $c_i \cdot L'_i = L_i$. Moreover, we can find an element $c_k \in H_1(S; \mathbb{Z}_N)$ such that $c_k \cdot L'_k = \bar{L}_0 \circ \cdots \circ \bar{L}_n$, thanks to Lemma 3.5. Denoting by $\prod_{i=0}^{n} M_i$ the composition $M_0 \circ \cdots \circ M_n$, from relation 19 we can deduce

$$\begin{aligned} id &= \prod_{i=0}^{k-1} L_i \circ \prod_{j=0}^{n} \bar{L}_j \circ \prod_{l=k+1}^{h+1} L_l \circ \prod_{l=0}^{h-k} (L'_{h+1-l})^{-1} \circ (L'_k)^{-1} \circ \prod_{i=0}^{k-1} (L'_{k-1-i})^{-1} \\ &= \prod_{i=0}^{k-1} (c_i \cdot L'_i) \circ (c_k \cdot L'_k) \circ \prod_{l=k+1}^{h+1} (c_l \cdot L'_l) \circ \prod_{l=0}^{h-k} (L'_{h+1-l})^{-1} \circ (L'_k)^{-1} \circ \prod_{i=0}^{k-1} (L'_{k-1-i})^{-1} \\ &= \left( \sum_{i=0}^{h+1} c_i \right) \cdot id \hfill \text{Relation 18} \end{aligned}$$

Because the action of $H_1(S; \mathbb{Z}_N)$ on $\mathscr{L}^{\rho\rho}_{\lambda\lambda}$ is free, we must have $\sum_{i=0}^{h+1} c_h = 0$. We have shown that the actions are independent from the choices of the sections, so given

$$\begin{aligned} \pi_1 : \quad H_1(S; \mathbb{Z}_N)^{h+n+2} &\longrightarrow H_1(S; \mathbb{Z}_N) \\ (d_i)_i &\longmapsto \sum_{i=0}^{h+n+1} d_i \end{aligned}$$

$$\begin{aligned} \pi_2 : \quad H_1(S; \mathbb{Z}_N)^{h+2} &\longrightarrow H_1(S; \mathbb{Z}_N) \\ (d_i)_i &\longmapsto \sum_{i=0}^{h+1} d_i \end{aligned}$$

the maps through which the actions are defined, we can choose

$$s_1(d) = (d, 0, \ldots, 0) \in H_1(S; \mathbb{Z}_N)^{h+n+2}$$
$$s_2(d) = (d, 0, \ldots, 0) \in H_1(S; \mathbb{Z}_N)^{h+2}$$

as sections. Then the actions of $d$, defined through these two sections, give respectively

$$(d \cdot L)_1 = (d \cdot L_0) \circ \cdots \circ L_{k-1} \circ \bar{L}_0 \circ \cdots \circ \bar{L}_n \circ L_{k+1} \circ \cdots \circ L_h \circ L_{h+1}$$
$$(d \cdot L)_2 = (d \cdot L'_0) \circ \cdots \circ L'_{k-1} \circ L'_k \circ L'_{k+1} \circ \cdots \circ L'_h \circ L'_{h+1}$$

where $(c \cdot L)_1$ denotes the action of $d$ on $L$ defined with the first sequence and the section $s_1$, and $(c \cdot L)_1$ denotes the action of $d$ on $L$ defined with the second sequence and the section $s_2$.

Now, by virtue of the presentations of $(d \cdot L)_1$ and $(d \cdot L)_2$ just given and of the relation 18, we can rewrite the isomorphism $(d \cdot L)_1 \circ (d \cdot L)_2^{-1}$ as the following composition:

$$= (d \cdot L_0) \circ \prod_{i=1}^{k-1} L_i \circ \prod_{j=0}^{n} \bar{L}_j \circ \prod_{l=k+1}^{h+1} L_l \circ \prod_{l=0}^{h-k} (L'_{h+1-l})^{-1} \circ (L'_k)^{-1} \circ \prod_{i=0}^{k-2} (L'_{k-1-i})^{-1} \circ$$
$$\circ (d \cdot L'_0)^{-1}$$

$$= ((d + c_0) \cdot L'_0) \circ \prod_{i=1}^{k-1} (c_i \cdot L'_i) \circ (c_k \cdot L'_k) \circ \prod_{l=k+1}^{h+1} (c_l \cdot L'_l) \circ \prod_{l=0}^{h-k} (L'_{h+1-l})^{-1} \circ (L'_k)^{-1} \circ$$
$$\circ \prod_{i=0}^{k-2} (L'_{k-1-i})^{-1} \circ (d \cdot L'_0)^{-1}$$

$$= \left( d + \sum_{i=0}^{h+1} c_i - d \right) \cdot id$$
$$= id$$

and this finally proves the independence of the action from the chosen path of ideal triangulations. Hence the objects $(\mathcal{L}^{\rho\rho'}_{\lambda\lambda'}, \psi^{\rho\rho'}_{\lambda\lambda'})$ are well defined.

In order to prove the transitivity of the action, fix a certain path of ideal triangulations and two elements $L = L_0 \circ \cdots \circ L_{h+1}$ and $L' = L'_0 \circ \cdots \circ L'_{h+1}$ in $\mathcal{L}^{\rho\rho'}_{\lambda\lambda'}$. Because the actions are transitive in the elementary cases, for every $i = 0, \ldots, h+1$ there exists a $c_i \in H_1(S; \mathbb{Z}_N)$ such that $c_i \cdot L_i = L'_i$. So, choosing a section $s$ such that $s(\sum_i c_i) = (c_0, \ldots, c_n)$ (because $s$ is not required to be a homomorphism, we can always do so), we observe that

$$c \cdot (L_0 \circ \cdots \circ L_{h+1}) = L'_0 \circ \cdots \circ L'_{h+1}$$

hence the action is transitive. Now suppose that $c \cdot L = L$, then the following relation must hold

$$((c_0 \cdot L_0) \circ \cdots \circ (c_{h+1} \cdot L_{h+1})) \circ (L_0 \circ \cdots \circ L_{h+1})^{-1} = id$$

Now, applying the relation 18 we obtain

$$\left( \sum_{i=0}^{h+1} c_i \right) \cdot id = id$$

By freeness of the action of $H_1(S; \mathbb{Z}_N)$ on $\mathscr{L}_{\lambda\lambda}^{\rho\rho}$, we deduce $\sum_i c_i = 0$ and then $c = 0$, which proves the freeness of the action in the generic case.

Now we have defined all the elements needed to deal with the proof of the Existence Theorem:

**Theorem 4.1** (Existence Theorem). *For every surface $S$ there exists a collection $\{(\mathscr{L}_{\lambda\lambda'}^{\rho\rho'}, \psi_{\lambda\lambda'}^{\rho\rho'})\}$, indexed by couples of isomorphic local representations $\rho, \rho'$ of the quantum Teichmüller space $\mathcal{T}_S^q$ and by couples of ideal triangulations $\lambda, \lambda' \in \Lambda(S)$ such that*

INTERTWINING: for every couple of isomorphic local representations

$$\rho = \{\rho_\lambda \colon \mathcal{T}_\lambda^q \to \mathrm{End}(V_\lambda)\}_{\lambda \in \Lambda(S)} \quad \rho' = \{\rho'_\lambda \colon \mathcal{T}_\lambda^q \to \mathrm{End}(V'_\lambda)\}_{\lambda \in \Lambda(S)}$$

and for every $\lambda, \lambda' \in \Lambda(S)$, $\mathscr{L}_{\lambda\lambda'}^{\rho\rho'}$ is a set of linear isomorphisms $L_{\lambda\lambda'}^{\rho\rho'}$ from $V'_{\lambda'}$ to $V_\lambda$, considered up to scalar multiplication, verifying

$$(\rho_\lambda \circ \Phi_{\lambda\lambda'}^q)(X') = L_{\lambda\lambda'}^{\rho\rho'} \circ \rho'_{\lambda'}(X') \circ (L_{\lambda\lambda'}^{\rho\rho'})^{-1} \qquad \forall X' \in \mathcal{T}_{\lambda'}^q$$

ACTION: every set $\mathscr{L}_{\lambda\lambda'}^{\rho\rho'}$ is endowed with a transitive and free action $\psi_{\lambda\lambda'}^{\rho\rho'}$ of $H_1(S; \mathbb{Z}_N)$;

FUSION PROPERTY: let $R$ be a surface and $S$ obtained by fusion from $R$. Fix

$$\eta = \{\eta_\mu \colon \mathcal{T}_\mu^q \to \mathrm{End}(W_\mu)\}_{\mu \in \Lambda(R)} \quad \eta' = \{\eta'_\mu \colon \mathcal{T}_\mu^q \to \mathrm{End}(W'_\mu)\}_{\mu \in \Lambda(R)}$$

two isomorphic local representations of $\mathcal{T}_R^q$ and

$$\rho = \{\rho_\lambda \colon \mathcal{T}_\lambda^q \to \mathrm{End}(V_\lambda)\}_{\lambda \in \Lambda(S)} \quad \rho' = \{\rho'_\lambda \colon \mathcal{T}_\lambda^q \to \mathrm{End}(V'_\lambda)\}_{\lambda \in \Lambda(S)}$$

two isomorphic local representations of $\mathcal{T}_S^q$, with $\rho$ and $\rho'$ obtained by fusion from $\eta$ and $\eta'$, respectively. Then for every $\mu, \mu' \in \Lambda(R)$, if $\lambda, \lambda' \in \Lambda(S)$ are the corresponding ideal triangulations on $S$, there exists a natural inclusion $j \colon \mathscr{L}_{\mu\mu'}^{\eta\eta'} \to \mathscr{L}_{\lambda\lambda'}^{\rho\rho'}$ (in fact the inclusion as subsets of the quotient $\mathrm{Hom}(V'_{\lambda'}, V_\lambda)/\doteq$) such that, for every $L \in \mathscr{L}_{\mu\mu'}^{\eta\eta'}$, the following holds

$$(j \circ \psi_{\mu\mu'}^{\eta\eta'})(c, L) = \psi_{\lambda\lambda'}^{\rho\rho'}(\pi_*(c), j(L))$$

for every $c \in H_1(R; \mathbb{Z}_N)$, where $\pi \colon R \to S$ is the projection map;

COMPOSITION PROPERTY: for every $\rho, \rho', \rho''$ isomorphic local representations of $\mathcal{T}_S^q$ and for every $\lambda, \lambda', \lambda'' \in \Lambda(S)$ the composition map

$$\begin{array}{ccc} \mathscr{L}_{\lambda\lambda'}^{\rho\rho'} \times \mathscr{L}_{\lambda'\lambda''}^{\rho'\rho''} & \longrightarrow & \mathscr{L}_{\lambda\lambda''}^{\rho\rho''} \\ (L, M) & \longmapsto & L \circ M \end{array}$$

is well defined and it verifies

$$(c \cdot L) \circ (d \cdot M) = (c + d) \cdot (L \circ M)$$

*Proof.* We need to verify that Fusion and Composition properties hold.

FUSION PROPERTY: Let $R$ and $S$ be surfaces like in the assertion, each one endowed with a couple of ideal triangulations, $\mu, \mu' \in \Lambda(R)$ and $\lambda, \lambda' \in \Lambda(S)$ respectively, where the fusions of the firsts are the seconds, and with a couple of isomorphic local representations, $\eta, \eta'$ of $\mathcal{T}_R^q$ and $\rho, \rho'$ of $\mathcal{T}_S^q$, where the fusions of the firsts are the seconds. Firstly observe that
$W_\mu = V_\lambda$ and $W'_{\mu'} = V'_{\lambda'}$ because $\rho$ is fusion of $\eta$ and $\rho'$ is fusion of $\eta'$, so the sets $\mathscr{L}_{\lambda\lambda'}^{\rho\rho'}$ and $\mathscr{L}_{\mu\mu'}^{\eta\eta'}$ are both subsets of $\mathrm{Hom}(V'_{\lambda'}, V_\lambda)/\doteq$.

$$\begin{array}{ccc} \nu\,\nu' & R & \mu\,\mu' \\ & \downarrow & \\ \lambda\,\lambda' & S & \rho\,\rho' \end{array}$$

Now we choose a sequence of ideal triangulations in $R$ from $\mu$ to $\mu'$
$$\mu = \mu^{(0)}, \mu^{(1)}, \ldots, \mu^{(k)}, \mu^{(k+1)} = \mu'$$
which induces a corresponding sequence in $S$ from $\lambda$ to $\lambda'$ by fusion
$$\lambda = \lambda^{(0)}, \lambda^{(1)}, \ldots, \lambda^{(k)}, \lambda^{(k+1)} = \lambda'$$
The set $\mathscr{L}_{\mu\mu'}^{\eta\eta'}$ can be realized as the image under the composition map of the set
$$\mathscr{L}_{\mu\mu^{(1)}}^{\eta\eta} \times \cdots \times \mathscr{L}_{\mu^{(k)}\mu'}^{\eta\eta} \times \mathscr{L}_{\mu'\mu'}^{\eta\eta'}$$
Lemma 3.1 tells us that this product is contained, through a natural map that we still denote by $j$, in
$$\mathscr{L}_{\lambda\lambda^{(1)}}^{\rho\rho} \times \cdots \times \mathscr{L}_{\lambda^{(k)}\lambda'}^{\rho\rho} \times \mathscr{L}_{\lambda'\lambda'}^{\rho\rho'}$$
The image of this last set under the composition map is equal to $\mathscr{L}_{\lambda\lambda'}^{\rho\rho'}$. This clearly shows $\mathscr{L}_{\mu\mu'}^{\eta\eta'} \subseteq \mathscr{L}_{\lambda\lambda'}^{\rho\rho'}$ and the map is just the inclusion as subsets of $\mathrm{Hom}(V'_{\lambda'}, V_\lambda)/\doteq$, so it does not depend on the chosen sequence of ideal triangulations. Now we have to prove the relation between the actions: we fix $c \in H_1(R; \mathbb{Z}_N)$ and an element $L \in \mathscr{L}_{\mu\mu'}^{\eta\eta'}$, we take a presentation $L = L_0 \circ \cdots \circ L_{k+1}$ as compositions of a $(k+2)$-tuple in
$$\mathscr{L}_{\mu\mu^{(1)}}^{\eta\eta} \times \cdots \times \mathscr{L}_{\mu^{(k)}\mu'}^{\eta\eta} \times \mathscr{L}_{\mu'\mu'}^{\eta\eta'}$$
Then, recalling Lemma 3.1 and what just seen, we have
$$\begin{aligned} j(c \cdot (L_0 \circ \cdots \circ L_{k+1})) &= j((c_0 \cdot L_0) \circ \cdots \circ (c_{k+1} \cdot L_{k+1})) \\ &= j(c_0 \cdot L_0) \circ \cdots \circ j(c_{k+1} \cdot L_{k+1}) \\ &= (\pi_*(c_0) \cdot j(L_0)) \circ \cdots \circ (\pi_*(c_{k+1}) \cdot j(L_{k+1})) \\ &= \pi_*(c) \cdot (j(L_0) \circ \cdots \circ j(L_{k+1})) \\ &= \pi_*(c) \cdot j(L_0 \circ \cdots \circ L_{k+1}) \end{aligned}$$
where $\pi \colon R \to S$ is the quotient map. So the Fusion property holds.

COMPOSITION PROPERTY: Given $\lambda, \lambda', \lambda'' \in \Lambda(S)$, we choose paths of elementary moves on ideal triangulations
$$\lambda = \lambda^{(0)}, \ldots, \lambda^{(k)}, \lambda^{(k+1)} = \lambda', \lambda^{(k+2)}, \ldots, \lambda^{(h)}, \lambda^{(h+1)} = \lambda''$$
Then the sets $\mathscr{L}_{\lambda\lambda'}^{\rho\rho'}$ and $\mathscr{L}_{\lambda'\lambda''}^{\rho'\rho''}$ can be constructed as images of the composition maps as follows
$$\mathscr{L}_{\lambda\lambda'}^{\rho\rho'} = \mathrm{Im}(\mathscr{L}_{\lambda\lambda^{(1)}}^{\rho\rho} \times \cdots \times \mathscr{L}_{\lambda^{(k)}\lambda'}^{\rho\rho} \times \mathscr{L}_{\lambda'\lambda'}^{\rho\rho'} \longrightarrow \mathrm{Hom}(V'_{\lambda'}, V_\lambda)/\doteq)$$
$$\mathscr{L}_{\lambda'\lambda''}^{\rho'\rho''} = \mathrm{Im}(\mathscr{L}_{\lambda'\lambda^{(k+2)}}^{\rho'\rho'} \times \cdots \times \mathscr{L}_{\lambda^{(h)}\lambda''}^{\rho'\rho'} \times \mathscr{L}_{\lambda''\lambda''}^{\rho'\rho''} \longrightarrow \mathrm{Hom}(V''_{\lambda''}, V'_{\lambda'})/\doteq)$$

On the other hand, the set $\mathscr{L}_{\lambda\lambda''}^{\rho\rho''}$ can be presented as

$$\mathscr{L}_{\lambda\lambda''}^{\rho\rho''} = \mathrm{Im}(\mathscr{L}_{\lambda\lambda^{(1)}}^{\rho\rho} \times \cdots \times \mathscr{L}_{\lambda^{(h)}\lambda''}^{\rho\rho} \times \mathscr{L}_{\lambda''\lambda''}^{\rho\rho''} \longrightarrow \mathrm{Hom}(V_{\lambda''}'', V_\lambda)/\doteq)$$

Firstly we need to show that this set is equal to

$$\mathrm{Im}\left(\prod_{i=0}^{k} \mathscr{L}_{\lambda^{(i)}\lambda^{(i+1)}}^{\rho\rho} \times \mathscr{L}_{\lambda'\lambda'}^{\rho\rho'} \times \prod_{i=k+1}^{h} \mathscr{L}_{\lambda^{(i)}\lambda^{(i+1)}}^{\rho'\rho'} \times \mathscr{L}_{\lambda''\lambda''}^{\rho'\rho''} \to \mathrm{Hom}(V_{\lambda''}'', V_\lambda)/\doteq\right) \tag{20}$$

and then to explain how the actions are related. The first fact is an easy consequence of Lemma 3.6. Indeed, by virtue of that result, the following holds

$$\mathrm{Im}(\mathscr{L}_{\lambda'\lambda'}^{\rho\rho'} \times \mathscr{L}_{\lambda'\lambda^{(k+2)}}^{\rho'\rho'} \longrightarrow \mathrm{Hom}(V_{\lambda^{(k+2)}}', V_{\lambda'})/\doteq) =$$
$$= \mathrm{Im}(\mathscr{L}_{\lambda'\lambda^{(k+2)}}^{\rho\rho} \times \mathscr{L}_{\lambda^{(k+2)}\lambda^{(k+2)}}^{\rho\rho'} \longrightarrow \mathrm{Hom}(V_{\lambda^{(k+2)}}', V_{\lambda'})/\doteq)$$

where the maps are the obvious compositions. This implies that the set in relation 20 is equal to the image of the composition on the set

$$\prod_{i=0}^{k+1} \mathscr{L}_{\lambda^{(i)}\lambda^{(i+1)}}^{\rho\rho} \times \mathscr{L}_{\lambda^{(k+2)}\lambda^{(k+2)}}^{\rho\rho'} \times \prod_{i=k+2}^{h} \mathscr{L}_{\lambda^{(i)}\lambda^{(i+1)}}^{\rho'\rho'} \times \mathscr{L}_{\lambda''\lambda''}^{\rho'\rho''}$$

Now, iterating this process we conclude that the set in 20 is equal to the image of the composition on

$$\prod_{i=0}^{h} \mathscr{L}_{\lambda^{(i)}\lambda^{(i+1)}}^{\rho\rho} \times \mathscr{L}_{\lambda''\lambda''}^{\rho\rho'} \times \mathscr{L}_{\lambda''\lambda''}^{\rho'\rho''}$$

and by applying one last time Lemma 3.6 on $\mathscr{L}_{\lambda''\lambda''}^{\rho\rho'} \times \mathscr{L}_{\lambda''\lambda''}^{\rho'\rho''}$ we obtain the equality we are looking for. Now it remains to prove the relation between the actions. Given $L \in \mathscr{L}_{\lambda\lambda'}^{\rho\rho'}$ and $M \in \mathscr{L}_{\lambda'\lambda''}^{\rho'\rho''}$, we write $L$ as composition of a certain

$$(L_0, \ldots, L_{k+1}) \in \mathscr{L}_{\lambda\lambda^{(1)}}^{\rho\rho} \times \cdots \times \mathscr{L}_{\lambda^{(k)}\lambda'}^{\rho\rho} \times \mathscr{L}_{\lambda'\lambda'}^{\rho\rho'}$$

and analogously $M$ as composition of

$$(M_{k+2}, \ldots, M_{h+1}) \in \mathscr{L}_{\lambda'\lambda^{(k+2)}}^{\rho'\rho'} \times \cdots \times \mathscr{L}_{\lambda^{(h)}\lambda''}^{\rho'\rho'} \times \mathscr{L}_{\lambda''\lambda''}^{\rho'\rho''}$$

The element $L_{k+1} \circ M_{k+2}$ belongs to $\mathscr{L}_{\lambda'\lambda^{(k+2)}}^{\rho\rho'}$, which is equal to the image of compositions of isomorphisms in

$$\mathscr{L}_{\lambda'\lambda^{(k+2)}}^{\rho\rho} \times \mathscr{L}_{\lambda^{(k+2)}\lambda^{(k+2)}}^{\rho\rho'}$$

by virtue of Lemma 3.6. Then there exists a couple $(L'_{k+1}, \overline{M}_{k+2})$ in the set $\mathscr{L}_{\lambda'\lambda^{(k+2)}}^{\rho\rho} \times \mathscr{L}_{\lambda^{(k+2)}\lambda^{(k+2)}}^{\rho\rho'}$ such that $L_{k+1} \circ M_{k+2} = L'_{k+1} \circ \overline{M}_{k+2}$. In this way we have written $L \circ M$ as composition of an element in

$$\prod_{i=0}^{k+1} \mathscr{L}_{\lambda^{(i)}\lambda^{(i+1)}}^{\rho\rho} \times \mathscr{L}_{\lambda^{(k+2)}\lambda^{(k+2)}}^{\rho\rho'} \times \prod_{i=k+2}^{h} \mathscr{L}_{\lambda^{(i)}\lambda^{(i+1)}}^{\rho'\rho'} \times \mathscr{L}_{\lambda''\lambda''}^{\rho'\rho''}$$

Iterating this process, we rewrite $L \circ M$ as follows

$$\begin{aligned} L \circ M &= L_0 \circ \cdots \circ L_k \circ L_{k+1} \circ M_{k+1} \circ M_{k+2} \circ \cdots \circ M_{h+1} \\ &= L_0 \circ \cdots \circ L_k \circ L'_{k+1} \circ \overline{M}_{k+2} \circ M_{k+3} \circ \cdots \circ M_{h+1} \\ &= L_0 \circ \cdots \circ L_k \circ L'_{k+1} \circ M'_{k+2} \circ \overline{M}_{k+3} \circ \cdots \circ M_{h+1} \\ &\vdots \\ &= L_0 \circ \cdots \circ L_k \circ L'_{k+1} \circ M'_{k+2} \circ M'_{k+3} \circ \cdots \circ M'_{h+1} \end{aligned}$$

where

$$(L_0, \ldots, L_k, L'_{k+1}, M'_{k+2}, \ldots, M'_{h+1}) \in \mathscr{L}^{\rho\rho}_{\lambda\lambda^{(1)}} \times \cdots \times \mathscr{L}^{\rho\rho}_{\lambda^{(h)}\lambda''} \times \mathscr{L}^{\rho\rho''}_{\lambda''\lambda''}$$

so we have found a decomposition of $L \circ M$ as element of $\mathscr{L}^{\rho\rho''}_{\lambda\lambda''}$, described by the path $\lambda = \lambda^{(0)}, \ldots, \lambda^{(h+1)} = \lambda''$. The image of $L$ under the action of an element $c \in H_1(S; \mathbb{Z}_N)$ will have the form

$$c \cdot L = (c_0 \cdot L_0) \circ \cdots \circ (c_{k+1} \cdot L_{k+1})$$

with $\sum_i c_i = c$, and analogously the image of $M$ under the action of $d \in H_1(S; \mathbb{Z}_N)$ will appear like

$$d \cdot M = (d_{k+2} \cdot M_{k+2}) \circ \cdots \circ (d_{h+1} \cdot M_{h+1})$$

Recalling the relation between the actions in Lemma 3.6 we see that

$$\begin{aligned} (c_{k+1} \cdot L_{k+1}) \circ (d_{k+2} \circ M_{k+2}) &= (c_{k+1} + d_{k+2}) \cdot (L_{k+1} \circ M_{k+2}) \\ &= (c_{k+1} + d_{k+2}) \cdot (L'_{k+1} \circ \overline{M}_{k+2}) \\ &= (c_{k+1} \cdot L'_{k+1}) \circ (d_{k+2} \circ \overline{M}_{k+2}) \end{aligned}$$

This implies the following relation

$$\begin{aligned} (c \cdot L) \circ (d \cdot M) &= (c_0 \cdot L_0) \circ \cdots \circ (c_{k+1} \cdot L_{k+1}) \circ (d_{k+2} \cdot M_{k+2}) \circ \cdots \circ M_{h+1} \\ &= (c_0 \cdot L_0) \circ \cdots \circ (c_{k+1} \cdot L'_{k+1}) \circ (d_{k+2} \cdot \overline{M}_{k+2}) \circ \cdots \circ (d_{h+1} \circ M_{h+1}) \end{aligned}$$

By iterating this process as before we obtain that

$$(c \cdot L) \circ (d \cdot M) = (c_0 \cdot L_0) \circ \cdots \circ (c_{k+1} \cdot L'_{k+1}) \circ (d_{k+2} \cdot M'_{k+2}) \circ \cdots \circ (d_{h+1} \cdot M'_{h+1})$$

Now observe that the second member is equal to $(c + d) \cdot (L \circ M) \in \mathscr{L}^{\rho\rho''}_{\lambda\lambda''}$ by definition of $(\mathscr{L}^{\rho\rho''}_{\lambda\lambda''}, \psi^{\rho\rho''}_{\lambda\lambda''})$ and because

$$\sum_{i=0}^{k+1} c_i + \sum_{j=k+2}^{h+1} d_j = c + d \in H_1(S; \mathbb{Z}_N)$$

This concludes the proof. $\square$

**Theorem 4.2** (Uniqueness Theorem). *Suppose that $\{\mathscr{M}^{\rho\rho'}_{\lambda\lambda'}\}$ is a collection indexed by couples of isomorphic local representations $\rho, \rho'$ of the quantum Teichmüller space $\mathcal{T}^q_S$ and by couples of ideal triangulations $\lambda, \lambda' \in \Lambda(S)$ such that*

INTERTWINING: for every couple of isomorphic local representations

$$\rho = \{\rho_\lambda \colon \mathcal{T}_\lambda^q \to \mathrm{End}(V_\lambda)\}_{\lambda \in \Lambda(S)} \quad \rho' = \{\rho'_\lambda \colon \mathcal{T}_\lambda^q \to \mathrm{End}(V'_\lambda)\}_{\lambda \in \Lambda(S)}$$

and for every $\lambda, \lambda' \in \Lambda(S)$, $\mathscr{M}_{\lambda\lambda'}^{\rho\rho'}$ is a non-empty set of linear isomorphisms $M_{\lambda\lambda'}^{\rho\rho'}$ from $V'_{\lambda'}$ to $V_\lambda$, considered up to scalar multiplication, such that

$$(\rho_\lambda \circ \Phi_{\lambda\lambda'}^q)(X') = M_{\lambda\lambda'}^{\rho\rho'} \circ \rho'_{\lambda'}(X') \circ (M_{\lambda\lambda'}^{\rho\rho'})^{-1} \qquad \forall X' \in \mathcal{T}_{\lambda'}^q$$

WEAK FUSION PROPERTY: let $R$ be a surface and $S$ obtained by fusion from $R$. Fix

$$\eta = \{\eta_\mu \colon \mathcal{T}_\mu^q \to \mathrm{End}(W_\mu)\}_{\mu \in \Lambda(R)} \quad \eta' = \{\eta'_\mu \colon \mathcal{T}_\mu^q \to \mathrm{End}(W'_\mu)\}_{\mu \in \Lambda(R)}$$

two isomorphic local representations of $\mathcal{T}_R^q$ and

$$\rho = \{\rho_\lambda \colon \mathcal{T}_\lambda^q \to \mathrm{End}(V_\lambda)\}_{\lambda \in \Lambda(S)} \quad \rho' = \{\rho'_\lambda \colon \mathcal{T}_\lambda^q \to \mathrm{End}(V'_\lambda)\}_{\lambda \in \Lambda(S)}$$

two isomorphic local representations of $\mathcal{T}_S^q$, with $\rho$ and $\rho'$ obtained by fusion from $\eta$ and $\eta'$, respectively. Then for every $\mu, \mu' \in \Lambda(R)$, if $\lambda, \lambda' \in \Lambda(S)$ are the corresponding ideal triangulations on $S$, the inclusion map $j \colon \mathscr{M}_{\mu\mu'}^{\eta\eta'} \to \mathscr{M}_{\lambda\lambda'}^{\rho\rho'}$ as subset of $\mathrm{Hom}(V'_{\lambda'}, V_\lambda)/\doteq$ is well defined;

WEAK COMPOSITION PROPERTY: for every $\rho, \rho', \rho''$ isomorphic local representations of $\mathcal{T}_S^q$ and for every $\lambda, \lambda', \lambda'' \in \Lambda(S)$ the composition map

$$\begin{array}{ccc} \mathscr{M}_{\lambda\lambda'}^{\rho\rho'} \times \mathscr{M}_{\lambda'\lambda''}^{\rho'\rho''} & \longrightarrow & \mathscr{M}_{\lambda\lambda''}^{\rho\rho''} \\ (M, N) & \longmapsto & M \circ N \end{array}$$

is well defined.

Then, for every $\rho$ and $\rho'$ isomorphic local representations and for every $\lambda, \lambda' \in \Lambda(S)$ we have

$$\mathscr{L}_{\lambda\lambda'}^{\rho\rho'} \subseteq \mathscr{M}_{\lambda\lambda'}^{\rho\rho'}$$

where $\{\mathscr{L}_{\lambda\lambda'}^{\rho\rho'}\}$ is the family previously constructed.

*Proof.* Thanks to the Weak Composition property and to the surjectivity of the composition maps for the $(\mathscr{L}_{\lambda\lambda'}^{\rho\rho'}, \psi_{\lambda\lambda'}^{\rho\rho'})$, it is sufficient to show the inclusion $\mathscr{L}_{\lambda\lambda'}^{\rho\rho'} \subseteq \mathscr{M}_{\lambda\lambda'}^{\rho\rho'}$ in the elementary cases, in which the triangulations differ by a diagonal exchange or a re-indexing. Let $S$ be a surface and take $\lambda = \lambda' \in \Lambda(S)$, the other situation is analogous. Denote by $S_0$ the surface obtained by splitting $S$ along all the edges of $\lambda$ and by $\lambda_0$ the ideal triangulation induced on $S_0$. Moreover, we fix $\rho$ and $\rho'$ two isomorphic local representations of $\mathcal{T}_S^q$ and we choose two isomorphic representatives $\zeta_{\lambda_0}$ and $\zeta'_{\lambda_0}$ of $\rho_\lambda$ and $\rho'_\lambda$ respectively. The representations $\zeta_{\lambda_0}$ and $\zeta'_{\lambda_0}$ can be thought as local representations $\zeta$ and $\zeta'$ of the whole quantum Teichmüller space of $\mathcal{T}_{S_0}^q$, because $S_0$ admits the only triangulation $\lambda_0$, being a disjoint union of ideal triangles. The element $\zeta'_{\lambda_0}$ belongs by construction to the set $\mathscr{F}_{S_0}(\rho'_\lambda)$, so we can consider $c \cdot \zeta'_{\lambda_0}$ for every $c \in H_1(S; \mathbb{Z}_N)$. In this way, for every $c \in H_1(S; \mathbb{Z}_N)$ we obtain a local representation $c \cdot \zeta'$ of $\mathcal{T}_{S_0}^q$ isomorphic to $\zeta'$ and that still leads by fusion to $\rho'_\lambda$. Because $\zeta$ and $c \cdot \zeta'$ are isomorphic irreducible representations of $\mathcal{T}^q_{\lambda_0}$, there exists a linear isomorphism $L^{\zeta,c\cdot\zeta'}_{\lambda\lambda'} : V'_{\lambda'} \to V_\lambda$, unique up to multiplicative scalar, verifying

$$L^{\zeta,c\cdot\zeta'}_{\lambda\lambda} \circ \zeta'_{\lambda_0}(X) \circ (L^{\zeta,c\cdot\zeta'}_{\lambda\lambda})^{-1} = \zeta_{\lambda_0}(X) \qquad \forall X \in \mathcal{T}^q_{\lambda_0}$$

by virtue of Proposition 1.14. Then, because $\mathscr{M}^{\zeta\zeta'}_{\lambda_0\lambda_0}$ must be non-empty, the isomorphism $L^{\zeta,c\cdot\zeta'}_{\lambda_0\lambda_0}$ necessarily belongs to $\mathscr{M}^{\zeta,c\cdot\zeta'}_{\lambda_0\lambda_0}$. The weak Fusion property tells us that $\mathscr{M}^{\zeta,c\cdot\zeta'}_{\lambda_0\lambda_0}$ is contained in $\mathscr{M}^{\rho\rho'}_{\lambda\lambda}$, so the isomorphism $L^{\zeta,c\cdot\zeta'}_{\lambda_0\lambda_0}$ must belong to $\mathscr{M}^{\rho\rho'}_{\lambda\lambda}$ for every $c \in H_1(S;\mathbb{Z}_N)$. By definition $L^{\zeta\zeta'}_{\lambda_0\lambda_0}$ is an element of $\mathscr{L}^{\rho\rho'}_{\lambda\lambda}$, and $L^{\zeta,c\cdot\zeta'}_{\lambda_0\lambda_0}$ coincides with $c\cdot L^{\zeta\zeta'}_{\lambda_0\lambda_0}$, where $c \in H_1(S;\mathbb{Z}_N)$ is acting on $\mathscr{L}^{\rho\rho'}_{\lambda\lambda}$. This means that, by transitivity of the action $\psi^{\rho\rho'}_{\lambda\lambda}$, the whole set $\mathscr{L}^{\rho\rho'}_{\lambda\lambda}$ is contained in $\mathscr{M}^{\rho\rho'}_{\lambda\lambda}$, which is what we were looking for. $\square$

**Lemma 4.3.** Let $\lambda, \lambda' \in \Lambda(S)$ be two ideal triangulations and

$$\rho = \{\rho_\lambda \colon \mathcal{T}^q_\lambda \to \mathrm{End}(V_\lambda)\}_{\lambda \in \Lambda(S)} \qquad \rho' = \{\rho'_\lambda \colon \mathcal{T}^q_\lambda \to \mathrm{End}(V'_\lambda)\}_{\lambda \in \Lambda(S)}$$

two isomorphic local representations of the quantum Teichmüller space $\mathcal{T}^q_S$ of $S$. Then for every $c \in H_1(S;\mathbb{Z}_N)$ there exists an automorphism $B(c)$ of $V'_{\lambda'}$ with $\det B(c) = 1$, uniquely determined up to scalar multiplication by an $N$-th root of unity, such that

$$c \cdot L \doteq L \circ B(c)^{-1}$$

for every $L \in \mathscr{L}^{\rho\rho'}_{\lambda\lambda'}$.

*Proof.* Firstly we observe that we can assume $\lambda = \lambda'$. Indeed, every $L \in \mathscr{L}^{\rho\rho'}_{\lambda\lambda'}$ can be written as $L_0 \circ L_1$, with

$$(L_0, L_1) \in \mathscr{L}^{\rho\rho}_{\lambda\lambda'} \times \mathscr{L}^{\rho\rho'}_{\lambda'\lambda'}$$

and the element $c \cdot L$ is equal to $L_0 \circ (c \cdot L_1)$. So, by showing that the condition holds for $\mathscr{L}^{\rho\rho'}_{\lambda'\lambda'}$, we will conclude the general case.

An element $L \in \mathscr{L}^{\rho\rho'}_{\lambda'\lambda'}$ corresponds to a certain class $[\zeta, \zeta'] \in \mathscr{A}^{\rho\rho'}_{\lambda'\lambda'}$, where $\zeta$ and $\zeta'$ are representatives of $\rho_{\lambda'}$ and $\rho'_{\lambda'}$, respectively, and the following holds

$$L \circ \zeta'(X) \circ L^{-1} = \zeta(X) \qquad \forall X \in \mathcal{T}^q_{\lambda'_0}$$

where $\lambda'_0$ is the triangulation on the surface $S'_0$, obtained by splitting $S$ along the triangulation $\lambda'$. On the other hand, for every $c \in H_1(S;\mathbb{Z}_N)$ the element $c \cdot L$ corresponds to the class $[\zeta, c \cdot \zeta']$, where $c \cdot \zeta'$ is the action of $c$ on $\zeta' \in \mathscr{F}_{S'_0}(\rho'_{\lambda'})$. In Remark 1.19 we observed that $\zeta'$ and $c \cdot \zeta'$ are isomorphic and that there exists a linear isomorphism $D(c)$, described explicitly, such that

$$D(c) \circ \zeta'(X) \circ D(c)^{-1} = (c \cdot \zeta')(X) \qquad \forall X \in \mathcal{T}^q_{\lambda'_0}$$

These two relations imply immediately that

$$(L \circ D(c)^{-1}) \circ (c \cdot \zeta')(X) \circ (L \circ D(c)^{-1})^{-1} = \zeta(X) \qquad \forall X \in \mathcal{T}^q_{\lambda'_0}$$

and this proves that $c \cdot L \doteq L \circ D(c)^{-1}$. By asking that $\det D(c) = 1$, we obtain a linear isomorphism $B(c)$, uniquely determined up to scalar multiplication by an $N$-th root of unity, verifying the requests.

More precisely, in Remark 1.19 we have found $D(c)$ as conjugated to linear isomorphisms that are tensor products in $\mathrm{GL}(\bigotimes_i (\mathbb{C}^N)_i)$ of isomorphisms of $\mathbb{C}^N$ obtained as compositions of the applications $B_i$ for $i = 1, 2, 3$ and their inverses. It is immediate to see that every $B_i$ verifies $\det(B_i) = (-1)^{N+1}$, where we are using the fact that $q^N = (-1)^{N+1}$. Furthermore, the following relation holds

$$\det(L_1 \otimes \cdots \otimes L_m) = \det(L_1)^{N^{m-1}} \cdots \det(L_1)^{N^{m-1}} \tag{21}$$

for every $L_1, \ldots, L_k \in \mathrm{GL}_n(\mathbb{C})$, where $m$ is the number of triangles in an ideal triangulation of $S$. So, if $D(c)$ is equal to $A_1 \otimes \cdots \otimes A_m$, with $A_i$ composition of the $B_i$, then $\det(A_i)$ is a certain power of $(-1)^{N+1}$ and consequently $\det(A_i)^{N^{m-1}}$ is a power of $(-1)^{N^{m-1}(N+1)}$. If $m > 1$, then $N^{m-1}(N+1)$ is even, and consequently $\det(A_i)^{N^{m-1}} = 1$, which proves that $\det(D(c)) = 1$ thanks to the relation 21. So actually we do not need to rescale $D(c)$ in order to obtain the additional property $\det = 1$ if $m > 1$. □

# 5 Invariants of pseudo-Anosov diffeomorphisms

Now we adapt the construction of invariants of pseudo-Anosov surface diffeomorphisms described in [BBL07] by replacing [BBL07, Theorem 20] with our reformulation. In the first Subsection we will recall the known classification Theorems exposed in [BBL07] for local representations and in the second one we will give the details of the construction of such invariants.

## 5.1 Classification of local representations

Let $S$ be a surface and $\lambda \in \Lambda(S)$ be an ideal triangulation. Given $\rho_\lambda \colon \mathcal{T}_\lambda^q \to \mathrm{End}(V)$ a local representation, denote by $x_i \in \mathbb{C}^*$, for varying $i = 1, \ldots, n$, the numbers such that

$$\rho_\lambda(X_i^N) = x_i \, id_V$$

The algebra $\mathcal{T}_\lambda^1$ is just the commutative $\mathbb{C}$-algebra freely generated by the elements $X_i^{\pm 1}$, so where can define, starting from $\rho_\lambda$, a (irreducible) representation $\rho_\lambda^1$ of $\mathcal{T}_\lambda^1$ in $\mathbb{C}$ simply by defining

$$\rho_\lambda^1(X_i) := x_i \in \mathrm{End}(\mathbb{C})$$

The representation $\rho_\lambda^1$ is called the *non-quantum shadow* of $\rho_\lambda$.

If $\rho = \{\rho_\lambda \colon \mathcal{T}_\lambda^q \to \mathrm{End}(V_\lambda)\}$ is a local representation of the whole quantum Teichmüller space $\mathcal{T}_S^q$, we define the *non-quantum shadow* of $\rho$ as the collection

$$\rho^1 := \{\rho_\lambda^1 \colon \mathcal{T}_\lambda^1 \to \mathrm{End}(\mathbb{C})\}_{\lambda \in \Lambda(S)}$$

where $\rho_\lambda^1 \colon \mathcal{T}_\lambda^1 \to \mathrm{End}(\mathbb{C})$ is the non-quantum shadow of $\rho_\lambda$ for every $\lambda \in \Lambda(S)$.

**Theorem 5.1.** Given $\rho = \{\rho_\lambda \colon \mathcal{T}_\lambda^q \to \mathrm{End}(V_\lambda)\}_{\lambda \in \Lambda(S)}$ a local representation of the quantum Teichmüller space $\mathcal{T}_S^q$, the non-quantum shadow $\rho^1 = \{\rho_\lambda^1 \colon \mathcal{T}_\lambda^q \to \mathrm{End}(\mathbb{C})\}_{\lambda \in \Lambda(S)}$ of $\rho$ is a representation of the non-quantum space $\mathcal{T}_S^1$.

Let $\mathscr{R}_{loc}(\mathcal{T}_S^q)$ denote the set of the isomorphism classes of local representations of the quantum Teichmüller space $\mathcal{T}_S^q$.

**Theorem 5.2.** Let $S$ be a surface and let $q \in \mathbb{C}^*$ be a primitive $N$-th root of $(-1)^{N+1}$. Then, the application
$$\begin{array}{ccc} \mathscr{R}_{loc}(\mathcal{T}_S^q) & \longrightarrow & \mathrm{Repr}(\mathcal{T}_S^1, \mathbb{C}) \\ [\rho] & \longmapsto & \rho^1 \end{array}$$
that sends an isomorphism class of a local representation $\rho$ in its non-quantum shadow $\rho^1$ is well defined and onto. Moreover, the fibre on every element of $\mathrm{Repr}(\mathcal{T}_S^1, \mathbb{C})$ is composed of $N$ classes in $\mathscr{R}_{loc}(\mathcal{T}_S^q)$ and each element of the fibre on $\rho^1$ is determined by the choice an $N$-th root of the $x_1 x_2 \cdots x_n$, where $x_i = \rho_\lambda^1(X_i)$, for a certain $\lambda \in \Lambda(S)$.

We refer to [BBL07] for the definitions of *pleated surfaces* and its corresponding *enhanced homomorphism*.

**Theorem 5.3.** There exists a bijection between the set of conjugation classes of peripherally generic enhanced homomorphisms $(r, \{z_\pi\}_{\pi \in \Pi})$, from $\pi_1(S)$ to $\mathrm{Isom}^+(\mathbb{H}^3)$, and the set of non-quantum representations $\rho^1 = \{\rho_\lambda^1 \colon \mathcal{T}_\lambda^q \to \mathrm{End}(\mathbb{C})\}_{\lambda \in \Lambda(S)}$ of the non-quantum Teichmüller space $\mathcal{T}_S^1$, which sends the conjugation class of a peripherally generic enhanced homomorphism $(r, \{z_\pi\}_{\pi \in \Pi})$, in the non-quantum representation of $\mathcal{T}_S^1$ in which, for every $\lambda \in \Lambda(S)$, $\rho_\lambda^1$ is defined by the relation
$$\rho_\lambda^1(X_i) = x_i \, id_\mathbb{C}$$
for every $X_i$ generator of $\mathcal{T}_\lambda^q$, where $x_i$ denotes the shear-bend coordinates associated with $\lambda_i$ for a certain pleated surface $(\widetilde{f}, r)$ whose corresponding enhanced homomorphism is $(r, \{z_\pi\}_{\pi \in \Pi})$.

**Theorem 5.4.** Let $S$ be a surface and let $q \in \mathbb{C}^*$ be a primitive $N$-th root of $(-1)^{N+1}$. Then, the application
$$\begin{array}{ccc} \mathscr{R}_{loc}(\mathcal{T}_S^q) & \longrightarrow & \mathscr{E}\mathscr{H}(S) \\ [\rho] & \longmapsto & [r, \{z_\pi\}_\pi] \end{array}$$
that sends an isomorphism class of a local representation $\rho$ in its hyperbolic shadow $[r, \{z_\pi\}_\pi]$ is well defined and onto. Moreover, the fibre on every element of $\mathscr{E}\mathscr{H}(S)$ is composed of $N$ classes in $\mathscr{R}_{loc}(\mathcal{T}_S^q)$. Fixed $\lambda \in \Lambda(S)$, each element of the fibre on $\rho^1$ is determined by the choice an $N$-th root of the $x_1 x_2 \cdots x_n$, where the $x_i$ are the shear-bend coordinates associated with a certain pleated surface $(\widetilde{f}_\lambda, r)$ with pleating locus $\lambda$ realizing $[r, \{z_\pi\}_\pi]$ as enhanced homomorphism.

## 5.2 The definition of the invariants

We firstly need to define actions of the mapping class group $\mathcal{MCG}(S)$ of $S$ on the sets $\mathscr{E}\mathscr{H}(S)$ and $\mathrm{Repr}_{loc}(\mathcal{T}_S^q)$, which will be very useful hereinafter.

Let $[r, \{z_\pi\}_{\pi \in \Pi}]$ be a conjugation class of a peripherally generic enhanced homomorphism and let $\varphi \colon S \to S$ be a diffeomorphism. We define
$$[r, \{z_\pi\}_{\pi \in \Pi}] \cdot \varphi$$
as the conjugation class of a peripherally generic enhanced homomorphism $(s, \{\xi_\pi\}_{\pi \in \Pi})$ defined as follows:

- $s$ is equal to the composition $r \circ \varphi_*$, where $\varphi_* \colon \pi_1(S) \to \pi_1(S)$ is the isomorphism induced by $\varphi$ for an arbitrary choice of a path joining the base point of $S$ to its image under $\varphi$;

- for every $\pi \in \Pi$ $\xi_\pi$ is equal to $z_{\varphi_*(\pi)}$.

The conjugation class $[s, \{\xi_\pi\}_{\pi \in \Pi}]$ does not depend on the choices of the representative $(r, \{z_\pi\}_{\pi \in \Pi})$ and the path joining the base point of $S$ to its image under $\varphi$, so this construction defines a right action of $\mathcal{MCG}(S)$ on the set $\mathscr{EH}(S)$

$$\begin{array}{rcl} \mathscr{EH}(S) \times \mathcal{MCG}(S) & \longrightarrow & \mathscr{EH}(S) \\ ([r, \{\xi_\pi\}_\pi], [\varphi]) & \longmapsto & [r, \{z_\pi\}_\pi] \cdot \varphi \end{array}$$

This concludes the definition of the action on $\mathscr{EH}(S)$.

If $\varphi$ is a diffeomorphism of $S$, then for every $\lambda \in \Lambda(S)$ we denote by $\varphi(\lambda) \in \Lambda(S)$ the ideal triangulation defined by $\varphi(\lambda)_i := \varphi(\lambda_i)$ for every $i = 1, \ldots, n$. If $S_0$ and $S_0'$ are respectively the surfaces obtained from $S$ by splitting it along the triangulations $\lambda$ and $\varphi(\lambda)$, then $\varphi$ induces also a diffeomorphism $\overline{\phi}$ from $S_0$ to $S_0'$. $S_0$ and $S_0'$ can be endowed with unique ideal triangulations $\lambda_0$ and $\lambda_0'$, which clearly verify $\overline{\varphi}(\lambda_0) = \lambda_0'$. Consequently there is a natural algebra isomorphism $\overline{\varphi}_{\lambda_0}^q \colon \mathcal{T}_{\lambda_0}^q \to \mathcal{T}_{\lambda_0'}^q$ sending, for every $i$, the variable in $\mathcal{T}_{\lambda_0}^q$ corresponding to $(\lambda_0)_j$ in the variable in $\mathcal{T}_{\lambda_0}^q$ corresponding via $\overline{\phi}$ to $(\lambda_0')_j$. This implies that every local representation $\rho_{\varphi(\lambda)} \colon \mathcal{T}_{\varphi(\lambda)}^q \to \mathrm{End}(W)$ induces a local representation $\rho_\lambda' \colon \mathcal{T}_\lambda^q \to \mathrm{End}(W)$ defined as follows: fixed a representation $\zeta_0 \in \mathscr{F}_{S_0'}(\rho_{\varphi(\lambda)})$, we take the local representation $\rho_\lambda'$ represented by

$$\zeta_0 \circ \overline{\varphi}_\lambda^q \colon \mathcal{T}_{\lambda_0}^q \longrightarrow \mathrm{End}(W)$$

The local representation $\rho_\lambda'$ just constructed does not depend on the choice of the representative $\zeta_0$, we will label it as $\rho_{\varphi(\lambda)} \cdot \varphi$.

Moreover, the diffeomorphism $\varphi$ induces also, for every $\lambda \in \Lambda(S)$, an algebra isomorphism $\varphi_\lambda^q \colon \mathcal{T}_\lambda^q \to \mathcal{T}_{\varphi(\lambda)}^q$ sending, for every $i$, the variable in $\mathcal{T}_\lambda^q$ corresponding to $\lambda_i$ in the variable in $\mathcal{T}_{\varphi(\lambda)}^q$ corresponding to $\varphi(\lambda_i)$. It is clear that these isomorphisms have a good behaviour with respect to the coordinates changes $\Phi_{\lambda\lambda'}^q$. More precisely, they induce isomorphisms also on the fraction rings $\widehat{\varphi}_\lambda^q \colon \widehat{\mathcal{T}}_\lambda^q \to \widehat{\mathcal{T}}_{\varphi(\lambda)}^q$ in the obvious way and the following relation holds

$$\widehat{\varphi}_\lambda^q \circ \Phi_{\lambda\lambda'}^q = \Phi_{\varphi(\lambda)\varphi(\lambda')}^q \circ \widehat{\varphi}_{\lambda'}^q \tag{22}$$

for every $\lambda, \lambda' \in \Lambda(S)$.

Given $\rho = \{\rho_\lambda \colon \mathcal{T}_\lambda^q \to \mathrm{End}(V_\lambda)\}_{\lambda \in \Lambda(S)}$ a local representation of the quantum Teichmüller space $\mathcal{T}_S^q$, we define a collection of local representations

$$\rho \cdot \varphi = \{(\rho \cdot \varphi)_\lambda \colon \mathcal{T}_\lambda^q \to \mathrm{End}(V_{\varphi(\lambda)})\}_{\lambda \in \Lambda(S)}$$

by setting $(\rho \cdot \varphi)_\lambda := \rho_{\varphi(\lambda)} \cdot \varphi$ for every $\lambda \in \Lambda(S)$, where $\rho_{\varphi(\lambda)} \cdot \varphi$ is constructed with the procedure described above. It is simple to see that the relation 22 implies that the collection $\rho \cdot \varphi$ is indeed a local representation of the quantum Teichmüller space $\mathcal{T}_S^q$. Hence we have described also a right action

$$\begin{array}{rcl}
\mathrm{Repr}_{loc}(\mathcal{T}_S^q) \times \mathcal{MCG}(S) & \longrightarrow & \mathrm{Repr}_{loc}(\mathcal{T}_S^q) \\
(\rho, [\varphi]) & \longmapsto & \rho \cdot \varphi
\end{array}$$

In addition, observe that the isomorphisms $\varphi_\lambda^q$ send the central element $H_\lambda$ of $\mathcal{T}_\lambda^q$ in the central element $H_{\varphi(\lambda)}$ of $\mathcal{T}_{\varphi(\lambda)}^q$. Therefore, the central load of the representation $\rho \cdot \varphi$ is the same of the one of $\rho$.

Recall that in Theorem 5.4 we have shown the existence of a surjective map

$$\begin{array}{rcl}
\Theta: \; \mathscr{R}_{loc}(\mathcal{T}_S^q) & \longrightarrow & \mathscr{EH}(S) \\
[\rho] & \longmapsto & [r, \{z_\pi\}_\pi]
\end{array}$$

that sends each isomorphism class $[\rho]$ in the conjugation class of its hyperbolic shadow $[r, \{z_\pi\}_\pi]$ and that has the fibre on $[r, \{z_\pi\}_\pi]$ composed of the $N$ isomorphisms classes of representations, one for every possible central load, which are the $N$-roots of $x_1 \cdots x_n$.

**Lemma 5.5.** The following relation holds

$$\Theta([\rho] \cdot \varphi) = \Theta([\rho]) \cdot \varphi$$

where $\varphi$ is acting on $\mathscr{R}_{loc}(\mathcal{T}_S^q)$ in the first member and on $\mathscr{EH}(S)$ in the second.

*Proof.* Let $\widetilde{\varphi}: \widetilde{S} \to \widetilde{S}$ a certain lift of $\varphi$ on the universal covering. Recalling Theorem 5.4, if $\rho$ is equal to a collection $\{\rho_\lambda \colon \mathcal{T}_\lambda^q \to \mathrm{End}(V_\lambda)\}_{\lambda \in \Lambda(S)}$ and $[r, \{z_\pi\}_{\pi_\Pi}]$ is its hyperbolic shadow, then for every $\lambda \in \Lambda(S)$ we can find a pleated surface with bending locus $\lambda$ associated with $\rho_\lambda$ that looks like $(\widetilde{f}_\lambda, r)$. Then, by inspection of the shear-bend coordinates, the pleated surface $(\widetilde{f}_{\varphi(\lambda)} \circ \widetilde{\varphi}, r \circ \varphi_*)$ is a pleated surface with bending locus $\lambda$ associated with the representation $\rho_{\varphi(\lambda)} \cdot \varphi$. Following the definitions, this fact implies the assertion. □

Let $\varphi \colon S \to S$ be a diffeomorphism of the surface $S$. Denote by $M_\varphi$ the mapping torus of $\varphi$, which is the 3-manifold obtained as quotient of $S \times \mathbb{R}$ by the group of diffeomorphisms generated by

$$\begin{array}{rcl}
\tau_\varphi \colon \; S \times \mathbb{R} & \longrightarrow & S \times \mathbb{R} \\
(p, t) & \longmapsto & (\varphi(p), t+1)
\end{array}$$

We define also the inclusion

$$\begin{array}{rcl}
i \colon \; S & \longrightarrow & M_\varphi \\
p & \longmapsto & [p, 0]
\end{array}$$

Observe that the homomorphism $i_* \colon \pi_1(S) \to \pi_1(M_\varphi)$ induced by $i$ is injective. Indeed, assume that $\gamma \colon S^1 \to S$ is a closed path such that $i \circ \gamma$ it homotopically trivial in $M_\varphi$. Then there exists a map $f \colon D^2 \to M_\varphi$ such that $f|_{S^1} = i \circ \gamma$. Because $D^2$ is simply connected, we can lift $f$ to an application $\widetilde{f} \colon D^2 \to S \times \mathbb{R}$ such that

- $\pi \circ \widetilde{f} = f$, where $\pi \colon S \times \mathbb{R} \to M_\varphi$ is the projection map;

- $\widetilde{f}(S^1) \subset S \times \{0\}$.

By construction we have $\widetilde{f}|_{S^1} = \gamma \times \{0\}\colon S^1 \to S \times \{0\}$. If $p\colon S \times \mathbb{R} \to S \times \{0\}$ is the obvious projection on the first component, then $p \circ \widetilde{f}$ is a continuous map from $D^2$ to $S \times \{0\}$ whose restriction to $S^1$ is equal to $\gamma \times \{0\}$. This proves that $\gamma$ was homotopically trivial also in $S$, hence that $i_*$ is injective. Moreover, by construction of the mapping torus, the maps $i$ and $i \circ \varphi$ are homotopic, hence $i_* = i_* \circ \varphi_*$.

By virtue of the Thurston's Hyperbolization Theorem, the mapping torus $M_\varphi$ admits a complete finite-volume hyperbolic structure if and only if $\varphi$ is isotopic to a pseudo-Anosov diffeomorphism (see [Thu83] and [Ota96]). Assume that $\varphi$ is pseudo-Anosov and let $r\colon \pi_1(M_\varphi) \to \mathbb{P}\mathrm{SL}_2(\mathbb{C})$ be the holonomy of the complete structure on $M_\varphi$, which is unique up to conjugation in $\mathbb{P}\mathrm{SL}_2(\mathbb{C})$ because of the Mostow Rigidity Theorem. From $r$ we obtain the homomorphism $r_\varphi\colon \pi_1(S) \to \mathbb{P}\mathrm{SL}_2(\mathbb{C})$, defined as the composition $r \circ i_*$. Thanks to the injectivity of $i_*$, we have that $r_\varphi$ is also injective. Moreover, because $r$ is the holonomy of a complete finite-volume structure, $r_\varphi$ leads every peripheral subgroup $\pi$ of $\pi_1(S)$ in a parabolic subgroup of $\mathbb{P}\mathrm{SL}_2(\mathbb{C})$. This means that $r_\varphi$ admits a unique enhancement $\{z_\pi\}_{\pi \in \Pi}$. Observe also that $(r_\varphi, \{z_\pi\}_{\pi \in \Pi})$ is peripherally generic, because of [BBL07, Lemma 14]. Therefore a pseudo-Anosov diffeomorphism $\varphi$ produces a unique conjugation class of enhanced homomorphisms from $\pi_1(S)$ to $\mathbb{P}\mathrm{SL}_2(\mathbb{C})$ peripherally generic $[r_\varphi, \{z_\pi\}_{\pi \in \Pi}] \in \mathscr{EH}(S)$. By virtue of Theorem 5.4, there exists a unique non-quantum representation $\rho^1_\varphi$ of $\mathcal{T}^1_S$ related to $[r_\varphi, \{z_\pi\}_{\pi \in \Pi}]$.

**Lemma 5.6.** *The non-quantum representation $\rho^1_\varphi$ verifies*

$$(\rho^1_\varphi)_\lambda(H) = (\rho^1_\varphi)_\lambda(P_i) = id_\mathbb{C}$$

*for every $\lambda \in \Lambda(S)$ and for every $i = 1, \ldots, p$.*

*Proof.* See [BL07, Lemma 39]. □

As consequence of this lemma, the isomorphism classes of local representations in $\Theta^{-1}([r_\varphi, \{z_\pi\}_{\pi \in \Pi}])$ are classified by their central loads, which are all the possible $N$-th roots of unity. So, for every $k \in \mathbb{Z}_N$, we have found a unique isomorphism class of local representations $[\rho^k_\varphi]$ having $[r_\varphi, \{z_\pi\}_{\pi \in \Pi}]$ as hyperbolic shadow and $q^{2k}$ as central load. Because $i_* \circ \varphi_* = i_*$ and because $r_\varphi$ admits a unique enhancement, we have

$$[r_\varphi, \{z_\pi\}_{\pi \in \Pi}] \cdot \varphi = [r_\varphi, \{z_\pi\}_{\pi \in \Pi}]$$

In other words, the element $[r_\varphi, \{z_\pi\}_{\pi \in \Pi}]$ in $\mathscr{EH}(S)$ is fixed by the right action of $\varphi$. Hence we deduce, recalling Lemma 5.5, that the action of $\varphi$ on $\mathscr{R}_{loc}(S)$ preserves the set $\Theta^{-1}([r_\varphi, \{z_\pi\}_{\pi \in \Pi}])$. Furthermore, we have seen that the isomorphisms $\varphi^q_\lambda$ preserve the central element $H$, so the action of $\varphi$ on $\Theta^{-1}([r_\varphi, \{z_\pi\}_{\pi \in \Pi}])$ is necessarily the trivial one. We have proved

**Theorem 5.7.** *Let $S$ be a closed surface with punctures. Fix $q$ a $N$-th primitive root of $(-1)^{N+1}$, $k \in \mathbb{Z}_N$ and $\varphi\colon S \to S$ a pseudo-Anosov diffeomorphism. Then there exists a conjugation class of peripherally generic enhanced homomorphism $[r_\varphi, \{z_\pi\}_{\pi \in \Pi}]$, uniquely determined by $\varphi$, which is fixed by the action of $\varphi$ on* $\mathscr{EH}(S)$. Moreover, the set $\Theta^{-1}([r_\varphi, \{z_\pi\}_{\pi \in \Pi}])$ is composed of exactly such isomorphisms classes $[\rho_\varphi^k]$ of local representations of $\mathcal{T}_S^q$ having $[r_\varphi, \{z_\pi\}_{\pi \in \Pi}]$ as hyperbolic shadow and $q^{2k}$ as central load, and each $[\rho_\varphi^k]$ is fixed by the action of $\varphi$ on $\mathscr{R}_{loc}(S)$.

Henceforth, the numbers $q \in \mathbb{C}^*$ and $k \in \mathbb{Z}_N$ will be fixed. Choosing a representative $\rho$ of the class $[\rho_\varphi^k]$, Theorem 5.7 implies that the representations $\rho$ and $\rho \cdot \varphi$ are isomorphic as local representations of $\mathcal{T}_S^q$. Therefore, we can consider the family $\{(\mathscr{L}_{\lambda\lambda'}^{\rho \cdot \varphi\, \rho}, \psi_{\lambda\lambda'}^{\rho \cdot \varphi\, \rho})\}$ of intertwining operators between $\rho \cdot \varphi$ and $\rho$.

In particular, fixing $\lambda \in \Lambda(S)$, we can consider the element

$$(\mathscr{L}_{\lambda\,\varphi(\lambda)}^{\rho \cdot \varphi\, \rho}, \psi_{\lambda\,\varphi(\lambda)}^{\rho \cdot \varphi\, \rho})$$

Observe that both the representations $(\rho \cdot \varphi)_\lambda$, $\rho_{\varphi(\lambda)}$ have values in $\mathrm{End}(V_{\varphi(\lambda)})$, so the set $\mathscr{L}_{\lambda\,\varphi(\lambda)}^{\rho \cdot \varphi\, \rho}$ is contained in $\mathrm{GL}(V_{\varphi(\lambda)})$.

The element $(\mathscr{L}_{\lambda\,\varphi(\lambda)}^{\rho \cdot \varphi\, \rho}, \psi_{\lambda\,\varphi(\lambda)}^{\rho \cdot \varphi\, \rho})$ depends on the chosen representative $\rho$ of the isomorphism class $[\rho_\varphi^k]$ and on the ideal triangulation $\lambda$. We want to produce a more intrinsic object.

**Definition 5.8.** Let $\rho = \{\rho_\lambda \colon \mathcal{T}_\lambda^q \to \mathrm{End}(V_\lambda)\}_{\lambda \in \Lambda(S)}$ be a local representation of $\mathcal{T}_S^q$. Fix $\mu \in \Lambda(S)$ an ideal triangulation and a tensor-split linear isomorphism

$$M = M_1 \otimes \cdots \otimes M_m \colon \bigotimes_{j=1}^m V_{\mu,j} = V_\mu \longrightarrow W = \bigotimes_{j=1}^m W_j$$

Then we denote by $M \bullet_\mu \rho$ the local representation defined as follows

$$(M \bullet_\mu \rho)_\lambda := \begin{cases} \rho_\lambda \colon \mathcal{T}_\lambda^q \to \mathrm{End}(V_\lambda) & \text{if } \lambda \neq \mu \\ M \circ \rho_\mu(\cdot) \circ M^{-1} \colon \mathcal{T}_\mu^q \to \mathrm{End}(W) & \text{if } \lambda = \mu \end{cases}$$

Observe that $M \circ \rho_\mu(\cdot) \circ M^{-1} \colon \mathcal{T}_\mu^q \to \mathrm{End}(W)$ is a local representation because $M$ is tensor splitting.

**Lemma 5.9.** Let $\rho = \{\rho_\lambda \colon \mathcal{T}_\lambda^q \to \mathrm{End}(V_\lambda)\}_{\lambda \in \Lambda(S)}$ be a local representation of the quantum Teichmüller space $\mathcal{T}_S^q$. The following hold:

1. let $\mu \in \Lambda(S)$ be an ideal triangulation and $M \colon V_\mu \to W$ a tensor-split isomorphism. Then $M$ belongs to $\mathscr{L}_{\mu\mu}^{M \bullet_\mu \rho\; \rho}$;

2. let $\varphi \colon S \to S$ be a diffeomorphism, $\lambda \in \Lambda(S)$ be an ideal triangulation and $M \colon V_{\varphi(\lambda)} \to W$ a tensor-split isomorphism. Then

$$M \bullet_\lambda (\rho \cdot \varphi) = (M \bullet_{\varphi(\lambda)} \rho) \cdot \varphi$$

*Proof.* We will focus on one point at time.

1. As usual, we denote by $S_0$ the surface obtained by splitting $S$ along $\mu$ and by $\mu_0$ the ideal triangulation on $S_0$ induced by $\mu$. Fixed an element $\zeta \in \mathscr{F}_{S_0}(\rho_\mu)$, a representative of $(M \bullet_\mu \rho)_\mu$ is given by $\zeta' = M \circ \zeta(\cdot) \circ M^{-1}$ (it is a well defined representation of $\mathcal{T}_{\mu_0}^q$ because $M$ is tensor-split). Clearly $M$ sends $\zeta$ in $\zeta'$. Hence we have that the couple $(\zeta', \zeta)$ belongs to $\mathscr{F}_{S_0}((M \bullet_\mu \rho)_\mu) \times \mathscr{F}_{S_0}(\rho_\mu)$, so we have found an element $[\zeta', \zeta] \in \mathscr{A}_{\mu\mu}^{M\bullet_\mu\rho\ \rho}$ that corresponds to $M$, which clearly implies $M \in \mathscr{L}_{\mu\mu}^{M\bullet_\mu\rho\ \rho}$.

2. Denote by $S'_0$ and $S''_0$ respectively the surfaces obtained by splitting $S$ along the ideal triangulations $\lambda$ and $\varphi(\lambda)$, endowed with the ideal triangulations $\lambda_0$ and $\varphi(\lambda_0)$. Recall that the diffeomorphism $\varphi$ induces an isomorphism $\overline{\varphi}_{\lambda_0}^q$ from $\mathcal{T}_{\lambda_0}^q$ to $\mathcal{T}_{\varphi(\lambda_0)}^q$. Fixed $\zeta \in \mathscr{F}_{S''_0}(\rho_{\varphi(\lambda)})$, a representative of the representation $(M \bullet_\lambda (\rho \cdot \varphi))_\lambda$ is given by

$$M \circ (\zeta \cdot \varphi)(\cdot) \circ M^{-1} = M \circ (\zeta \circ \overline{\varphi}_{\lambda_0}^q)(\cdot) \circ M^{-1}$$
$$= M \circ (\zeta(\overline{\varphi}_{\lambda_0}^q(\cdot)) \circ M^{-1}$$
$$= (M \circ \zeta(\cdot) \circ M^{-1}) \cdot \varphi$$

Now observe that $M \circ \zeta(\cdot) \circ M^{-1}$ is a representative of $(M \bullet_{\varphi(\lambda)} \rho)_{\varphi(\lambda)}$ and consequently $(M \circ \zeta(\cdot) \circ M^{-1}) \cdot \varphi$ is a representative of

$$(M \bullet_{\varphi(\lambda)} \rho)_{\varphi(\lambda)} \cdot \varphi = ((M \bullet_{\varphi(\lambda)} \rho) \cdot \varphi)_\lambda$$

Hence we have proved

$$(M \bullet_\lambda (\rho \cdot \varphi))_\lambda = ((M \bullet_{\varphi(\lambda)} \rho) \cdot \varphi)_\lambda$$

It is simple to see that on all the other triangulations these representations are obviously equal, because they coincide with $\rho \cdot \varphi$ on them, so the proof of the second assertion is done. □

Recall that we are interested in the construction of an intrinsic object starting from

$$(\mathscr{L}_{\lambda\,\varphi(\lambda)}^{\rho\cdot\varphi\ \rho}, \psi_{\lambda\,\varphi(\lambda)}^{\rho\cdot\varphi\ \rho})$$

In order to do so, we choose a tensor-split isomorphism $M$ from $V_{\varphi(\lambda)}$ to a fixed vector space $W = \bigotimes_i W_i$. For example, a natural choice of $W$ is

$$W := \bigotimes_{i=1}^m (\mathbb{C}^N)_i$$

Now we modify the representations $\rho \cdot \varphi$ and $\rho$ by taking $M \bullet_\lambda (\rho \cdot \varphi)$ and $M \bullet_{\varphi(\lambda)} \rho$. Considering now the object

$$i(\varphi; \lambda, \rho, M) := (\mathscr{L}_{\lambda\,\varphi(\lambda)}^{M\bullet_\lambda(\rho\cdot\varphi)\ M\bullet_{\varphi(\lambda)}\rho}, \psi_{\lambda\,\varphi(\lambda)}^{M\bullet_\lambda(\rho\cdot\varphi)\ M\bullet_{\varphi(\lambda)}\rho})$$

we obtain a set of isomorphisms in $GL(W)$ and so automorphisms of a vector space that is independent from the choices done. The following study will be focused on the dependence of $i(\varphi; \lambda, \rho, M)$ on

- the representative $\rho$ of the isomorphism class $[\rho_\varphi^k]$;
- the tensor split isomorphism $M \colon V_{\varphi(\lambda)} \to W$;
- the ideal triangulation $\lambda \in \Lambda(S)$.

Observe that the set $\mathscr{L}_{\lambda\lambda'}^{\rho\rho'}$ does depend only on the local representations $\rho_\lambda$ and $\rho'_{\lambda'}$. In particular, in the case of $i(\varphi;\lambda,\rho,M)$, the only local representations involved are $M \circ \rho_{\varphi(\lambda)}(\cdot) \circ M^{-1}$ and $M \circ (\rho_{\varphi(\lambda)} \circ \varphi_\lambda^q)(\cdot) \circ M^{-1}$. Therefore, fixed $\lambda \in \Lambda(S)$, $i(\varphi;\lambda,\rho,M)$ depends only on $M$ and $\rho_{\varphi(\lambda)}$. By virtue of the first point of Lemma 5.9, the representation $M \bullet_{\varphi(\lambda)} \rho$ is isomorphic to $\rho$, in particular it belongs to $[\rho_\varphi^k]$, and it has the property that $(M \bullet_{\varphi(\lambda)} \rho)_{\varphi(\lambda)}$ has values in $\mathrm{End}(W)$. Moreover, if $\rho'$ is a representation of $\mathcal{T}_S^q$ isomorphic to $\rho$ such that $\rho'_{\varphi(\lambda)}$ has values in $\mathrm{End}(W)$, then $\rho'_{\varphi(\lambda)}$ is equal to $\rho_{\varphi(\lambda)} \bullet_{\varphi(\lambda)} M'$ for some $M': V_{\varphi(\lambda)} \to W$ (it is sufficient to take $M' \in \mathscr{L}_{\varphi(\lambda)\,\varphi(\lambda)}^{\rho'\rho}$). The second point of Lemma 5.9 tells us that the representation $M \bullet_\lambda (\rho \cdot \varphi)$ is just the image under the action of $\varphi$ of $M \bullet_{\varphi(\lambda)} \rho$. Putting together these observations, we have proved that every $i(\varphi;\lambda,\rho,M)$ is equal to $i(\varphi;\lambda, M \bullet_{\varphi(\lambda)} \rho, id)$. Moreover, for every $\rho'$ and $\rho''$ representations in $[\rho_\varphi^k]$ such that $\rho'_{\varphi(\lambda)}$ and $\rho''_{\varphi(\lambda)}$ are equal and have values in $\mathrm{End}(W)$, we have

$$i(\varphi;\lambda,\rho',id) = i(\varphi;\lambda,\rho'',id)$$

Henceforth, fixed $\lambda \in \Lambda(S)$ and $W$, the study of the objects $i(\varphi;\lambda,\rho,M)$ can be reduced to the investigations of the elements

$$j(\varphi;\lambda,\rho) := i(\varphi;\lambda,\rho,id) = (\mathscr{L}_{\lambda\,\varphi(\lambda)}^{\rho\cdot\varphi\,\rho}, \psi_{\lambda\,\varphi(\lambda)}^{\rho\cdot\varphi\,\rho})$$

where $\rho$ is a local representation of $\mathcal{T}_S^q$ having $\rho_{\varphi(\lambda)}$ with values in $\mathrm{End}(W)$.

**Lemma 5.10.** *Let $\lambda, \lambda' \in \Lambda(S)$ be two ideal triangulations of $S$ and $\rho, \rho'$ two local representations of $\mathcal{T}_S^q$. Then for every $L \in \mathscr{L}_{\varphi(\lambda')\,\varphi(\lambda)}^{\rho'\rho}$ the application*

$$\begin{array}{rcl} f_L: \mathscr{L}_{\lambda\,\varphi(\lambda)}^{\rho\cdot\varphi\,\rho} & \longrightarrow & \mathscr{L}_{\lambda'\,\varphi(\lambda')}^{\rho'\cdot\varphi\,\rho'} \\ A & \longmapsto & L \circ A \circ L^{-1} \end{array}$$

*is well defined and bijective. Moreover, it respects the actions, i.e. $f_L(c \cdot A) = c \cdot f_L(A)$ for every $A \in \mathscr{L}_{\lambda\,\varphi(\lambda)}^{\rho\cdot\varphi\,\rho}$ and $c \in H_1(S; \mathbb{Z}_N)$.*

*Proof.* Because $L$ is an element of $\mathscr{L}_{\varphi(\lambda')\,\varphi(\lambda)}^{\rho'\rho}$, the composition $A \circ L^{-1}$ belongs to $\mathscr{L}_{\lambda\,\varphi(\lambda')}^{\rho\cdot\varphi\,\rho'}$. If $L$ belongs also $\mathscr{L}_{\lambda'\lambda}^{\rho'\cdot\varphi\,\rho\cdot\varphi}$, then the assertion is an obvious corollary of the Composition property of the family $\{(\mathscr{L}_{\lambda\lambda'}^{\rho\rho'}, \psi_{\lambda\lambda'}^{\rho\rho'})\}$. It is sufficient to show that $L \in \mathscr{L}_{\lambda'\lambda}^{\rho'\cdot\varphi\,\rho\cdot\varphi}$ when $\lambda, \lambda'$ differ by an elementary move. We will investigate the diagonal exchange case $\lambda' = \Delta_i(\lambda)$, the other possibilities are simpler to study.

Firstly we need to introduce the notations: we denote by $S_0$ the surface obtained by splitting $S$ along all the edges of $\lambda$ except for $\lambda_i$ and we label its ideal triangulations as $\lambda_0$ and $\lambda'_0$, the first corresponding to $\lambda \in \Lambda(S)$ and the second to $\lambda' \in \Lambda(S)$. Moreover, $\varphi(S_0)$ will be the surface obtained by splitting $S$ along all the edges of $\varphi(\lambda)$ except for $\varphi(\lambda_i)$ and $\varphi(\lambda_0), \varphi(\lambda'_0)$ its ideal triangulations induced respectively by $\varphi(\lambda)$ and $\varphi(\lambda')$.

Because $L$ is in $\mathscr{L}_{\varphi(\lambda')\,\varphi(\lambda)}^{\rho'\rho}$, there exists a couple

$$(\zeta', \zeta) \in \mathscr{F}_{\varphi(S_0)}(\rho'_{\varphi(\lambda')}) \times \mathscr{F}_{\varphi(S_0)}(\rho_{\varphi(\lambda)})$$

such that
$$L \circ \zeta(X) \circ L^{-1} = (\zeta' \circ \Phi^q_{\varphi(\lambda'_0)\varphi(\lambda_0)})(X)$$
for every $X \in \mathcal{T}^q_{\varphi(\lambda_0)}$. Because $\zeta$ represents the local representation $\rho_{\varphi(\lambda)}$, by composing it with the natural algebra isomorphism $\overline{\varphi}^q_{\lambda_0}\colon \mathcal{T}^q_{\lambda_0} \to \mathcal{T}^q_{\varphi(\lambda_0)}$ we obtain a representative of the local representation $(\rho \cdot \varphi)_\lambda$. The same argument shows that $\zeta' \circ \overline{\varphi}^q_{\lambda'_0}$ is a representative of $(\rho' \cdot \varphi)_{\lambda'}$. This means that the couple $(\zeta' \circ \overline{\varphi}^q_{\lambda'_0}, \zeta \circ \overline{\varphi}^q_{\lambda_0})$ belongs to

$$\mathscr{F}_{S_0}((\rho' \cdot \varphi)_{\lambda'}) \times \mathscr{F}_{S_0}((\rho \cdot \varphi)_\lambda)$$

Moreover, the following holds
$$\begin{aligned}L \circ (\zeta \circ \overline{\varphi}^q_{\lambda_0})(Y) \circ L^{-1} &= L \circ \zeta(\overline{\varphi}^q_{\lambda_0}(Y)) \circ L^{-1} \\ &= (\zeta' \circ \Phi^q_{\varphi(\lambda'_0)\varphi(\lambda_0)})(\overline{\varphi}^q_{\lambda_0}(Y)) \\ &= (\zeta' \circ \overline{\varphi}^q_{\lambda'_0} \circ \Phi^q_{\lambda'_0\lambda_0})(Y)\end{aligned}$$

So the equivalence class $[\zeta' \circ \overline{\varphi}^q_{\lambda'_0}, \zeta \circ \overline{\varphi}^q_{\lambda_0}]$ belongs to $\mathscr{A}^{\rho'\cdot\varphi\ \rho\cdot\varphi}_{\lambda'\lambda}$ and corresponds in $\mathscr{L}^{\rho'\cdot\varphi\ \rho\cdot\varphi}_{\lambda'\lambda}$ exactly to the isomorphism $L$. This shows that $L \in \mathscr{L}^{\rho'\cdot\varphi\ \rho\cdot\varphi}_{\lambda'\lambda}$ and so the assertion. $\square$

For every $\lambda \in \Lambda(S)$ and for every local representation $\rho$ such that $\rho_{\varphi(\lambda)}$ has values in $\mathrm{End}(W)$, the couples $(\mathscr{L}^{\rho\cdot\varphi\ \rho}_{\lambda\,\varphi(\lambda)}, \psi^{\rho\cdot\varphi\ \rho}_{\lambda\,\varphi(\lambda)})$ verifies

- the set $\mathscr{L}^{\rho\cdot\varphi\ \rho}_{\lambda\,\varphi(\lambda)}$ in contained in $\mathrm{GL}(W)$;

- $H_1(S; \mathbb{Z}_N)$ acts on $\mathscr{L}^{\rho\cdot\varphi\ \rho}_{\lambda\,\varphi(\lambda)}$ by right multiplication via certain tensor-split isomorphisms of $W$, uniquely determined up to multiplication by an $N$-th root of unity.

We will say that two such couples are *conjugate* if there exists an element $L$ in $\mathrm{GL}(W)$ such that

- the conjugation map $f_L\colon \mathrm{GL}(W) \to \mathrm{GL}(W)$ that sends $A$ in $LAL^{-1}$ is a bijection between the sets;

- $f_L$ commutes with the actions.

Two couples are tensor-split conjugate if there exists an automorphism $L$ like above that is tensor-split.

**Theorem 5.11.** *The conjugacy class of a set $i(\varphi; \lambda, \rho, M)$ does not depend on the tensor-split isomorphism $M\colon V_{\varphi(\lambda)} \to W$, on the representation $\rho \in [\rho^k_\varphi]$ and on $\lambda \in \Lambda(S)$.*

*Proof.* We have seen that it is sufficient to study the sets $j(\varphi; \lambda, \rho)$ with $\rho_{\varphi(\lambda)}$ with values in $\mathrm{End}(W)$. Fixing $\lambda$ and varying $\rho$ in the sets of local representations in $[\rho^k_\varphi]$ with values in $\mathrm{End}(W)$, Lemma 5.10 with $\lambda = \lambda'$ shows that the tensor-split conjugacy class does not change (when $\lambda = \lambda'$ the isomorphism $L \in \mathscr{L}^{\rho'\rho}_{\varphi(\lambda)\varphi(\lambda)}$ is tensor-split). So the tensor-split conjugacy class of $j(\varphi; \lambda, \rho)$ depends only on $\lambda$.

Now choose two different ideal triangulations and two representations $\rho$, $\rho'$ such that $\rho_{\varphi(\lambda)}$ and $\rho'_{\varphi(\lambda')}$ have values in $\mathrm{End}(W)$. Because $i(\varphi; \lambda, \rho)$ does not depend on $\rho_{\varphi(\lambda')}$, we can assume that also $\rho_{\varphi(\lambda')}$ has values in $\mathrm{End}(W)$, by replacing $\rho$ with $M \bullet_{\varphi(\lambda')} \rho$ for a certain $M$. Analogously we can assume that $\rho'_{\varphi(\lambda)}$ has values in $\mathrm{End}(W)$. Now we have that each element $L$ of $\mathscr{L}^{\rho'\rho}_{\varphi(\lambda')\,\varphi(\lambda)}$ belongs to $\mathrm{GL}(W)$. Hence, by applying Lemma 5.10 on a fixed $L$ in $\mathscr{L}^{\rho'\rho}_{\varphi(\lambda')\,\varphi(\lambda)}$, we have that $j(\varphi; \lambda, \rho)$ and $j(\varphi; \lambda', \rho')$ are conjugated (not necessarily tensor-split conjugated) and this finally proves the announced result. $\square$

We will denote by $I(q, k, \varphi)$ this conjugacy class, depending only on the primitive $N$-th root of unity $q$, the number $k \in \mathbb{Z}_N$ and the pseudo-Anosov diffeomorphism $\varphi$.

Explicitly, in order to obtain the invariant of the diffeomorphism $\varphi$ with $q$ and $k$ fixed, we can proceed as follows

1. we fix an ideal triangulation $\lambda \in \Lambda(S)$ and a local representation $\rho \in [\rho^k_\varphi]$;

2. we possibly replace $\rho$ with a representation $\rho'$ such that $\rho_{\varphi(\lambda)}$ has values in $\mathrm{End}(W)$. We can also assume that $\rho_{\varphi(\lambda)}$ is in standard form. More precisely, we can choose $\rho$ such that every representative of $\rho_{\varphi(\lambda)}$ is the tensor product of triangle representations that are in standard form, as described in Remark 1.5. Observe that, if $\rho_{\varphi(\lambda)}$ is in standard position, the same holds for the representation $(\rho \cdot \varphi)_\lambda$;

3. we fix a sequence of ideal triangulations $\lambda = \lambda^{(0)}, \ldots, \lambda^{(k)} = \varphi(\lambda)$ leading from $\lambda$ to $\varphi(\lambda)$ and we find through it an element $L$ of $\mathscr{L}^{\rho\cdot\varphi\,\rho}_{\lambda\,\varphi(\lambda)}$. The other elements of $\mathscr{L}^{\rho\cdot\varphi\,\rho}_{\lambda\,\varphi(\lambda)}$ can be produced as $L \circ B(c)^{-1}$ for varying $c \in H_1(S; \mathbb{Z}_N)$, where $B(c)$ is an element of $\mathrm{GL}(W)$ as described in Lemma 4.3;

4. we take the conjugacy class of this set, which will have the following form
$$\{C \circ (H_1(S; \mathbb{Z}_N) \cdot L) \circ C^{-1} \mid C \in \mathrm{GL}(W)\}$$
where
$$C \circ (H_1(S; \mathbb{Z}_N) \cdot L) \circ C^{-1} = \{C \circ L \circ B(c)^{-1} \circ C^{-1} \mid c \in H_1(S; \mathbb{Z}_N)\}$$

This is the resulting invariant for the pseudo-Anosov diffeomorphism, having chosen $q$ a primitive $N$-th root of unity and $k$ a certain element of $\mathbb{Z}_N$.

# A   Proof of Proposition 1.14

The result stated in Proposition 1.14 has been announced firstly in [BBL07, Lemma 21] without proof. In order to prove this result we will just trace the procedure used in [BL07] to prove the much more intricate classification of irreducible representations for punctured closed surfaces.

In this Section, we will assume $S$ is an ideal polygon with $p$ vertices, i. e. a surface obtained from $D = D^2$ by removing $p$ punctures in $\partial D$, with $p \geq 3$. Let $\lambda \in \Lambda(S)$ be an ideal triangulation and let $\Gamma$ be the dual graph of $\lambda$. $\Gamma$ is a deformation retract of $S$, hence it is a tree by virtue of the simply connectedness of $S$. Moreover, the leaves of $\Gamma$ exactly correspond to those triangles having two edges lying in $\partial S$. Our purpose is to find a presentation of $(\mathcal{H}(\lambda; \mathbb{Z}), \sigma)$ analogous to the one in [BL07, Proposition 5], in order to simplify the study of $\mathcal{T}_\lambda^q$. Firstly we must deal with the simplest case, in which $S$ is just an ideal triangle:

*Remark* A.1. Let $T$ be an oriented ideal triangle, with a fixed indexing of the edges that proceeds in the opposite way of the one given by the orientation, as in Figure 1. In these notations, the bilinear form $\sigma$ is represented by the matrix

$$\begin{pmatrix} 0 & 1 & -1 \\ -1 & 0 & 1 \\ 1 & -1 & 0 \end{pmatrix}$$

By taking the basis $e_1' := e_1$, $e_2' := e_2$, $e_3' = e_1 + e_2 + e_3$, the matrix representing $\sigma$ becomes

$$\begin{pmatrix} 0 & 1 & 0 \\ -1 & 0 & 0 \\ 0 & 0 & 0 \end{pmatrix}$$

The surface $T$ clearly admits a unique ideal triangulation and its Chekhov-Fock algebra is

$$\frac{\mathbb{C}[X_1^{\pm 1}, X_2^{\pm 1}, X_3^{\pm 1}]}{(X_i X_{i+1} = q^2 X_{i+1} X_i \mid i \in \mathbb{Z}_3)}$$

Denote by $H$ the element $\underline{X}^{(1,1,1)} = q^{-1} X_1 X_2 X_3 \in \mathcal{T}_T^q$. Then $H$ belongs to the monomial center $\mathcal{Z}_T^q$ and, if $q^2$ is a primitive $N$-th root of unity, the same holds for $X_1^N, X_2^N, X_3^N$. From the expression of $\sigma$ with respect to the basis $(e_j')_j$, we observe that $\mathcal{T}_\lambda^q$ is isomorphic, through the isomorphism given by

$$\begin{array}{rcl} X_1 & \longmapsto & X_1' \\ X_2 & \longmapsto & X_2' \\ X_3 & \longmapsto & q(X_2')^{-1}(X_1')^{-1} \otimes H' \end{array}$$

to the algebra $\mathcal{W}^q[X_1', X_2'] \otimes \mathbb{C}[(H')^{\pm 1}]$ (see [BL07, Lemma 17] for the definition of $\mathcal{W}^q[U, V]$).

Going back to the generic case of an ideal polygon, with $p \geq 4$, with simple calculations we can show that the following relations hold

$$m = p - 2$$
$$n = 2p - 3$$

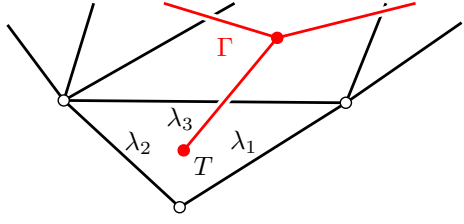

Figure 9: A leaf of $\Gamma$

where $n$ is the number of 1-cells of $\lambda$ and $m$ is the number of triangles in $\lambda$. Now, let $T = T_h$ be a triangle of $\lambda$ corresponding to a leaf in $\Gamma$. In order to simplify the notations, assume that the edges of $T$ are the 1-cells $\lambda_1$, $\lambda_2$, $\lambda_3$ of $\lambda$, ordered as in Figure 9. Because $\lambda_1$ and $\lambda_2$ belong to the only triangle $T$, the following holds

$$\begin{aligned} \sigma(e_1, e_j) &= 0 \\ \sigma(e_2, e_j) &= 0 \\ \sigma(e_1, e_2) &= 1 \end{aligned} \tag{23}$$

$$\sigma(e_1 + e_2 + e_3, e_j) = \sigma(e_3, e_j) \qquad \forall j \geq 4 \tag{24}$$

Now define a new basis $(e'_j)_j$ of $\mathbb{Z}^n \cong \mathcal{H}(\lambda; \mathbb{Z})$ given by

$$\begin{aligned} e'_1 &:= e_1 \\ e'_2 &:= e_2 \\ e'_3 &:= e_1 + e_2 + e_3 \\ e'_j &:= e_j \qquad\qquad\qquad \forall j \geq 4 \end{aligned}$$

The block diagonal matrix that represents $\sigma$ in this new basis is

$$\left( \begin{array}{cc|c} 0 & 1 & \\ -1 & 0 & \\ \hline & & \sigma' \end{array} \right)$$

thanks to relations 23. The equation 24 tells us that $\sigma'$ coincides with the bilinear form associated with the surface $S'$, obtained from $S$ by removing the triangle $T$, with the obvious induced triangulation $\lambda'$. Because $S'$ is an ideal polygon with $p-1$ punctures, we can reiterate this procedure until the $(m-1)$-th step, obtaining a surface $S^{(m-1)}$ composed of a single triangle. Analogously to what we have seen in Remark A.1, we construct an isomorphism between $(\mathcal{H}(\lambda; \mathbb{Z}), \sigma)$ and $(\mathbb{Z}^n, A)$, where $A$ is a bilinear skew-symmetric form, represented in the canonical basis of $\mathbb{Z}^n$ by the block diagonal matrix

$$\begin{pmatrix} 0 & 1 & & & & & \\ -1 & 0 & & & & & \\ & & \ddots & & & & \\ & & & 0 & 1 & \\ & & & -1 & 0 & \\ & & & & & 0 \end{pmatrix}$$

with $m = p - 2$ blocks $2 \times 2$ and a block $1 \times 1$ equal to zero. By inspection of the iterative procedure, we can see that the vector in $\mathcal{H}(\lambda; \mathbb{Z})$ corresponding to $e_n \in \mathbb{Z}^n$ is the element $(1, \ldots, 1) \in \mathcal{H}(\lambda; \mathbb{Z})$.

**Theorem A.2.** Let $S$ be an ideal polygon with $p \geq 3$ vertices and let $\lambda \in \Lambda(S)$ be an ideal triangulation. Then the Chekhov-Fock algebra $\mathcal{T}_\lambda^q$ of $S$ associated with $\lambda$ is isomorphic to

$$\mathcal{W}_{m,0,1}^q = \bigotimes_{i=1}^m \mathcal{W}^q \otimes \mathbb{C}[Z^{\pm 1}]$$

where $m = p - 2$. Moreover, the element $1 \otimes \cdots \otimes 1 \otimes Z \in \mathcal{W}_{m,0,1}^q$ corresponds to $H = \underline{X}^{(1,\ldots,1)} \in \mathcal{T}_\lambda^q$.

**Proposition A.3.** The following facts hold:

- if $q$ is not an $N$-th root of unity, then $\mathcal{Z}_\lambda^q$ is isomorphic to the direct sum of $\mathbb{C}^*$ and of the abelian subgroup generated by $H = \underline{X}^{(1,\ldots,1)} \in \mathcal{T}_\lambda^q$;

- if $q^2$ is a primitive $N$-th root of unity, then $\mathcal{Z}_\lambda^q$ is isomorphic to the direct sum of $\mathbb{C}^*$ and of the abelian subgroup generated by the elements $X_i^N$ with $i = 1, \ldots, n$, and $H$, endowed with the relation

$$H^N = q^{-N^2 \sum_{i<j} \sigma_{ij}} X_1^N \cdots X_n^N$$

*Proof.* Analogous to [BL07, Proposition 14] and [BL07, Proposition 15] using Theorem A.2 instead of [BL07, Theorem 12]. □

**Theorem A.4.** Let $S$ be an ideal polygon with $p \geq 3$ vertices and let $\lambda \in \Lambda(S)$ be an ideal triangulation of $S$. Fix $q \in \mathbb{C}^*$ such that $q^2$ is a primitive $N$-th root of unity. Then every irreducible representation $\rho \colon \mathcal{T}_\lambda^q \to \mathrm{End}(V)$ has dimension $N^{p-2}$ and it is uniquely determined, up to isomorphism, by the induced evaluation homomorphism $\rho \colon \mathcal{Z}_\lambda^q \to \mathbb{C}^*$ on the monomial center. Moreover, every homomorphism $\rho \colon \mathcal{Z}_\lambda^q \to \mathbb{C}^*$ is realized by a certain irreducible representation $\rho \colon \mathcal{T}_\lambda^q \to \mathrm{End}(V)$, unique up to isomorphism.

*Proof.* Analogous to the proof of [BL07, Theorem 20], using Theorem A.2, [BL07, Lemma 17], [BL07, Lemma 18] and [BL07, Lemma 19]. □

**Theorem A.5.** Let $S$ be an ideal polygon with $p \geq 3$ vertices and let $\lambda \in \Lambda(S)$ be an ideal triangulation of it, with $n$ 1-cells $\lambda_1, \ldots, \lambda_n$. Fix $q \in \mathbb{C}^*$ such that $q^2$ is a primitive $N$-th root of unity. Then $\mathscr{R}_{irr}(\mathcal{T}_\lambda^q)$ is in bijection with the set of elements $((x_i)_i; h) \in (\mathbb{C}^*)^n \times \mathbb{C}^*$, where $h$ is an $N$-th root of

$$q^{-N^2 \sum_{s<t} \sigma_{st}} x_1 \cdots x_n$$

and the correspondence associates, with the isomorphism class of an irreducible representation $\rho \colon \mathcal{T}_\lambda^q \to \mathrm{End}(V)$, the element $((x_i)_i; h)$ defined by the relations

$$\rho(X_i^N) = x_i \, id_V$$
$$\rho(H) = h \, id_V$$

for $i = 1, \ldots, n$.

*Proof.* Analogous to the proof of [BL07, Theorem 21], using Theorem A.4 and Proposition A.3. □

*Proof of Proposition 1.14.* It is sufficient to observe that the representation $\rho$ has dimension $N^m$, just as any irreducible representation, as seen in Theorem A.5. Therefore, if $\rho$ had a proper invariant subspace $0 \subsetneq W \subsetneq V_1 \otimes \cdots \otimes V_m$, then we would find an irreducible representation of $\mathcal{T}_\lambda^q$ with dimension strictly lower than $N^m$, which is absurd by virtue of Theorem A.5. □

# B  Proof of Theorem 1.22

In this Section we will give a proof of Theorem 1.22 using the Fusion property in order to simplify the case-by-case discussion of [Liu09], as suggested by the authors of [BBL07] in Subsection 4.1.

*Proof of Theorem 1.22.* In what follows, we will ignore the re-indexing property, because from the construction it is quite clear that this relation holds, but it would be very annoying to carry on all the indices in order to verify it.

Given $S$ a surface like above, we will firstly define the isomorphisms $\Phi_{\lambda\lambda'}^q$ in the case in which $\lambda$ and $\lambda'$ differ by an elementary move. We need to introduce some notations: we denote by $Q'$ the square in $S$ in which there is the diagonal exchange and we designate $T_1, T_2$ and $T_1', T_2'$ the triangles in $\lambda$ and $\lambda'$, respectively, that compose the square $Q'$. Furthermore, let $S_0$ be the surface obtained from $S$ by splitting it along all the edges of $\lambda$ (or $\lambda'$) except for the diagonal $\lambda_i$ of $Q'$ (or $\lambda_i'$). Then $S_0$ is the disjoint union of an embedded square $Q$ and triangles $T_i = T_i'$ for $i > 2$. $S_0$ is endowed with two triangulations $\lambda_0$, $\lambda_0'$, and $\lambda, \lambda'$ are obtained from them, respectively, by fusion on $S$. We clearly have that $\lambda_0$ and $\lambda_0'$ coincide on all the triangles $T_i = T_i'$ for $i > 2$, except on $Q$, where they coincide with the only two possible triangulations $\lambda_Q$ and $\lambda_Q'$, respectively, that $Q$ admits. By the Naturality, Disjoint Union and Diagonal Exchange properties, the isomorphism $\Phi_{\lambda_0\lambda_0'}$ is forced to be the extension to the quotient ring $\widehat{\mathcal{T}}_{\lambda_0}^q$ of the following injective map

$$\mathcal{T}_{\lambda_0}^q = \mathcal{T}_{\lambda_Q}^q \bigotimes_{j \neq 1,2} \mathcal{T}_{T_j'}^q \longrightarrow \widehat{\mathcal{T}}_{\lambda_Q'}^q \bigotimes_{j \neq 1,2} \widehat{\mathcal{T}}_{T_j'}^q \xrightarrow{\Phi_{\lambda_Q\lambda_Q'}^q \otimes_j id_{T_j}} \widehat{\mathcal{T}}_{\lambda_Q}^q \bigotimes_{j \neq 1,2} \widehat{\mathcal{T}}_{T_j}^q \longrightarrow \widehat{\mathcal{T}}_{\lambda_0}^q$$

It can be easily verified that $(\Phi_{\lambda_Q\lambda_Q'}^q)^{-1} = \Phi_{\lambda_Q'\lambda_Q}$ by explicit calculations on the formulae expressed in the Diagonal Exchange property. Hence we deduce that $\Phi_{\lambda_0\lambda_0'}^q$ is an isomorphism and $(\Phi_{\lambda_0\lambda_0'}^q)^{-1} = \Phi_{\lambda_0'\lambda_0}^q$. We would like to define $\Phi_{\lambda\lambda'}^q$ as $\hat{\imath}_{\lambda_0\lambda}^{-1} \circ \Phi_{\lambda_0\lambda_0'} \circ \hat{\imath}_{\lambda_0'\lambda'}$ but, in the first place, we must verify that the image of $\Phi_{\lambda_0\lambda_0'} \circ \hat{\imath}_{\lambda_0'\lambda'}$ is contained in $\hat{\imath}_{\lambda_0\lambda}(\widehat{\mathcal{T}}_\lambda^q)$. In order to prove this assertion, we have to discuss all the possible configurations of the square $Q'$ in $S$. Denote by $X_i$, $\overline{X}_i$, $X_i'$, $\overline{X}_i'$ the generators of the algebras $\mathcal{T}_\lambda^q$, $\mathcal{T}_{\lambda_0}^q$, $\mathcal{T}_{\lambda'}^q$, $\mathcal{T}_{\lambda_0'}^q$, respectively, and assume that the edges in $\lambda_Q \subseteq \lambda_0$ and $\lambda_Q' \subseteq \lambda_0'$ are indexed as their identifications in $\lambda$ and $\lambda'$, respectively. We refer to the cases of [Liu09] in the following discussion:

CASE 1 $\lambda_j, \lambda_k, \lambda_l, \lambda_m$ are all distinct.

Suppose that the edge $\lambda_l \in \lambda$ is the result of the identification of the edge $(\lambda_0)_l \in \lambda_Q$ and of an edge $(\lambda_0)_n$, belonging to a certain triangle different from $T_1$ and $T_2$ in $\lambda$. Then observe

$$(\Phi^q_{\lambda_0 \lambda'_0} \circ \hat{\iota}_{\lambda'_0 \lambda'})(X'_l) = \Phi_{\lambda_0 \lambda'_0}(\overline{X}'_l \overline{X}'_n) = (1 + q\overline{X}_i)\overline{X}_l \overline{X}_n$$
$$= \hat{\iota}_{\lambda_0 \lambda}((1 + qX_i)X_l)$$

Notice that the polynomial $(1 + qX_i)X_l$ coincides with the image of $X'_l$ in the case of an embedded square. This situation arise for every external edge of the square that is not identified to an other side of it. In the following, we will always omit the calculations for these cases and we will focus on the identified couples of sides of the square, if there is any. In conclusion, when the edges $\lambda_j$, $\lambda_k$, $\lambda_l$, $\lambda_m$ are all distinct, the expressions of the images of the elements $X'_s$ under $\Phi^q_{\lambda_0 \lambda'_0} \circ \hat{\iota}_{\lambda'_0 \lambda'}$ are the following

$$\Phi^q_{\lambda \lambda'}(X'_i) = \hat{\iota}_{\lambda_0 \lambda}(X_i^{-1})$$
$$\Phi^q_{\lambda \lambda'}(X'_j) = \hat{\iota}_{\lambda_0 \lambda}((1 + qX_i)X_j)$$
$$\Phi^q_{\lambda \lambda'}(X'_k) = \hat{\iota}_{\lambda_0 \lambda}((1 + qX_i^{-1})^{-1} X_k)$$
$$\Phi^q_{\lambda \lambda'}(X'_l) = \hat{\iota}_{\lambda_0 \lambda}((1 + qX_i)X_l)$$
$$\Phi^q_{\lambda \lambda'}(X'_m) = \hat{\iota}_{\lambda_0 \lambda}((1 + qX_i^{-1})^{-1} X_m)$$

CASE 2 $\lambda_j = \lambda_k$ and $\lambda_l \neq \lambda_m$.

Studying the case of $\lambda_j = \lambda_k$, we obtain

$$(\Phi^q_{\lambda_0 \lambda'_0} \circ \hat{\iota}_{\lambda'_0 \lambda'})(X'_j) = \Phi_{\lambda_0 \lambda'_0}(q^{-1} \overline{X}'_j \overline{X}'_k)$$
$$= q^{-1}(1 + q\overline{X}_i)\overline{X}_j(1 + q\overline{X}_i^{-1})^{-1} \overline{X}_k$$
$$= q^{-1}(1 + q\overline{X}_i)(1 + q^{-1}\overline{X}_i^{-1})^{-1} \overline{X}_j \overline{X}_k$$
$$= q^{-1} q \overline{X}_i (\overline{X}_j \overline{X}_k)$$
$$= \hat{\iota}_{\lambda_0 \lambda}(X_i X_j)$$

The image of the other elements have the same appearance of the ones in the Case 1.

CASE 3 $\lambda_j = \lambda_m$ and $\lambda_k \neq \lambda_l$.

In same spirit as in the previous case, we have

$$(\Phi^q_{\lambda_0 \lambda'_0} \circ \hat{\iota}_{\lambda'_0 \lambda'})(X'_j) = \hat{\iota}_{\lambda_0 \lambda}(X_i X_j)$$

The image of the other elements have the same appearance of the ones in the Case 1.

CASE 4 $\lambda_j = \lambda_l$ and $\lambda_k \neq \lambda_m$.

We observe

$$(\Phi^q_{\lambda_0 \lambda'_0} \circ \hat{\iota}_{\lambda'_0 \lambda'})(X'_j) = \hat{\iota}_{\lambda_0 \lambda}((1 + qX_i)(1 + q^3 X_i)X_j)$$

The image of the other elements have the same appearance of the ones in the Case 1.

CASE 5  $\lambda_k = \lambda_m$ and $\lambda_j \neq \lambda_l$.

We observe
$$(\Phi^q_{\lambda_0 \lambda'_0} \circ \hat\iota_{\lambda'_0 \lambda'})(X'_k) = \hat\iota_{\lambda_0 \lambda}((1+qX_i^{-1})^{-1}(1+q^3 X_i^{-1})^{-1} X_k)$$

The image of the other elements have the same appearance of the ones in the Case 1.

CASE 6  $\lambda_j = \lambda_k$ and $\lambda_m = \lambda_l$.

We observe
$$(\Phi^q_{\lambda_0 \lambda'_0} \circ \hat\iota_{\lambda'_0 \lambda'})(X'_j) = \hat\iota_{\lambda_0 \lambda}(X_i X_j)$$
$$(\Phi^q_{\lambda_0 \lambda'_0} \circ \hat\iota_{\lambda'_0 \lambda'})(X'_l) = \hat\iota_{\lambda_0 \lambda}(X_i X_l)$$

CASE 7  $\lambda_j = \lambda_m$ and $\lambda_k = \lambda_l$.

We observe
$$(\Phi^q_{\lambda_0 \lambda'_0} \circ \hat\iota_{\lambda'_0 \lambda'})(X'_j) = \hat\iota_{\lambda_0 \lambda}(X_i X_j)$$
$$(\Phi^q_{\lambda_0 \lambda'_0} \circ \hat\iota_{\lambda'_0 \lambda'})(X'_k) = \hat\iota_{\lambda_0 \lambda}(X_i X_k)$$

CASE 8  $\lambda_j = \lambda_l$ and $\lambda_k = \lambda_m$.

We observe
$$(\Phi^q_{\lambda_0 \lambda'_0} \circ \hat\iota_{\lambda'_0 \lambda'})(X'_j) = \hat\iota_{\lambda_0 \lambda}((1+qX_i)(1+q^3 X_i) X_j)$$
$$(\Phi^q_{\lambda_0 \lambda'_0} \circ \hat\iota_{\lambda'_0 \lambda'})(X'_j) = \hat\iota_{\lambda_0 \lambda}((1+qX_i^{-1})^{-1}(1+q^3 X_i^{-1})^{-1} X_k)$$

The above discussion allows us to define $\Phi^q_{\lambda \lambda'}$ as $(\hat\iota_{\lambda_0 \lambda})^{-1} \circ \Phi_{\lambda_0 \lambda'_0} \circ \hat\iota_{\lambda'_0 \lambda'}$, for every $S$ and for every $\lambda, \lambda' \in \Lambda(S)$ that differ by a diagonal exchange. By definition, we have

$$\Phi^q_{\lambda_0 \lambda'_0} \circ \hat\iota_{\lambda'_0 \lambda'} = \hat\iota_{\lambda_0 \lambda} \circ \Phi^q_{\lambda \lambda'} \tag{25}$$

Now, let $R$ be a surface with $\mu, \mu' \in \Lambda(R)$ triangulations that differ by a diagonal exchange along $\mu_i$. Suppose that $S$ is obtained from $R$ by fusion and that $\lambda, \lambda'$ are the triangulations of $S$ induced by $\mu, \mu'$ respectively. We want to prove that, in this situation, the following holds

$$\Phi^q_{\mu \mu'} \circ \hat\iota_{\mu' \lambda'} = \hat\iota_{\mu \lambda} \circ \Phi^q_{\lambda \lambda'} \tag{26}$$

Denoting by $\mu_0$ and $\mu'_0$ the triangulations of $R_0$ and $R'_0$ that appear in the definition of $\Phi^q_{\mu \mu'}$, as above, then we clearly have $R_0 = S_0$ and $\mu_0 = \lambda_0$, $\mu'_0 = \lambda'_0$, by construction. Because of the injectivity of $\hat\iota_{\mu_0 \mu}$, it is sufficient to prove

$$\hat\iota_{\mu_0 \mu} \circ \hat\iota_{\mu \lambda} \circ \Phi^q_{\lambda \lambda'} = \hat\iota_{\mu_0 \mu} \circ \Phi^q_{\mu \mu'} \circ \hat\iota_{\mu' \lambda'}$$

Now observe

$$\begin{aligned}
\hat\iota_{\mu_0 \mu} \circ \hat\iota_{\mu \lambda} \circ \Phi^q_{\lambda \lambda'} &= \hat\iota_{\mu_0 \lambda} \circ \Phi_{\lambda \lambda'} && \text{Relation 5}\\
&= \hat\iota_{\lambda_0 \lambda} \circ \Phi^q_{\lambda \lambda'} && \mu_0 = \lambda_0\\
&= \Phi^q_{\lambda_0 \lambda'_0} \circ \hat\iota_{\lambda'_0 \lambda'} && \text{Relation 25}\\
&= \Phi^q_{\mu_0 \mu'_0} \circ \hat\iota_{\mu'_0 \lambda'} && \mu_0 = \lambda_0 \text{ and } \mu'_0 = \lambda'_0\\
&= \Phi^q_{\mu_0 \mu'_0} \circ \hat\iota_{\mu'_0 \mu'} \circ \hat\iota_{\mu' \lambda'} && \text{Relation 5}\\
&= \hat\iota_{\mu_0 \mu} \circ \Phi^q_{\mu \mu'} \circ \hat\iota_{\mu' \lambda'} && \text{Relation 25}
\end{aligned}$$

and hence the relation 26, that is a "baby" version of the general Fusion property, holds whenever we are in the above situation.

Now we are going to define the $\Phi^q_{\lambda\lambda'}$ in the general case and to prove that the Composition relation holds. In order to do this, it is necessary to show that the isomorphisms, defined in the elementary cases, respect the Pentagon relation (the other relations between ideal triangulations in Theorem 1.21 are easier and can be verified in the same way).

Select in $S$ a triangulation $\lambda \in \Lambda(S)$ and $\lambda_i$, $\lambda_j$ two diagonals of a certain pentagon in $\lambda$. Designate also with $\lambda^{(0)}, \lambda^{(1)}, \ldots, \lambda^{(5)}$ the following sequence of triangulations

$$\lambda,\ \Delta_i(\lambda),\ \Delta_j\Delta_i(\lambda),\ \Delta_i\Delta_j\Delta_i(\lambda),\ \Delta_j\Delta_i\Delta_j\Delta_i(\lambda),\ \Delta_i\Delta_j\Delta_i\Delta_j\Delta_i(\lambda) = \tau_{ij}(\lambda)$$

Then we have to prove that

$$\Phi^q_{\lambda\lambda^{(5)}} \circ \Phi^q_{\lambda^{(5)}\lambda^{(4)}} \circ \cdots \circ \Phi^q_{\lambda^{(1)}\lambda} = \Phi^q_{\lambda\lambda}$$

Assuming for a moment that this relation holds for every $\lambda, \lambda' \in \Lambda(S)$, we can select a sequence $\lambda = \lambda^{(0)}, \ldots, \lambda^{(k)} = \lambda'$, by virtue of the Theorem 1.20, and define $\Phi^q_{\lambda\lambda'}$ as

$$\Phi^q_{\lambda\lambda^{(1)}} \circ \cdots \circ \Phi^q_{\lambda^{(k-1)}\lambda'}$$

Now, by virtue of the Theorem 1.21 and the assumed Pentagon relation (together with the others, on which we will not focus), it easy to verify that this is a good definition and that the Composition relation naturally holds.

Let $S$ be a surface and $\lambda \in \Lambda(S)$ a certain triangulation. Select also $\lambda_i$ and $\lambda_j$ diagonals of a pentagon in $S$. Let $R$ be the surface obtained by splitting $S$ along all the edges of $\lambda$ except for $\lambda_i$ and $\lambda_j$. Hence $R$ is the disjoint union of an embedded pentagon $P$ and some triangles. Let $\mu = \mu^{(0)}, \ldots, \mu^{(5)}$ the triangulations of $R$ such that their fusions induce the triangulations $\lambda = \lambda^{(0)}, \ldots, \lambda^{(5)}$ of $S$. Suppose that the following holds

$$\Phi^q_{\mu\mu^{(1)}} \circ \cdots \circ \Phi^q_{\mu^{(5)}\mu} = id_{\widetilde{\mathcal{T}}^q_\mu}$$

and observe

$$\hat{\iota}_{\mu\lambda} \circ \Phi_{\lambda\lambda^{(1)}} \circ \cdots \circ \Phi_{\lambda^{(1)}\lambda} = \Phi_{\mu\mu^{(5)}} \circ \hat{\iota}_{\mu^{(5)}\lambda^{(5)}} \circ \Phi_{\lambda^{(1)}\lambda^{(2)}} \circ \cdots \circ \Phi_{\lambda^{(5)}\lambda} \quad \text{Rel. 25}$$

$$\vdots$$

$$= \Phi^q_{\mu\mu^{(1)}} \circ \cdots \circ \Phi^q_{\mu^{(5)}\mu} \circ \hat{\iota}_{\mu\lambda} \quad \text{Rel. 25}$$

$$= \hat{\iota}_{\mu\lambda}$$

Then, because of the injectivity of $\hat{\iota}_{\mu\lambda}$, the assumption of $\Phi^q_{\mu\mu^{(1)}} \circ \cdots \circ \Phi^q_{\mu^{(5)}\mu} = id_{\widetilde{\mathcal{T}}^q_\mu}$ implies that the Pentagon relation holds in the general case. From the definition given of the $\Phi_{\mu^{(i+1)}\mu^{(i)}}$ it is clear that the identity $\Phi^q_{\mu\mu^{(1)}} \circ \cdots \circ \Phi^q_{\mu^{(5)}\mu} = id_{\widetilde{\mathcal{T}}^q_\mu}$ follows from the proof of the pentagon relation in case in which $S$ is an embedded pentagon. For the proof of this case we refer to [Liu09, Proposition 9].

Finally we have defined the isomorphisms $\Phi^q_{\lambda\lambda'}$ in the general case and we have proved that the Composition relation holds. It remains to verify that, with this definition, all the properties hold. The validity of the Naturality property is clear from the definition, just as the Disjoint Union property in case of $\lambda_i$ and $\lambda'_i$ that differs by a diagonal exchange. Now, by applying the Composition property, it is straightforward to prove the general case of the Disjoint Union property.

In the matter of the fusion property, we we have already shown it when $\mu$ and $\mu'$ differ by an elementary move in 26. In what follows, we will see how to deduce the general case from 26 and from the Composition relation. Suppose that $S$ is obtained by fusing a surface $R$ and that $\lambda, \lambda' \in \Lambda(S)$ are constructed as fusion of $\mu, \mu' \in \Lambda(R)$, respectively. Connect the triangulations $\mu$ and $\mu'$ with a sequence $\mu = \mu^{(0)}, \mu^{(1)}, \ldots, \mu^{(k)} = \mu'$, in which $\mu^{(l+1)}$ is obtained from $\mu^{(l)}$ by a diagonal exchange. Then we can define an induced sequence $\lambda = \lambda^{(0)}, \lambda^{(1)}, \ldots, \lambda^{(k)} = \lambda'$, where $\lambda^{(l)} \in \Lambda(S)$ is obtained by fusion of $\mu^{(l)}$. Now, using the Composition property, we see

$$\hat{\iota}_{\mu\lambda} \circ \Phi^q_{\lambda\lambda'} = \hat{\iota}_{\mu\lambda} \circ \Phi^q_{\lambda\lambda^{(1)}} \circ \cdots \circ \Phi^q_{\lambda^{(k-1)}\lambda'}$$
$$\Phi^q_{\mu\mu'} \circ \hat{\iota}_{\mu'\lambda'} = \Phi^q_{\mu\mu^{(1)}} \circ \cdots \circ \Phi^q_{\mu^{(k-1)}\mu'} \circ \hat{\iota}_{\mu'\lambda'}$$

Applying the relation 26 to $\lambda^{(i)}, \lambda^{(i+1)}, \mu^{(i)}$ and $\mu^{(i+1)}$ for every $i = 0, \ldots, k-1$ we observe

$$\begin{aligned}
\hat{\iota}_{\mu\lambda} \circ \Phi^q_{\lambda\lambda'} &= \hat{\iota}_{\mu\lambda} \circ \Phi^q_{\lambda\lambda^{(1)}} \circ \cdots \circ \Phi^q_{\lambda^{(k-1)}\lambda'} \\
&= \Phi^q_{\mu\mu^{(1)}} \circ \hat{\iota}_{\mu^{(1)}\lambda^{(1)}} \circ \Phi^q_{\lambda^{(1)}\lambda^{(2)}} \circ \cdots \circ \Phi^q_{\lambda^{(k-1)}\lambda'} \\
&\vdots \\
&= \Phi^q_{\mu\mu^{(1)}} \circ \cdots \circ \Phi^q_{\mu^{(k-1)}\mu'} \circ \hat{\iota}_{\mu'\lambda'} \\
&= \Phi^q_{\mu\mu'} \circ \hat{\iota}_{\mu'\lambda'}
\end{aligned}$$

as desired.

As said in the very beginning of this construction, the definition of the isomorphisms $\Phi^q_{\lambda_0\lambda'_0}$, when $S_0$ is a disjoint union of triangles and a square, is obliged by the Naturality, the Disjoint Union and the Diagonal Exchange properties. Furthermore, retracing the above discussion, we see that from the uniqueness of this base case, follows the uniqueness of the $\Phi^q_{\lambda\lambda'}$ in the general case and so we conclude the proof of the assertion. $\square$